\documentclass[final,10pt]{autart} 

\usepackage{graphicx}
\usepackage{lmodern}
\usepackage{cite}

\usepackage{amsfonts,amssymb,amsmath}
\usepackage{subcaption}
\usepackage[american]{circuitikz}

\allowdisplaybreaks 
\usepackage{graphics,epsfig,color}
\usepackage{wrapfig}
\usepackage[
  pdftex,
  bookmarks=false,
  hyperfootnotes=false
]{hyperref}

\usepackage{enumitem}
\usepackage{balance}
\usepackage{stackengine}

  \newif\ifitsdraft
\def\itsdraft{\global\itsdrafttrue}

\newenvironment{proof}{{\it Proof. }}{\hfill $\blacksquare$ } 

\newtheorem{exe}{Example}
\newtheorem{corol}{Corollary}
\newtheorem{ass}{Assumption}
\newtheorem{defin}{Definition}

\newenvironment{lemma}{\begin{lem}}{\hfill $\square$ \end{lem}}
\newenvironment{proposition}{\begin{prop}}{\hfill $\square$ \end{prop}}

\newenvironment{example}{\begin{exe}}{\hfill $\blacktriangle $ \end{exe}}
\newenvironment{remark}{\begin{rem} \rm}{\hfill $\blacktriangle $ \end{rem}}
\newenvironment{assumption}{\begin{ass}}{\hfill $\bullet$ \end{ass}}
\newenvironment{theorem}{\begin{thm}}{\hfill $\square$ \end{thm}}
\newenvironment{definition}{\begin{defin}}{\hfill $\bullet$ \end{defin}}

\definecolor{olivegreen}{rgb}{0.14,0.29,0}

\usepackage[utf8]{inputenc}
\usepackage{chngcntr}
\usepackage{apptools}
\newcommand{\comment}[1]{} 
\usepackage{tikz} 
\usepackage{mathtools}
\mathtoolsset{showonlyrefs} 
\usepackage[normalem]{ulem}
\usepackage{soul}
\usepackage{cancel}

\newlength{\myeqskip}  
\AtBeginDocument{%
\setlength\abovedisplayskip{\myeqskip}%
\setlength\belowdisplayskip{\myeqskip}%
\setlength\abovedisplayshortskip{\myeqskip-\baselineskip}%
\setlength\belowdisplayshortskip{\myeqskip}}

\itsdraft  

\begin{document}

\begin{frontmatter}
\title{Nagumo Theorem For Constrained Differential Inclusions: Transversal
Intersection and Set-Map Regularity
\thanksref{footnoteinfo}}

\thanks[footnoteinfo]{This paper was not presented at any IFAC meeting. This work has been supported by the LabEx EnergyAlps 2025. Corresponding author M. Maghenem.}

\vspace{-0.8cm}

\author[Grenoble]{Olayo Reynaud}\ead{olayo.reynaud@grenoble-inp.fr},
\author[Grenoble]{Mohamed Maghenem}\ead{mohamed.maghenem@gipsa-lab.grenoble-inp.fr},
\author[Benguerir]{Sadek Belamfedel Alaoui}\ead{Sadek.BELAMFEDEL@um6p.ma},
\author[Benguerir]{Adnane Saoud}\ead{Adnane.SAOUD@um6p.ma},
\author[Lyon]{Ahmad Hably}\ead{ahmad.hably@univ-lyon1.fr},

\address[Grenoble]{Universit\'e Grenoble Alpes, CNRS, Grenoble INP, GIPSA-lab, Grenoble, France}
\address[Benguerir]{College of computing, University Mohammed VI, Polytechnic (UM6P), Benguerir, Morocco.}
\address[Lyon]{Universit\'e Claude Bernard Lyon 1 (UCBL), LAGEPP laboratory, Lyon, France}

\maketitle

\vspace{-0.5cm}

\begin{abstract} 
We study a solution-independent condition certifying forward invariance of closed sets for constrained continuous-time systems. The proposed condition is only necessary in general, we establish its sufficiency under three assumptions. These assumptions, related to critical boundary points, restrict the dynamics and enforce mild transversality between the set constraining the dynamics and the set to render forward invariant. When the boundary of the latter set is within the interior of the former, we recover the well-known Nagumo characterization of forward invariance. We further simplify the proposed framework when the considered sets are regular; namely, by  considering \textit{finitely-represented} 
(resp., \textit{polytopic}) sets---defined as the intersection of sublevel sets of smooth (resp., linear) functions---yielding conditions directly checkable from the defining functions. These conditions are shown to reduce to finite \textit{linear matrix inequalities} (LMI)s under polytopic dynamics.   
These results carry important impacts on safety verification/falsification, reachability analysis, and compositional analysis and design for both constrained and hybrid systems.
\end{abstract}

\end{frontmatter}
 
\section{Introduction}
Forward invariance of sets for dynamical systems refers to the property whereby all solutions starting in a given set remain in it over their entire domain of definition. Set invariance plays a central role in control theory, due its 
particular importance in safety-critical applications, including transportation networks~\cite{wang2023multi}, power systems~\cite{alaraifiDomainStabilityCharacterization2024}, and robotics~\cite{tanner2007flocking, glotfelter2019hybrid}.  Furthermore, forward invariance provides a cornerstone for establishing more elaborate dynamical behaviors. In particular, 
invariance analysis is instrumental for establishing asymptotic stability using weak Lyapunov inequalities, as in the well-known invariance principle~\cite{sanfeliceInvariancePrinciplesHybrid2007} and Matrosov’s theorem~\cite{sanfeliceAsymptoticStabilityHybrid2009}. 
Beyond safety and stability, extensions of invariance theory have been developed to encompass concepts of quasi-invariance, conditional quasi-invariance~\cite{laddeAnalysisInvariantSets1972}, contractivity~
\cite{blanchini1999survey}, as well as signal temporal logic \cite{9894077}.

Characterizing forward invariance typically requires conditions that guarantee the property without explicitly computing the system's trajectories. A relevant context where this task  is not fully covered yet concerns dynamical systems modeled by differential inclusions subject to state constraints. That is, the solutions are well defined in this case as long as they are contained within the set where the differential inclusion is defined. 
The main motivation behind studying this class of systems is that they represent a key component of hybrid systems
\cite{goebel2012hybrid}. Constrained systems are also common when modeling physical power systems~\cite{chiangStabilityRegionsConstrained2015}, and they appear naturally in the contexts of singular systems~\cite{graciaGeneralizedGeometricFramework1992} and sweeping processes~\cite{brogliato2020dynamical}.

\subsection{Background}

The earliest solution-independent characterization of forward invariance is due to Nagumo~\cite{nagumo1942lage}, who established a necessary-and-sufficient condition for the existence of a solution that remains in a closed set $K$ from every initial condition in $K$. This property, sometimes referred to as \textit{weak forward invariance}~
\cite{clarke1998nonsmooth}, coincides with forward invariance when the solutions are unique. Specifically, Nagumo's theorem states that a closed set $K \subset \mathbb{R}^n$ is weakly forward invariant for the system $\dot{x} = f(x)$, with $f$ continuous, if and only if $f(x) \in T_{K}(x)$ for all $x \in K$, where $T_{K}(x)$ is the contingent cone to $K$ at $x$; see Section~\ref{sec_preliminaries}. The same condition is shown to ensure forward invariance when $f$ satisfies the one-sided
Lipschitz-continuity property \cite{redheffer1972theorems}. Subsequent works generalized Nagumo's condition to stronger invariance notions and extended the result to differential inclusions \cite{Aubin1991Viability} and impulse differential inclusions \cite{aubin2002impulse}. For instance, it is  shown 
in~\cite[Theorem 5.3.4]{Aubin1991Viability} that, for a 
locally-Lipschitz 
set-valued map $F: \mathbb{R}^n \rightrightarrows \mathbb{R}^n$ with convex and compact images, the closed set $K$ is forward invariant for the differential inclusion
$\dot{x} \in F(x)$ 
if and only if
\begin{equation} \label{eqTCone}
  F(x) \subset T_{K}(x) \qquad  \forall x \in K.
\end{equation}

More recently, the development of hybrid-systems theory~\cite{goebel2012hybrid} motivated analogous studies. Such a general modeling framework captures the system's continuous-time behavior through a constrained differential inclusion of the form 
\begin{equation} \label{eq_constrained_DI_intro}
    \Sigma: \quad \dot{x} \in F(x) \qquad x \in C \subset \mathbb{R}^n.
\end{equation}
In this context, the presence of the  constraint $x \in C$ renders condition~\eqref{eqTCone} overly conservative. For instance, when $x$ is within $\partial C \cap \partial K$, vectors within $F(x)$ `pointing' outside $C$ are, roughly speaking, not able to generate solutions. Hence, they should not be forced to be within $T_K(x)$. Less restrictive (compared to \eqref{eqTCone}) sufficient conditions for forward invariance in hybrid systems are proposed in~\cite{chaiForwardInvarianceSets2019, chaiForwardInvarianceSets2021} and \cite{maghenemSufficientConditionsForward2021}, circumventing the subtlety at the intersection of the boundaries of $K$ and $C$ by assuming unnecessary conditions on neighborhoods of the latter points. Similar results for polyhedral sets~\cite{t.doreaBInvarianceConditionsPolyhedral1999}, interconnected systems~\cite{saoudAssumeguaranteeContractsContinuoustime2021}, and subsystem composition \cite{sanfeliceInvariantsOpenHybrid2023} are also proposed. To the best of our knowledge, these conditions did not nuance the geometry of the intersection between the two boundaries so to derive necessary-and-sufficient conditions. 

In many applications, the sets $C$ and $K$ correspond to the intersection of sublevel sets of differentiable scalar functions (i.e., sets defined by $g_i(x) \leq 0$ for all $i \in \{1,2,... \}$). We refer to these sets as \textit{finitely-represented} sets,
or \textit{polytopic}  sets when the defining functions are linear. Forward invariance of finitely-represented
sets is widely studied, through their defining functions, sometimes called \textit{barrier functions}~
\cite{amesControlBarrierFunction2017,glotfelterNonsmoothBarrierFunctions2017, prajnaFrameworkWorstCaseStochastic2007, daiBarrierCertificatesRevisited2017, kongExponentialConditionBasedBarrierCertificate2013}. These tools are especially powerful in the context of input-affine systems since the invariance  conditions are enforced via solving quadratic programs~\cite{amesControlBarrierFunction2017, wielandConstructiveSafetyUsing2007}. 

That said, the challenging question of providing necessary-and-sufficient conditions in the context of constrained systems remains open even when the sets $K$ and $C$ are finitely represented or polytopic. 

Answering the aforementioned question offers immediate impact on methodologies that have not been applied to constrained and hybrid systems, to the best of our knowledge. This includes: 1- exact verification and falsification of a given safe set: sufficient conditions, such as barrier certificates \cite{prajnaOptimizationBasedMethodsNonlinear2005}, can only prove safety, whereas a necessary-and-sufficient condition can also disprove it, as illustrated in \cite{jouretSafetyVerificationNeuralNetwork2023}, which relies on the equivalence provided by Nagumo's theorem to develop a complete 
safety-verification algorithm, in the unconstrained case, for linear systems controlled by neural networks; 2- compositional characterization of forward invariance for large-scale interconnected systems using assume-guarantee contracts, where global invariance is deduced from local subsystem certificates through both implications of Nagumo's theorem \cite[Theorem~5.3.4]{Aubin1991Viability}, as in \cite[Theorem~3]{saoudAssumeguaranteeContractsContinuoustime2021}, extending this deduction to constrained subsystems requires the constrained counterpart of Nagumo's theorem, which is precisely what the present paper establishes.
3- characterizing the 
non-increase of a state-dependent function along  solutions, through analyzing forward invariance of the restricted epigraph $\text{epi}B \cap (C \times \mathbb{R})$, $\text{epi} B := 
 \{(x, r) \in \mathbb{R}^n
\times \mathbb{R} : r \geq B(x) \}$,  for the augmented constrained system \cite{clarke1998nonsmooth, maghenem2022necessary}
    \begin{center}
    $
    (\dot{x},\dot{r})
    \in 
    (F(x),0)
    \qquad 
    (x,r) \in C \times \mathbb{R}. 
    $ \end{center}

\subsection{Contributions}

In this paper, we establish a framework, in which, the relaxed Nagumo-inspired criterion
\begin{equation} \label{eqconstNag}
    F(x) \subset T_{K}(x) \qquad \forall x \in K \cap \mathrm{int}(C)
\end{equation}
is equivalent to forward invariance of $K$ for $\Sigma$. When the boundary of $K$ lies strictly within $\mathrm{int}(C)$, condition~\eqref{eqconstNag} reduces to~\eqref{eqTCone}. 
However, as we show via examples, condition~\eqref{eqconstNag} is only necessary in general. 
We, therefore, introduce three assumptions to conclude sufficiency of~\eqref{eqconstNag}. A particular focus is made on a specific subset, denoted $\mathcal{P} \subset \partial K \cap \partial C$, from where, trajectories can leave $K$ ``immediately"  even when~\eqref{eqconstNag} holds. In particular, 
\begin{itemize}
\item[1)] we \textit{propagate} the inclusion in~\eqref{eqconstNag} by imposing suitable structural assumptions so that it holds on $\mathcal{P}$. 
\item[2)] we impose a form of  \textit{transversal}, or \textit{non-tangent}, intersection between the sets $K$ and $C$. This assumption applies to general sets and is weaker than those commonly used in the literature~\cite{Aubin1991Viability, clarkeOptimizationNonsmoothAnalysis1990}. 
\item[3)] we restrict $F$, on $\mathcal{P}$, so its vectors are not tangent to $\partial K \cap \partial C$. Thereby, we are able to show that the first-order approximation of every solution, starting from $\mathcal{P}$, is \textit{decisive} in terms of indicating a specific side of $\partial K$ the solution must be leaning towards while leaving $K$. 
\end{itemize}

Then, the Nagumo-type condition~\eqref{eqconstNag} is used to exclude this unique possibility, and thus, the possibility of leaving $K$ immediately. Examples are proposed to illustrate  the assumptions and the loss of forward invariance when they are not verified.

 In the second part of the paper, 
 we restrict ourselves to finitely-represented sets $K$ and $C$. We show that the transversality condition can be simplified to a more tractable, although more restrictive, condition involving the gradients of the defining functions. The latter condition is not needed when the sets are polytopic, since the general transversality condition is shown to hold for free in the polytopic case. Furthermore, the propagation assumption holds under finitely-represented sets as a direct consequence of the new transversality assumption and the Nagumo condition \eqref{eqconstNag}. The latter Nagumo condition is enough to verify the propagation assumption under polytopic sets, thanks to convexity.
 
 Finally,  when the sets $K$ and $C$, and $F$ are polytopic, we derive a characterization of forward invariance in the form of finite linear matrix inequalities (LMI)s.

The remainder of this paper is organized as follows. Preliminaries related to differential inclusions and tangent cones are in Section~\ref{sec_preliminaries}.  The problem formulation, its motivation, and the rationale of the proposed assumptions are in Section~\ref{sec_problem_formulation}.
Our main results for general sets, finitely-represented, and polytopic sets are stated and proved in Sections~\ref{sec_main_result}, \ref{sec_prac}, and \ref{sec:app-lmi}, respectively. 

For clarity, some of the intermediate proofs and technical lemmas are relegated to the Appendix. Furthermore, some well-known notions, intermediate results, as well as some illustrative examples are omitted to the online version \cite{arXiv_version_auto}.

A preliminary version 
of this work is in~\cite{reynaudNagumoTypeCharacterizationForward2025}, where intermediate results and detailed proofs were not included. More importantly, the contexts of 
finitely-represented sets and finite LMI characterizations under polytopic sets were not considered.

\textbf{Notation.}
We denote by $\mathbb{N}$ and   $\mathbb{N}_{\geq 0}$ the sets of positive and nonnegative integers, respectively. 
For $n \in \mathbb{N}$, we let $[n] :=  \{1,2,\ldots, n\}$. For $K \subset \mathbb{R}^n$, $\text{cl}(K)$, $\text{int}(K)$, and $\partial K$ denote, respectively, the closure, the interior, and the boundary of $K$. For $x \in \mathbb{R}^n$, $\mathrm{Proj}_{K}(x) := \{ y \in K : |x-y| = \inf_{z \in K} |x - z| \}$ denotes the projection of  $x$ on $K$, where $|\cdot|$ is the Euclidean norm. 
The closed ball of radius $\epsilon$ centered at $x$ is denoted by $\mathbb{B}_\epsilon(x)$, and $\mathbb{B} := \mathbb{B}_1(0)$ represents the unit ball centered at the origin. 
For a map $F : \mathrm{dom}\,F \to \mathrm{Im}(F)$ (or $F : \mathrm{dom}\,F \rightrightarrows \mathrm{Im}(F)$ in the set-valued case),  $\mathrm{dom}\,F$ and $\mathrm{Im}(F)$ denote, respectively, its domain of definition and its image set. Given $A \subset \mathrm{dom}\,F$, $F(A) := \{F(x): x \in A\}$ denotes the image of $A$ under $F$. When $F$ is a single-valued map, $x \in \mathrm{dom}\,F$ means that $F(x) \in \mathrm{Im}(F)$ is well defined and finite and, when $f$ is a 
set-valued map, $x \in \mathrm{dom}\,F$ means that $F(x) \subset \mathrm{Im}(F)$ is a subset of $\mathrm{Im}(F)$, hence, $F(x)$ can possibly be empty or unbounded. We use $o: \mathbb{R} \to \mathbb{R}$ to denote any function satisfying $\lim_{h\to 0^+} |o(h)|/h=0$.  A set-valued map $F$, with $\mathrm{dom}\,F \subset \mathrm{Im}(F)$, is said to be \textit{one-sided locally Lipschitz} if, for each compact set $\mathcal{N} \subset \operatorname{dom} F$, there exists $L > 0$ such that, for each $x_1, x_2 \in \mathcal{N}$ and for each $f_1 \in F\left(x_1\right)$, there exists $f_2 \in F\left(x_2\right)$ such that 
$(f_1-f_2)^\top(x_1-x_2) \leq L|x_1-x_2|^2$. Finally, every vector inequality is interpreted 
componentwise. Throughout the paper, the symbol $\otimes$ denotes the Kronecker product of matrices. The notation $d_H(A,B)$ denotes the Hausdorff distance between $A$ and $B$.

We refer to \cite{Aubin1991Viability} for the definitions of upper and lower semicontinuity, and continuity  for set-valued maps.
\section{Preliminaries}
\label{sec_preliminaries}

Consider the constrained differential inclusion  $\Sigma$ in \eqref{eq_constrained_DI_intro}. The set $C$ and the map $F : C \rightrightarrows \mathbb{R}^n$ are assumed to
verify the following standing assumption.
\begin{itemize} 
    \item [(SA)] The set $C$ is closed, and $F$ is upper semicontinuous with  non-empty, compact, and convex images.  
\end{itemize}

We recall that $\phi : \operatorname{dom}\phi \rightarrow \operatorname{Im}(\phi)$ is a solution to $\Sigma$ if $\operatorname{Im}(\phi) \subset C$ and 
\begin{itemize}
    \item  $\operatorname{dom}\phi = [0,T]$ for some $T \geq 0$ or $\operatorname{dom}\phi = [0,T)$ for some $T \in \mathbb{R}_{\geq 0} \cup \{+\infty\}$;
    \item  $\phi$ is locally absolutely continuous; \item $\dot{\phi}(t) \in F(\phi(t))$ for almost all $t \in \operatorname{dom}\phi$.
\end{itemize}

Furthermore, given a closed set $K \subset \mathbb{R}^n$, a solution $\phi$, starting from $x \in K$, is said to leave $K$ \textit{immediately} if
\begin{equation}\label{eq_leave_inmediatelly}
    \phi((0,T]) \subset C \backslash K \quad \text{ for some } T > 0.
\end{equation}

Moreover,  the set $K \subset \mathbb{R}^n$ is forward invariant for $\Sigma$ if, for every solution $\phi$ to $\Sigma$, with $\phi(0) \in K$, $\phi(\text{dom}\, \phi) \subset K$.

 \ifitsdraft

We recall below  the well-established notions of \textit{contingent}, \textit{Clarke} tangent, \textit{adjacent} contingent, 
and \textit{Dubovitskiy} cones, some of their key properties and relationships, as well as an example illustrating their differences.

- The \textit{contingent} cone to the set $K \subset \mathbb{R}^n$ at $x \in \operatorname{cl}(K)$, denoted by $T_K(x)$, is defined by~\cite[Definition 5.1.1]{Aubin1991Viability}
\[
    T_K(x) := \left\{v \in \mathbb{R}^n : \exists t_i \to 0^+,~ \exists v_i \to v,~ x + t_i v_i \in K\right\}.
\]
Note that $T_K(x)$ is always closed. Furthermore, it holds that $T_{K}(x) = \mathbb{R}^n$ whenever $x \in \mathrm{int}(K)$. Furthermore, if the set-valued map $K \ni y \leadsto T_K(y)$ is lower semicontinuous at $x \in \text{cl}(K)$, then $K$ is said to be \textit{sleek} at $x$. We refer the reader to~\cite[Definition 1.4.2]{aubinSetValuedAnalysis2009} for the definition of lower semicontinuity for set-valued maps.
\\
- The \textit{Clarke} tangent cone to the set $K \subset \mathbb{R}^n$ at $x \in \operatorname{cl}(K)$, denoted by $C_K(x)$, is defined by~\cite[Proposition 3.7]{clarkeGeneralizedGradientsApplications1975}
\[
    \begin{split}
        C_K(x) := \{&v \in \mathbb{R}^n : \forall t_i \to 0^+,~\forall x_i \rightarrow_K x,\\ 
        &  ~ \exists v_i \rightarrow v, ~ x_i + t_i v_i \in K \},
    \end{split}
\]
where $x_i \rightarrow_K x$ means that the sequence $\{x_i\}_i \subset K$ satisfies $x_i \to x$. Note that the map $x \leadsto C_K(x)$ is lower semicontinuous, and that $C_K(x)$ is a closed and convex subset of $T_K(x)$. $K$ is said to be \textit{tangentially regular} at $x \in \operatorname{cl}(K)$ if $C_K(x) = T_K(x)$~\cite[Definition 4.1.5]{aubinSetValuedAnalysis2009}. It is proven in~\cite[Proposition 4.2.]{bigolinGeometricCharacterizationsC1C12012} that if $K$ is a $C^1$ manifold of $\mathbb{R}^n$ then it is tangentially regular at every $x \in K$. Furthermore, if $K$ is sleek at $x \in \operatorname{cl}(K)$, then it is tangentially regular at $x$~\cite[Theorem 4.1.8]{aubinSetValuedAnalysis2009}.\\
- The \textit{adjacent} contingent cone to the set $K \subset \mathbb{R}^n$ at $x \in \operatorname{cl}(K)$, denoted by $T^a_K(x)$, is defined by~\cite[Definition 4.1.5]{aubinSetValuedAnalysis2009}
\[
    T^a_K(x) := \left\{v \in \mathbb{R}^n : \forall t_i \to 0^+,~\exists v_i \to v,~ x + t_i v_i \in K \right\}.
\]
The adjacent contingent cone is also referred to in the literature as the derivable cone~\cite{rockafellarVariationalAnalysis2009}. Note that $T^a_K(x)$ is always closed and satisfies
\begin{equation}\label{eq_implications_Tcones}
    C_K(x) \subset T^a_K(x) \subset T_K(x) \qquad \forall x \in \operatorname{cl}(K).
\end{equation}
Moreover, the set $K$ is said to be \textit{derivable} at $x \in \operatorname{cl}(K)$ if $T_K(x) = T^a_K(x)$. From~\eqref{eq_implications_Tcones}, it follows that tangential regularity of $K$ at $x \in \operatorname{cl}(K)$ implies derivability of $K$ at $x \in \operatorname{cl}(K)$. \\
- The \textit{Dubovitskiy} cone to the set $K \subset \mathbb{R}^n$ at $x \in \text{cl}(K)$, denoted by $D_K(x)$, is given by~\cite[Definition 4.3.1]{aubinSetValuedAnalysis2009}
\[
    \begin{split}
        \hspace{-0.3cm} & D_K(x) := 
        \\ &
        \hspace{-0.3cm} \left\{  v \in \mathbb{R}^n :  \forall t_i \to 0^+,~ \forall v_i \rightarrow v,~  x + t_i v_i \in \text{int}(K) ~ \forall i~\text{large} \right\}.
    \end{split}
\]
It is also true that $v \in D_K(x)$ if and only if there exist $\alpha, \epsilon > 0$ such that
\[
    x + a (v + e) \subset \text{int}(K) \qquad \forall a \in (0,\alpha], \quad \forall e \in \epsilon \mathbb{B}
\]

Although  the aforementioned cones are sometimes introduced differently, the definitions we recall here are particularly useful for our proofs. 

The next example illustrates the latter constructs. 
\begin{example}
    Let 
    \begin{align*}
    K_1 & := \left\{x \in \mathbb{R}^2 : x_2 \leq \sqrt{|x_1|}\right\},
    \\
    K_2 & := \left\{x \in \mathbb{R}^2 :x_1 \neq 0, x_2 = x_1 \sin\left(\frac{1}{x_1}\right)\right\}.
    \end{align*}
    It holds that
    \[
        C_{K_1}(0) = \{x \in \mathbb{R}^2 : x_2 \leq 0\} \subsetneq T_{K_1}^a(0) = T_{K_1}(0) = \mathbb{R}^n,
    \]
    and
    \[
        \begin{split}
            C_{K_2}(0) = T_{K_2}^a(0) = \{0\} \subsetneq &T_{K_2}(0) =\\
            &\{x \in \mathbb{R}^2 : |x_2| \leq |x_1|\}.
        \end{split}
    \]
    In fact, the set $K_2$ is proposed in~\cite[Section 6.A]{rockafellarVariationalAnalysis2009} as an example of set that is not derivable. The Dubovitskiy cone at the origin is
    \[
        D_{K_1}(0) = \mathbb{R}^2 \backslash \{x \in \mathbb{R}^2 : x_1 >0,  x_2 = 0\}, ~ D_{K_2}(0) = \emptyset.
    \]   
\end{example}

\else

Finally, we refer the reader to the extended version \cite{arXiv_version_auto}, where we recall the well-established notions of
\textit{contingent}, \textit{Clarke} tangent, \textit{adjacent} contingent, 
and \textit{Dubovitskiy} cones, denoted by
$T_K$ $C_K$, $T^a_K$, and $D_K$, respectively, some of their key properties and relationships, as well as an example illustrating their differences.

\fi

\section{Problem Formulation} \label{sec_problem_formulation}
  
For general constrained systems $\Sigma$, with $F$ lower semicontinuous,  condition~\eqref{eqconstNag} is shown to be only necessary for forward invariance of $K$. It becomes sufficient in the following restrictive  situations.\\
- When $K \subset 
\text{int}(C)$, so that condition~\eqref{eqconstNag} reduces to 
\begin{equation}\label{equiv2}  F(x) \subset T_{K}(x) \qquad \forall x \in  K.
  \end{equation}
Sufficiency of~\eqref{equiv2} is a consequence of Nagumo's theorem for unconstrained systems~\cite[Theorem 5.3.4]{Aubin1991Viability}. 
 \\ 
- When $\partial K \cap \text{int}(C)$ is dense in $\partial (K \cap C)$ (i.e., every $x \in  \partial (K \cap C)$ is the limit of a sequence $\{x_i\}^{\infty}_{i = 1} \subset \partial K \cap \text{int}(C)$), $F$ is lower semicontinuous, and $T_{K \cap C}$ is upper semicontinuous on the same set. In this case, we use the density property to conclude that 
$\partial (K \cap C) = \left(\partial K \cap \text{int}(C) \right) \cup \left(\partial K \cap \partial C \right).  $
Furthermore,~\eqref{eqconstNag} entails
\begin{equation*}
    F(x) \subset T_{K \cap C}(x) \left( = T_{K}(x) \right) \quad \forall x \in \partial K \cap \text{int}(C).
\end{equation*} 
Hence, using lower semicontinuity of $F$ and upper semicontinuity of $T_{K \cap C}$ on $\partial K \cap \partial C$ imply that
\begin{equation}\label{eqtoshow} 
    F(x) \subset T_{K \cap C}(x) \qquad  \forall x \in \partial K \cap \partial C.  
\end{equation}
And thus
$F(x) \subset T_{K \cap C}(x)$ for all $x \in \partial (K \cap C)$. The latter guarantees forward invariance of $K \cap C$ by the aforementioned Nagumo's invariance theorem. 

For general sets $K$ and $C$, however, sufficiency of~\eqref{eqconstNag} requires additional assumptions, as we show later.   

\subsection{Rationale} \label{subsec_rationale}

Our approach is based on ensuring that every solution $\phi$ that leaves $K$ immediately satisfies:
\begin{itemize}
  \item [($\star$)] There exists $T > 0$ such that 
  \[
    \text{Proj}_{\partial K}(\phi(t)) \cap \text{int}(C) \neq \emptyset \qquad \forall t \in (0,T] \subset \text{dom}\, \phi.
  \]  
\end{itemize}
We then show that \eqref{eqconstNag} excludes the possibility of ($\star$), and thus the possibility of solutions leaving $K$ immediately.

We also show that the only solutions that can leave $K$ immediately and not verify ($\star$) must start from the set
\begin{equation*}
\begin{aligned}
    & \mathcal{P} := \left\{ x \in \text{cl}(C \backslash K) \cap \partial K \cap \partial C :   \right. \\ & \left.
   \hspace{1cm}
\overline{ \left[ F(x) \subset \mathbb{R}^n \backslash T_C(x) \right]} \land \overline{ \left[ F(x) \subset  D_K(x) \right] } \land
    \right.
    \\ &
   \left.
\hspace{-0.2cm} \overline{ \left[ \exists \epsilon > 0,   S \subset K : x \in S, ~ \widehat{F}(y) \subset T_{S}(y) ~~ \forall y \in \mathbb{B}_\epsilon(x) \cap S \right]}
   \land \right.
    \\ &
   \left.
    \overline{  
   \begin{bmatrix} \left( \text{$C$  convex at $x$} \right) \land 
   \\
   \left( \exists \epsilon> 0 : F(y) \subset T_{C \cap K}(y) ~~ \forall y \in  \mathbb{B}_\epsilon(x) \cap \partial K \cap \partial C \right)
\end{bmatrix} }
   \right\},
\end{aligned}   
\end{equation*}
where the set $S$ is closed and $\widehat{F}$ is any extension of $F$ to $S$ (when $S \not\subset C$) that preserves the regularities of $F$. \footnote{We say that $C$ is convex at $x$ when there exists $\delta > 0$ such that $\mathbb{B}_{\delta}(x) \cap C$ is convex. Moreover, $\overline{[P]}$ stands for the negation of the proposition $[P]$.}


When one of the negated properties in the definition of $\mathcal{P}$ holds at some $x \in \text{cl}(C \backslash K) \cap \partial K \cap \partial C$, we show, respectively, that  the solutions starting from $x$ are trivial, enter $\text{int}(K)$ immediately, cannot leave $S$ immediately, and cannot leave $K$ immediately if \eqref{eqconstNag} additionally  holds. Furthermore, when $x \notin \text{cl}(C \backslash K) \cap \partial K \cap \partial C$, we show that continuity of solutions plus \eqref{eqconstNag} are enough to exclude the possibility of leaving $K$ immediately.
We prove these facts using arguments that do not rely on the framework we establish here.


As a result, to guarantee ($\star$) for every solution $\phi$ starting from $\mathcal{P}$ that leaves $K$ immediately, we enforce 
\begin{equation} \label{eqMotiv1}
    \text{cl}(F_\phi(x)) \subset T_{\partial K \cap C}(x) \backslash T_{\partial K \backslash \text{int}(C)}(x) \quad \forall x \in \mathcal{P},
\end{equation}
where 
\begin{equation}
\label{eq_def_Fphi}
    \begin{aligned}
    & F_{\phi}(x)  :=  \left\{  v \in \mathbb{R}^n : \exists 
    \{ t_i \}^\infty_{i=1} \subset \mathbb{R}_{\geq 0} ~ \text{with} \right. 
    \\ &  
    \left. \lim_{i \rightarrow \infty} t_i = 0 \text{ and } \lim_{i \rightarrow \infty} \frac{\phi(t_i) - x}{t_i} = v  \right\} \in T_C(x)
    \end{aligned}
\end{equation}
represents the set of initial speeds of the solution $\phi$, and the latter inclusion holds by definition of $T_C$.
In words, $F_\phi(x)$ directs the solutions along the part of $\partial K$ that is in $\text{int}(C)$, so that, locally, the solution's projection on $\partial K$ lands in $\text{int}(C)$.

The proposed assumptions will restrict the sets $K$ and $C$ as well as $F$ at $\mathcal{P}$ so that \eqref{eqMotiv1} holds.

\ifitsdraft
   In the following example, we illustrate the set $\mathcal{P}$, and properties ($\star$) and~\eqref{eqMotiv1}.
    \begin{example}
        Consider the sets 
        \begin{align*}
            C := \{x \in \mathbb{R}^2 : x_2 \geq 0\}, ~~
            K  := \{x \in \mathbb{R}^2 : x_1 \geq 0\},
        \end{align*}
         and let $F(x) := \{[-x_2 ~~ 1]^\top\}$.  
         The solution $\phi(t) =
          \left[- t^2/2 ~~ t \right]^\top$ for all $t \geq 0$, starting from the point $x \in \{ [0 ~~ 0]^\top \} = \partial K \cap \partial C$, leaves $K$ immediately; hence, we necessary have $\mathcal{P} = \{ [0 ~~ 0]^\top \}$.

         Since $\phi$ is differentiable at $t=0$, the set of initial speeds defined in~\eqref{eq_def_Fphi} is $F_\phi(x)= \{\dot{\phi}(0)\} = \{[0 ~~ 1]^\top\}.$  Moreover,
        \[
            T_{\partial K \cap C}(x) = \{v \in \mathbb{R}^2 : v_1=0,\; v_2\geq 0\},
        \]
        and
        \[
            T_{\partial K \setminus \mathrm{int}(C)}(x) = \{v \in \mathbb{R}^2 : v_1=0,\; v_2\leq 0\}.
        \]
        Hence,
        \[
            F_\phi(x) \subset T_{\partial K \cap C}(x) \setminus T_{\partial K \setminus \mathrm{int}(C)}(x),
        \]
        which shows that~\eqref{eqMotiv1} holds. Finally, for every $t>0$,
        \[
            \phi(t)\in \{x\in\mathbb{R}^2:x_2>0,\;x_1<0\}.
        \]
        Consequently,
        \[
            \mathrm{Proj}_{\partial K}(\phi(t)) = \{[0 ~~ t]^\top\} \subset \partial K\cap\mathrm{int}(C) \qquad \forall t>0,
        \]
        and therefore property~($\star$) also holds.
    \end{example}
\else

An illustration of $\mathcal{P}$, ($\star$), and~\eqref{eqMotiv1} is in~\cite[Example 2]{arXiv_version_auto}.

 \fi

\subsection{Assumptions} \label{subsec_assumptions}

To enforce \eqref{eqMotiv1},
we need to guarantee that 
\begin{equation}\label{eqMotiv}
    F_{\phi}(x) \subset T_{\partial K \cap C}(x) \qquad \forall x \in \mathcal{P}.
\end{equation}
Hence, we consider the following assumption that propagates the restriction in \eqref{eqconstNag} to $\mathcal{P}$.
\begin{assumption}[Propagation] \label{ass_propagation}
  $F \cap T_C \subset T_{K}$ on $\mathcal{P}$.
\end{assumption}
 Under Assumption~\ref{ass_propagation}, we show that any solution $\phi$ starting from $x \nolinebreak\in\nolinebreak \mathcal{P}$ and leaving $K$ immediately must satisfy $F_\phi(x)\nolinebreak\subset\nolinebreak T_{\partial K}(x)$. Furthermore, since  $F_\phi(x) \nolinebreak\subset\nolinebreak T_C(x)$, a key step towards verifying~\eqref{eqMotiv} consists in guaranteeing 
\begin{equation}\label{eq_transversality_C_partialKe}
    T_{\partial K}(x) \cap T_C(x) = T_{\partial K \cap C}(x) \qquad \forall x \in \mathcal{P}.
\end{equation}

Having \eqref{eqMotiv}, to ensure (as required in~\eqref{eqMotiv1}) that 
\begin{align} \label{eqaddnew} 
    \text{cl}(F_\phi(x)) \cap T_{\partial K \backslash \text{int}(C)}(x) = \emptyset \qquad \forall x \in \mathcal{P}, 
\end{align}
 it is enough to show that
\begin{equation} \label{eq_transversality_C_partialKe--} 
        \hspace{-0.4cm} T_{\partial K \cap C}(x)  \cap T_{\partial K \backslash \text{int}(C)}(x) 
        = T_{\partial K \cap \partial C}(x) \quad   \forall x \in \mathcal{P}, 
\end{equation}
and use the following assumption.

\begin{assumption} [Decisiveness]
\label{ass_decisiveness}
 \hspace{-0.25cm} $F \hspace{-0.05cm} \cap  \hspace{-0.05cm} T_{\partial K \cap \partial C} = \emptyset$ on $\mathcal{P}$. 
\end{assumption}

\begin{remark}

Assumption \ref{ass_decisiveness} is named decisiveness because it allows $F$ (through the first-order approximation of solutions) to indicate the specific side of $\partial K$ (either $\partial K \cap \text{int}(C)$ or $\partial K \backslash C$) the solutions must be, locally, more leaning towards, while leaving $K$ immediately. In other words, Assumption \ref{ass_decisiveness} implies \eqref{eqaddnew}
knowing that, for all $x \in \mathcal{P}$,
\begin{align*} 
F_\phi(x) 
\subset T_{\partial K \cap C}(x)  = \left(\mathbb{R}^n \backslash T_{\partial K \backslash \text{int}(C)}(x) \right) \cup T_{\partial K \cap \partial C}(x). 
\end{align*}

\end{remark}

To verify~\eqref{eq_transversality_C_partialKe} and \eqref{eq_transversality_C_partialKe--}, we use the next assumption.

\begin{assumption}[Transversality] \label{ass_transversality}~~\\
    - For each $x \in \mathcal{P}$, $T_C(x) = T^a_C(x)$.\\
    - For each $x \in \mathcal{P}$,  $T_{\partial K}(x) = T^a_{\partial K}(x)$, and\\   
    - For each $x \in \mathcal{P}$,
    \\
    ($\dagger$) $\exists$ a neighborhood of $x$, denoted by $\mathcal{N}(x)$, $c> 0$, and $\alpha \in [0,1)$ such that
    \begin{align*}
            & x_2 - x_1 \in v_1 - v_2 + \alpha (|x_1 - x_2|) \mathbb{B} \\
            & \forall (x_1,x_2) \in 
             ((\partial K \backslash C) \cap \mathcal{N}(x)) \times
                ((\partial C \backslash \partial K) \cap \mathcal{N}(x)),
    \end{align*}
    $(v_1,v_2) \in  T_{\partial K} (x_1) \times T_{C} (x_2)$ : $|\left[v_1^\top ~ v_2^\top \right]| \leq c |x_1 - x_2|$.
\end{assumption}

Note that the first item in Assumption~\ref{ass_transversality} is not needed to verify~\eqref{eq_transversality_C_partialKe} nor~\eqref{eq_transversality_C_partialKe--}. It is rather used to ensure that $T^a_{\partial K \cap C} (x) = T_{\partial K \cap C} (x)$ for all $x \in \mathcal{P}$, a useful property to guarantee  that~\eqref{eqMotiv1} implies ($\star$). 

The second and third items are used to establish~\eqref{eq_transversality_C_partialKe} and~\eqref{eq_transversality_C_partialKe--}, i.e., the intersection of contingent cones coincides with the contingent cone of the intersection. Conditions ensuring such a property are termed ``transversality'' conditions in~\cite{aubinSetValuedAnalysis2009, Aubin1991Viability}, and are usually more restrictive than the last two items in Assumption \ref{ass_transversality}; see~\cite[Remark 8]{arXiv_version_auto}. Roughly speaking, near a point in $x \in \mathcal{P}$, 
($\dagger$) requires that each pair of nearby points $(x_1,x_2)$ on the mismatched boundary pieces can be associated with tangent directions $(v_1,v_2)$ such that the difference $v_1-v_2$ approximates the displacement $x_2-x_1$ up to a contraction error of size $\alpha|x_1-x_2|$, while the tangent vectors remain bounded proportionally to the separation, i.e., $|[v_1^\top ~ v_2^\top]| \leq c|x_1-x_2|$; see Fig.~\ref{fig_dagger}.

\begin{figure}[t]
    \centering
    \begin{subfigure}[t]{0.48\columnwidth}
        \centering
        \resizebox{\linewidth}{!}{%
            \tikzset{every picture/.style={line width=0.75pt}} 

\begin{tikzpicture}[x=0.75pt,y=0.75pt,yscale=-1,xscale=1]

\draw  [draw opacity=0][fill={rgb, 255:red, 126; green, 211; blue, 33 }  ,fill opacity=0.38 ] (290.89,26) -- (352.14,26) -- (352.14,148) -- (290.89,148) -- cycle ;
\draw  [draw opacity=0][fill={rgb, 255:red, 74; green, 144; blue, 226 }  ,fill opacity=0.61 ] (337.09,115.75) .. controls (336.28,115.25) and (351.09,119.25) .. (363.84,101) .. controls (376.59,82.75) and (380.09,55.5) .. (388.59,43) .. controls (397.09,30.5) and (411.89,31) .. (412.39,30.25) .. controls (412.89,29.5) and (412.64,146.5) .. (412.51,146.75) .. controls (412.38,147) and (336.78,149) .. (337.09,148) .. controls (337.39,147) and (337.89,116.25) .. (337.09,115.75) -- cycle ;
\draw    (337.09,115.75) .. controls (385.64,115) and (369.14,31) .. (412.39,30.25) ;
\draw  [draw opacity=0][fill={rgb, 255:red, 74; green, 144; blue, 226 }  ,fill opacity=0.61 ] (290.89,115.9) -- (337.52,115.9) -- (337.52,148) -- (290.89,148) -- cycle ;
\draw    (352.14,26) -- (352.14,148) ;
\draw    (291.14,115.75) -- (337.09,115.75) ;
\draw  [fill={rgb, 255:red, 0; green, 0; blue, 0 }  ,fill opacity=1 ] (388.96,40) .. controls (388.96,38.69) and (390.03,37.63) .. (391.34,37.63) .. controls (392.65,37.63) and (393.71,38.69) .. (393.71,40) .. controls (393.71,41.31) and (392.65,42.38) .. (391.34,42.38) .. controls (390.03,42.38) and (388.96,41.31) .. (388.96,40) -- cycle ;
\draw [color={rgb, 255:red, 144; green, 19; blue, 254 }  ,draw opacity=1 ]   (352.09,81) -- (380.75,51.19) -- (389.96,41.45) ;
\draw [shift={(391.34,40)}, rotate = 133.41] [color={rgb, 255:red, 144; green, 19; blue, 254 }  ,draw opacity=1 ][line width=0.75]    (10.93,-3.29) .. controls (6.95,-1.4) and (3.31,-0.3) .. (0,0) .. controls (3.31,0.3) and (6.95,1.4) .. (10.93,3.29)   ;
\draw [color={rgb, 255:red, 65; green, 117; blue, 5 }  ,draw opacity=1 ]   (351.84,43.5) -- (352.07,79) ;
\draw [shift={(352.09,81)}, rotate = 269.62] [color={rgb, 255:red, 65; green, 117; blue, 5 }  ,draw opacity=1 ][line width=0.75]    (10.93,-3.29) .. controls (6.95,-1.4) and (3.31,-0.3) .. (0,0) .. controls (3.31,0.3) and (6.95,1.4) .. (10.93,3.29)   ;
\draw  [fill={rgb, 255:red, 0; green, 0; blue, 0 }  ,fill opacity=1 ] (364.25,95.63) .. controls (364.25,94.31) and (365.31,93.25) .. (366.63,93.25) .. controls (367.94,93.25) and (369,94.31) .. (369,95.63) .. controls (369,96.94) and (367.94,98) .. (366.63,98) .. controls (365.31,98) and (364.25,96.94) .. (364.25,95.63) -- cycle ;
\draw [color={rgb, 255:red, 144; green, 19; blue, 254 }  ,draw opacity=1 ]   (352.39,122) -- (365.68,97.39) ;
\draw [shift={(366.63,95.63)}, rotate = 118.35] [color={rgb, 255:red, 144; green, 19; blue, 254 }  ,draw opacity=1 ][line width=0.75]    (10.93,-3.29) .. controls (6.95,-1.4) and (3.31,-0.3) .. (0,0) .. controls (3.31,0.3) and (6.95,1.4) .. (10.93,3.29)   ;
\draw [color={rgb, 255:red, 65; green, 117; blue, 5 }  ,draw opacity=1 ]   (351.86,93.34) -- (352.35,120) ;
\draw [shift={(352.39,122)}, rotate = 268.93] [color={rgb, 255:red, 65; green, 117; blue, 5 }  ,draw opacity=1 ][line width=0.75]    (10.93,-3.29) .. controls (6.95,-1.4) and (3.31,-0.3) .. (0,0) .. controls (3.31,0.3) and (6.95,1.4) .. (10.93,3.29)   ;
\draw  [fill={rgb, 255:red, 0; green, 0; blue, 0 }  ,fill opacity=1 ] (348.75,41.45) .. controls (349.1,40.01) and (350.56,39.13) .. (352,39.48) .. controls (353.44,39.83) and (354.32,41.28) .. (353.96,42.72) .. controls (353.61,44.16) and (352.16,45.04) .. (350.72,44.69) .. controls (349.28,44.34) and (348.4,42.89) .. (348.75,41.45) -- cycle ;
\draw  [fill={rgb, 255:red, 0; green, 0; blue, 0 }  ,fill opacity=1 ] (349.25,92.7) .. controls (349.6,91.26) and (351.06,90.38) .. (352.5,90.73) .. controls (353.94,91.08) and (354.82,92.53) .. (354.46,93.97) .. controls (354.11,95.41) and (352.66,96.29) .. (351.22,95.94) .. controls (349.78,95.59) and (348.9,94.14) .. (349.25,92.7) -- cycle ;
\draw [color={rgb, 255:red, 255; green, 255; blue, 255 }  ,draw opacity=1 ][line width=2.25]    (290.89,148) -- (413.39,148.5) ;
\draw [color={rgb, 255:red, 255; green, 255; blue, 255 }  ,draw opacity=1 ][line width=2.25]    (414.14,28.75) -- (413.39,148.5) ;

\draw (295.29,30.15) node [anchor=north west][inner sep=0.75pt]    {$K$};
\draw (334,37.15) node [anchor=north west][inner sep=0.75pt]  [font=\scriptsize]  {$x_{1}$};
\draw (391.34,42.4) node [anchor=north west][inner sep=0.75pt]  [font=\scriptsize]  {$x_{2}$};
\draw (370.5,90.65) node [anchor=north west][inner sep=0.75pt]  [font=\scriptsize]  {$\hat{x}_{2}$};
\draw (334,84.65) node [anchor=north west][inner sep=0.75pt]  [font=\scriptsize]  {$\hat{x}_{1}$};
\draw (354.96,62.56) node [anchor=north west][inner sep=0.75pt]  [font=\tiny,color={rgb, 255:red, 144; green, 19; blue, 254 }  ,opacity=1 ,rotate=-314.17]  {${\textstyle \textcolor[rgb]{0.56,0.07,1}{-v}\textcolor[rgb]{0.56,0.07,1}{_{\textcolor[rgb]{0.56,0.07,1}{2}}}}$};
\draw (338.86,55.03) node [anchor=north west][inner sep=0.75pt]  [font=\tiny,color={rgb, 255:red, 208; green, 2; blue, 27 }  ,opacity=1 ,rotate=-0.3]  {${\textstyle \textcolor[rgb]{0.25,0.46,0.02}{v}\textcolor[rgb]{0.25,0.46,0.02}{_{1}}}$};
\draw (336.61,100.28) node [anchor=north west][inner sep=0.75pt]  [font=\tiny,color={rgb, 255:red, 208; green, 2; blue, 27 }  ,opacity=1 ,rotate=-0.3]  {${\textstyle \textcolor[rgb]{0.25,0.46,0.02}{\hat{v}}\textcolor[rgb]{0.25,0.46,0.02}{_{1}}}$};
\draw (356.39,111.8) node [anchor=north west][inner sep=0.75pt]  [font=\tiny,color={rgb, 255:red, 144; green, 19; blue, 254 }  ,opacity=1 ,rotate=-358.78]  {${\textstyle \textcolor[rgb]{0.56,0.07,1}{-}\textcolor[rgb]{0.56,0.07,1}{\hat{v}}\textcolor[rgb]{0.56,0.07,1}{_{\textcolor[rgb]{0.56,0.07,1}{2}}}}$};
\draw (393.09,127.65) node [anchor=north west][inner sep=0.75pt]    {$C$};
\draw (294.09,132.15) node [anchor=north west][inner sep=0.75pt]  [font=\footnotesize]  {$K\ \cap C$};

\end{tikzpicture}
        }
        \caption{Intersection verifying ($\dagger$)}
        \label{fig_dagger_holds}
    \end{subfigure}
    \hfill
    \begin{subfigure}[t]{0.48\columnwidth}
        \centering
        \resizebox{\linewidth}{!}{%
            \tikzset{every picture/.style={line width=0.75pt}} 

\begin{tikzpicture}[x=0.75pt,y=0.75pt,yscale=-1,xscale=1]

\draw  [draw opacity=0][fill={rgb, 255:red, 74; green, 144; blue, 226 }  ,fill opacity=0.61 ] (369.04,108.9) .. controls (368.93,99.3) and (369.04,86.31) .. (371.89,70) .. controls (374.74,53.69) and (378.58,43.59) .. (382.64,38.75) .. controls (386.71,33.91) and (387.99,33.44) .. (394.64,31.25) .. controls (401.29,29.06) and (406.81,28.45) .. (410.89,28.25) .. controls (414.98,28.05) and (428.82,28.19) .. (428.64,28.25) .. controls (428.47,28.31) and (429.64,150) .. (429.39,149.25) .. controls (429.14,148.5) and (369.39,149) .. (369.14,149.15) .. controls (368.88,149.3) and (369.14,118.5) .. (369.04,108.9) -- cycle ;
\draw    (369.04,108.9) .. controls (369.36,60.63) and (377.83,40.08) .. (388.49,33.54) .. controls (399.14,27) and (412.47,28.13) .. (428.39,28.25) ;
\draw  [draw opacity=0][fill={rgb, 255:red, 126; green, 211; blue, 33 }  ,fill opacity=0.38 ] (306.03,27.25) -- (369.14,27.25) -- (369.14,150.15) -- (306.03,150.15) -- cycle ;
\draw  [draw opacity=0][fill={rgb, 255:red, 74; green, 144; blue, 226 }  ,fill opacity=0.61 ] (306.14,108.9) -- (369.04,108.9) -- (369.04,150.15) -- (306.14,150.15) -- cycle ;
\draw    (369.14,27.25) -- (368.96,120.15) -- (368.89,150.15) ;
\draw  [fill={rgb, 255:red, 0; green, 0; blue, 0 }  ,fill opacity=1 ] (380.04,38.11) .. controls (380.39,36.67) and (381.84,35.79) .. (383.28,36.14) .. controls (384.72,36.5) and (385.6,37.95) .. (385.25,39.39) .. controls (384.9,40.83) and (383.44,41.71) .. (382,41.36) .. controls (380.56,41) and (379.68,39.55) .. (380.04,38.11) -- cycle ;
\draw [color={rgb, 255:red, 65; green, 117; blue, 5 }  ,draw opacity=1 ]   (368.89,36.5) -- (368.67,54.5) ;
\draw [shift={(368.64,56.5)}, rotate = 270.7] [color={rgb, 255:red, 65; green, 117; blue, 5 }  ,draw opacity=1 ][line width=0.75]    (10.93,-3.29) .. controls (6.95,-1.4) and (3.31,-0.3) .. (0,0) .. controls (3.31,0.3) and (6.95,1.4) .. (10.93,3.29)   ;
\draw  [fill={rgb, 255:red, 0; green, 0; blue, 0 }  ,fill opacity=1 ] (366.14,36.5) .. controls (366.14,34.98) and (367.37,33.75) .. (368.89,33.75) .. controls (370.41,33.75) and (371.64,34.98) .. (371.64,36.5) .. controls (371.64,38.02) and (370.41,39.25) .. (368.89,39.25) .. controls (367.37,39.25) and (366.14,38.02) .. (366.14,36.5) -- cycle ;
\draw [color={rgb, 255:red, 144; green, 19; blue, 254 }  ,draw opacity=1 ]   (368.64,56.5) -- (381.4,40.32) ;
\draw [shift={(382.64,38.75)}, rotate = 128.26] [color={rgb, 255:red, 144; green, 19; blue, 254 }  ,draw opacity=1 ][line width=0.75]    (10.93,-3.29) .. controls (6.95,-1.4) and (3.31,-0.3) .. (0,0) .. controls (3.31,0.3) and (6.95,1.4) .. (10.93,3.29)   ;
\draw    (305.89,109.25) -- (369.29,109.4) ;
\draw  [fill={rgb, 255:red, 0; green, 0; blue, 0 }  ,fill opacity=1 ] (366.79,97.1) .. controls (367.14,95.66) and (368.59,94.78) .. (370.03,95.13) .. controls (371.47,95.48) and (372.35,96.93) .. (372,98.37) .. controls (371.65,99.81) and (370.2,100.69) .. (368.76,100.34) .. controls (367.32,99.99) and (366.44,98.54) .. (366.79,97.1) -- cycle ;
\draw [color={rgb, 255:red, 65; green, 117; blue, 5 }  ,draw opacity=1 ]   (368.48,87.88) -- (368.78,140.3) ;
\draw [shift={(368.79,142.3)}, rotate = 269.68] [color={rgb, 255:red, 65; green, 117; blue, 5 }  ,draw opacity=1 ][line width=0.75]    (10.93,-3.29) .. controls (6.95,-1.4) and (3.31,-0.3) .. (0,0) .. controls (3.31,0.3) and (6.95,1.4) .. (10.93,3.29)   ;
\draw  [fill={rgb, 255:red, 0; green, 0; blue, 0 }  ,fill opacity=1 ] (365.73,90.88) .. controls (365.73,89.36) and (366.96,88.13) .. (368.48,88.13) .. controls (370,88.13) and (371.23,89.36) .. (371.23,90.88) .. controls (371.23,92.39) and (370,93.63) .. (368.48,93.63) .. controls (366.96,93.63) and (365.73,92.39) .. (365.73,90.88) -- cycle ;
\draw [color={rgb, 255:red, 144; green, 19; blue, 254 }  ,draw opacity=1 ]   (368.79,142.3) -- (369.37,99.73) ;
\draw [shift={(369.4,97.74)}, rotate = 90.78] [color={rgb, 255:red, 144; green, 19; blue, 254 }  ,draw opacity=1 ][line width=0.75]    (10.93,-3.29) .. controls (6.95,-1.4) and (3.31,-0.3) .. (0,0) .. controls (3.31,0.3) and (6.95,1.4) .. (10.93,3.29)   ;
\draw [color={rgb, 255:red, 255; green, 255; blue, 255 }  ,draw opacity=1 ][line width=2.25]    (306.14,150.15) -- (430.14,150.5) ;
\draw [color={rgb, 255:red, 255; green, 255; blue, 255 }  ,draw opacity=1 ][line width=2.25]    (430.39,26) -- (430.14,150.5) ;

\draw (350.87,29.8) node [anchor=north west][inner sep=0.75pt]  [font=\scriptsize,xslant=0.01]  {$x_{1}$};
\draw (352.86,43.28) node [anchor=north west][inner sep=0.75pt]  [font=\tiny,color={rgb, 255:red, 208; green, 2; blue, 27 }  ,opacity=1 ,rotate=-0.3]  {${\textstyle \textcolor[rgb]{0.25,0.46,0.02}{v}\textcolor[rgb]{0.25,0.46,0.02}{_{1}}}$};
\draw (408.09,131.15) node [anchor=north west][inner sep=0.75pt]    {$C$};
\draw (311.03,31.65) node [anchor=north west][inner sep=0.75pt]    {$K$};
\draw (350.12,83.55) node [anchor=north west][inner sep=0.75pt]  [font=\scriptsize,xslant=0.01]  {$\hat{x}_{1}$};
\draw (370.39,115.3) node [anchor=north west][inner sep=0.75pt]  [font=\tiny,color={rgb, 255:red, 144; green, 19; blue, 254 }  ,opacity=1 ,rotate=-358.78]  {${\textstyle \textcolor[rgb]{0.56,0.07,1}{-}\textcolor[rgb]{0.56,0.07,1}{\hat{v}}\textcolor[rgb]{0.56,0.07,1}{_{\textcolor[rgb]{0.56,0.07,1}{2}}}}$};
\draw (355.11,115.78) node [anchor=north west][inner sep=0.75pt]  [font=\tiny,color={rgb, 255:red, 208; green, 2; blue, 27 }  ,opacity=1 ,rotate=-0.3]  {${\textstyle \textcolor[rgb]{0.25,0.46,0.02}{\hat{v}}\textcolor[rgb]{0.25,0.46,0.02}{_{1}}}$};
\draw (310.09,132.4) node [anchor=north west][inner sep=0.75pt]  [font=\footnotesize]  {$K\ \cap C$};
\draw (372.64,55.25) node [anchor=north west][inner sep=0.75pt]  [font=\tiny,color={rgb, 255:red, 144; green, 19; blue, 254 }  ,opacity=1 ,rotate=-316.18]  {${\textstyle \textcolor[rgb]{0.56,0.07,1}{-v}\textcolor[rgb]{0.56,0.07,1}{_{\textcolor[rgb]{0.56,0.07,1}{2}}}}$};
\draw (385.28,36.29) node [anchor=north west][inner sep=0.75pt]  [font=\scriptsize]  {$x_{2}$};
\draw (373.37,91.55) node [anchor=north west][inner sep=0.75pt]  [font=\scriptsize,xslant=0.01]  {$\hat{x}_{2}$};

\end{tikzpicture}
        }
        \caption{($\dagger$) is not verified}
        \label{fig_dagger_doesnt_hold}
    \end{subfigure}
    \caption{Each sub-figure depicts a pair of sets $(K,C)$ and an intersection point in $\partial K \cap \partial C$. In each case, two pair of points $(x_1,x_2), (\hat{x}_1,\hat{x}_2) \in \partial K \setminus C \times  \partial C \setminus \partial K$  are shown, with the pair $(\hat{x}_1,\hat{x}_2)$ being closer to the intersection point. We illustrate, in each case, the evolution of norm-smallest pairs $(v_1, v_2) \in T_{\partial K}(x_1) \times  T_C(x_2)$ and $(\hat{v}_1, \hat{v}_2) \in T_{\partial K}(\hat{x}_1) \times  T_C(\hat{x}_2)$ such that $x_2-x_1 = v_1-v_2$ and $\hat{x}_2-\hat{x}_1 = \hat{v}_1 - \hat{v}_2$. 
    In the left figure, we observe that the ratio 
    $|[v_1^\top ~~ v_2^\top]| / |x_1 - x_2|$ remains bounded while approaching the intersection point, which is not the case in the right sub-figure, due a tangent separation.}
    \label{fig_dagger}
\end{figure}

\ifitsdraft

\else

We include, in \cite[Example 3]{arXiv_version_auto}, an example illustrating the different assumptions. Furthermore, we gave in the conference version  \cite{reynaudNagumoTypeCharacterizationForward2025} examples where  Assumptions \ref{ass_propagation}-\ref{ass_decisiveness} (resp., Assumptions \ref{ass_decisiveness}-\ref{ass_transversality}) hold but Assumption \ref{ass_transversality} (resp., Assumption \ref{ass_propagation}) does not, and showed in each case that \eqref{eqconstNag} is not enough to guarantee forward invariance of $K$. These examples are recalled in the report version \cite{arXiv_version_auto}, in addition to a new example demonstrating the importance of Assumptions \ref{ass_decisiveness}; see \cite[Sections 5.1-5.3]{arXiv_version_auto}.

Before stating our result for general sets, the following remarks are in order. 

\fi

\begin{remark} \label{rem_prop_is_given}
    We show in
    \ifitsdraft
        Lemma ~\ref{lemAp} in the Appendix
    \else
       \cite[Lemma 21]{arXiv_version_auto} 
    \fi
that the second item in Assumption~\ref{ass_transversality} holds under the following property.
     \begin{itemize}
        \item[(Pr1)] The map $T_{K} : \partial K \rightrightarrows \mathbb{R}^n$ is continuous on $\mathcal{P}$, and 
        $$ x \in \text{cl}(\partial K \cap \text{int}(C)) \qquad  \forall x \in \mathcal{P}. $$ 
    \end{itemize} 
    Moreover, Assumption~\ref{ass_propagation} also holds under (Pr1) when additionally  \eqref{eqconstNag} holds and $F \cap T_C$ is continuous on $\mathcal{P}$.   
\end{remark}

\ifitsdraft
\else
\begin{remark}
According to Lemma~\ref{lemtranssimple} in the Appendix, ($\dagger$) in Assumption~\ref{ass_transversality} is verified, at a given $x \in \partial K \cap \partial C$, if 
\begin{itemize}
    \item[($\ddagger$)] $\exists \mathcal{N}(x)$, a neighborhood of $x$  such that 
    \begin{equation*} 
        \begin{aligned}
             x_2 - x_1 & \in 
            C_{\partial K} (x) - C_{C} (x) 
            \\ &
            := \{ v - w : (v, w) \in C_{\partial K} (x) \times C_{C} (x) \}  
          \\  \forall (x_1,x_2) & \in 
                ((\partial K \backslash C) \cap \mathcal{N}(x)) \times  ((\partial C \backslash \partial K) \cap \mathcal{N}(x)).
            \end{aligned}
    \end{equation*}
\end{itemize}
Although ($\ddagger$) is in general more restrictive than ($\dagger$), 
see ~\cite[Example 7]{arXiv_version_auto},
it will be used later to verify ($\dagger$) when the sets $(C,K)$ are finitely represented, using simpler and more tractable conditions.
\end{remark}
\fi 

\ifitsdraft
The following example illustrates the assumptions. 
\begin{example} \label{example_assumptions1}
    Consider the sets $K := \{ x \in \mathbb{R}^3 : x_1 \leq 1\}$ and $C := \{ x \in \mathbb{R}^3 : x_1^2 + x_2^2 \geq 1, x_1 \geq 0, x_2 \leq 0\}$, and let  
  $F(x):= [-2, 0] \times 
    [-x_1-2, -x_1] \times \{2\}$ for all $x \in C$.
           
The considered Nagumo  condition~\eqref{eqconstNag} holds since, for all $x\in K\cap \mathrm{int}(C)$,
\begin{equation}\label{eq_F_at_P_ex}
\begin{aligned}
 F(x) & \subset  [-1, 0] \times [-3, -1] \times \{ 2 \} \\ & \subset T_K(x) = \{v\in\mathbb{R}^3:v_1\leq 0\}.
        \end{aligned}
    \end{equation}    
    
Having 
    $\partial K \cap \partial C = 
       \{ x \in \mathbb{R}^3 : x_1=1,\;x_2=0 \}, $
       we can easily verify that $\partial K \cap \partial C = \partial K \cap \partial C \cap \text{cl}(C \backslash K) $. Moreover, we can show that $C$ is not convex at $x$, $F(x) \not\subset 
    D_{K}(x)$, and $F(x) \not\subset \mathbb{R}^n \backslash T_C(x)$ whenever $x \in \partial K \cap \partial C$. 
    Hence, we propose to verify the three assumptions on $\mathcal{P}_o = \partial K \cap \partial C$, which is an over approximation of the set $\mathcal{P}$, obtained by  ignoring the third negated property in its definition. 

    Since (Pr1) in Remark~\ref{rem_prop_is_given} holds, Assumption~\ref{ass_propagation} follows. We can also verify ~\eqref{eq_F_at_P_ex}  for all $x\in\mathcal P_o$. Moreover, since
    \[ T_{\partial K\cap\partial C}(x) 
        = 
\{v\in\mathbb{R}^3:v_1=0,\;v_2=0\} \qquad \forall x\in\mathcal P_o, \]
    \eqref{eq_F_at_P_ex}, verified on $\mathcal{P}_o$, implies that Assumption~\ref{ass_decisiveness} holds.

    The first two items of Assumption~\ref{ass_transversality} also hold, since
    \[
        \begin{aligned}
            T_{\partial K}(x) &= T_{\partial K}^{a}(x) = \{v\in\mathbb{R}^3:v_1=0\}, \\
            T_C(x) &= T_C^{a}(x) = \{v\in\mathbb{R}^3:v_1\geq 0,\;v_2\leq 0\},
        \end{aligned}
        ~~ \forall x\in\mathcal P_o.
    \]
    
    Finally, we verify  ($\dagger$) in Assumption~\ref{ass_transversality} for any $x \in \mathcal P_o$. Indeed, given $x \in \mathcal P_o$ for any neighborhood $\mathcal N(x)$,
    \[
        T_{\partial K}(x_1) = 
        \{ v \in \mathbb{R}^3 : v_1=0 \} \qquad \forall x_1\in(\partial K\setminus C)\cap\mathcal N(x).
    \]
    Moreover, for any $x_2\in(C\setminus\partial K)\cap\mathcal N(x)$, one of the following must hold.
    \begin{itemize}
        \item If $x_2\in K$, 
        $
        T_C(x_2) \supset \{v\in\mathbb{R}^3:v_1\geq 0,\;v_2=0\} $
         and $x_2-x_1=[u_1 ~~ u_2 ~~ u_3]^\top$ with $u_1\leq 0$, $u_2\leq 0$.    
        \item If $x_2\notin K$,  $T_C(x_2) = \{v\in\mathbb{R}^3:v_2\leq 0\}$ and $x_2-x_1=[u_1  ~~ u_2  ~~ 0]^\top$ with $u_1\geq 0$, $u_2\leq 0$ and $u_3\in \mathbb{R}$.
    \end{itemize}
    In both cases, condition ($\dagger$) holds with $v_1=[0 ~~ u_2  ~~ u_3]^\top$, $v_2=[-u_1 ~~ 0  ~~ 0]^\top$ and $\alpha=0$.
    
    Since Assumption~\ref{ass_transversality} holds, relation~\eqref{eq_transversality_C_partialKe} follows. Indeed, for all $x\in\mathcal P_o$,
    \[
        T_{\partial K}(x)\cap T_C(x) = T_{\partial K\cap C}(x) = \{v\in\mathbb{R}^3:v_1=0,\;v_2\leq 0\}.
    \]
        It is also straightforward~\eqref{eq_transversality_C_partialKe--}.
    
   Finally, we find an explicit solution $\phi$ to $\Sigma$ starting from $\phi(0)=[1 ~~ 0 ~~ 0]^\top\in \mathcal{P}_o$ verifying  $\text{cl}(F_\phi(x)) \cap T_{\partial K \backslash \text{int}(C)}(x) = \emptyset$, as stated in~\eqref{eqaddnew}. Indeed, we consider $\phi(t) = [1-t^2 ~~ -2t  ~~~ 2t\,]^\top$. Then $\phi(\operatorname{dom}\phi)\subset C$ and
    \[
        \dot\phi(t) = [-2t ~~ -2 ~~~ 2]^\top \in F(\phi(t)) \qquad \forall t\in\operatorname{dom}\phi.
    \]
    Hence, $\phi$ is a solution to $\Sigma$. Moreover, since $\phi$ is differentiable at $t=0$, then
        $F_\phi(x)= \{\dot{\phi}(0)\} = \{[0 ~~ -2 ~~ 2]^\top\}$, and thus
   $ F_\phi(x) \cap  T_{\partial K \backslash \text{int}(C)}(x) =
            \{[0 ~~ -2 ~~ 2]^\top\} \cap \{v\in\mathbb{R}^3:v_1=0,\;v_2\geq 0\} = \emptyset.$ 
\end{example}
\fi

\section{Result for general sets}\label{sec_main_result}

We are now ready to state a general version of  Nagumo's invariance theorem for constrained systems. 

\begin{theorem} \label{theo_main_result}
Let system $\Sigma$ and a closed set $K \subset \mathbb{R}^n$.  
\begin{itemize} 
\item When $F$ is lower semicontinuous, forward invariance of $K$ $\Rightarrow$~\eqref{eqconstNag}. 
\item Under Assumptions~\ref{ass_propagation}-\ref{ass_transversality} and for $F$ one-sided locally Lipschitz, 
\eqref{eqconstNag} $\Rightarrow$ forward invariance of $K$.
\end{itemize}
\end{theorem}

\begin{proof}
The proof of the first bullet is the same as in the unconstrained case \cite[Theorem 5.3.4]{Aubin1991Viability};
\ifitsdraft
hence, it is reported to the Appendix.
\else
see \cite{arXiv_version_auto}. 
\fi

To prove the second bullet, we proceed by contradiction. That is, we let a solution $\phi$ starting from $x \in \partial K \cap C$ that leaves $K$ immediately.  We next distinguish  complementary situations. 
\\
- When $x \in \mathcal{P}$. 
    We use Proposition~\ref{prop_clFphi_supersets}, under the three assumptions (propagation, decisiveness, and transversality), to conclude that the inclusion in \eqref{eqMotiv1} holds. As a result, Proposition~\ref{prop_leaving_implies_star}, under the transversality assumption (Assumption~\ref{ass_transversality}),  can be used to conclude that $\phi$ must satisfy ($\star$). Now, for $S = K$ and using \eqref{eqconstNag}, we conclude that condition  
    \eqref{eqS} in   Proposition~\ref{prop_star_implies_not_leaving} holds under ($\star$). However, the latter proposition guarantees $\phi([0,T]) \subset S = K$, for some $T > 0$, which leads to a contradiction. 
  \\  
- When $x \in  \text{cl}(C \backslash K) \cap \partial K \cap \partial C$ and $F(x) \cap T_C(x) = \emptyset$. Using \cite[Proposition 3.4.1.]{Aubin1991Viability}, we conclude that $\text{dom} \phi = \{0\}$, which leads to a contradiction.    
 \\
-  When  $ x \in  \text{cl}(C \backslash K) \cap \partial K \cap \partial C$ and $F(x) \in D_K(x)$. 
We use \cite[Theorem 4.3.4]{Aubin1991Viability} to conclude that every nontrivial solution, starting from $x$, must remain in $K$ for a nontrivial interval of time. 
 \\
- When  $ x \in  \text{cl}(C \backslash K) \cap \partial K \cap \partial C$ and 
    there exists $\epsilon > 0$ and a closed subset $S \subset K$ containing $x$ such that 
\begin{align} \label{eqnewnew}
    \widehat{F}(y) \subset T_{S}(y) \quad  \forall y \in \mathbb{B}_\epsilon(x) \cap S.
\end{align}    
  We use Proposition~\ref{prop_star_implies_not_leaving} applied to 
    $\widehat{\Sigma} : \dot{x} \in \widehat{F}(x)$ defined on $C \cup K$, instead of $\Sigma$.
    By continuity of solutions, for every solution $\phi$ to $\widehat{\Sigma}$ starting from $x$, there exists $T > 0$ such that $\phi([0,T]) \subset \mathbb{B}_\epsilon(x)$. Hence, 
    $ \text{Proj}_{S} (\phi([0,T])) \subset \mathbb{B}_\epsilon(x) \cap S$. Thus, condition \eqref{eqS} in Proposition \ref{prop_star_implies_not_leaving}, with $\widehat{F}$ instead of $F$, is verified under \eqref{eqnewnew}. Consequently, the solution $\phi$ to $\widehat{\Sigma}$ cannot immediately leave $S \subset K$; thus, so does every solution to $\Sigma$ starting from $x$.  
\\ 
- When  
$x \in \text{cl}(C \backslash K) \cap \partial K \cap \partial C$, $C$ is convex at $x$, and
$F(y) \subset T_{C \cap K}(y)$
for all $y \in \mathbb{B}_{\epsilon}(x) \cap \partial K \cap \partial C$. We apply Proposition \ref{prop_star_implies_not_leaving} with $S = K \cap C$. To do so, it enough to let a solution $\phi$ starting from $x$ that leaves $K$ immediately and show that  
$\text{Proj}_{K \cap C} (\phi([0,T])) \subset (\partial K \cap  C)$ for some $T > 0$. To find a contradiction, we let $t > 0$ such that $y \in  \text{int}(K) \cap C$ for some $y \in \text{Proj}_{K \cap C} (\phi(t))$. By convexity of $C$ at $x$, the line segment $[y, \phi(t)]$ is subset of $C$ (when $t > 0$ is small), with $y \in \text{int}(K)$ and $\phi(t) \in \mathbb{R}^n \backslash K$. Hence, there exists $\lambda \in (0,1)$
such that
$z := \lambda y + (1-\lambda) \phi(t) \in \partial K \cap C$. At the same time, $\phi(t) - z = \lambda (\phi(t) - y) $, which contradicts $y$ being a projection of $\phi(t)$ on $K \cap C$.   
\\
- When $x \not\in  \text{cl}(C \backslash K) \cap \partial K \cap \partial C$.  One of the following must be true in this case.
    \begin{itemize}
     \item $x \notin \partial K$, which contradicts $\phi$ leaving $K$ immediately, by continuity of solutions.    
     \item $x \notin \partial C$, which means that $x \in \partial K \cap \text{int}(C)$. In this case, continuity of solutions allows us to guarantee ($\star$), and thus Proposition \ref{prop_star_implies_not_leaving} can be used, under \eqref{eqS} and for $S = K$, to exclude  leaving $K$ immediately. 
     
     \item $x \notin \text{cl}(C \backslash K)$. 
     This means the existence of neighborhood $\mathcal{N}_x$ of $x$ such that $ \mathcal{N}_x \cap (C \backslash K) = \emptyset$.  Hence, a continuous solution starting from $\partial K \cap \partial C$ cannot flow immediately to $C \backslash K$. 
\end{itemize}  
\end{proof}

\subsection{The Intermediate Results}

We now prove the intermediate results.

\begin{proposition}\label{prop_clFphi_supersets}
    Consider system $\Sigma$ and a closed set 
    $K \subset \mathbb{R}^n$  such that Assumptions~\ref{ass_propagation}-\ref{ass_transversality} hold. 
    Let $\phi$ be a solution starting from $x \in \mathcal{P}$ that leaves the set $K$ immediately. Then, $ \text{cl}(F_\phi(x)) \subset T_{\partial K \cap C}(x) \backslash T_{\partial K \backslash \text{int}(C)}(x)$.
\end{proposition}

\begin{proof}
    We first prove, under (SA), that
    \begin{align} 
        \text{cl} (F_{\phi}(x)) & \subset F(x),  \label{eq_clFphi_in_F}\\ 
        \text{cl} (F_{\phi}(x)) & \subset T_{C \backslash K} (x) \subset T_{\mathbb{R}^n \backslash K}(x). \label{eq_clFphi_in_TCnotKe}
    \end{align}
    
    Indeed,  since $F(x)$ and $T_{C \backslash K}(x)$ are closed, we show that $v \in F(x)$ and $v \in T_{C \backslash K}(x)$ for all $v \in F_{\phi}(x)$.  To do so, we let a sequence $\{ t_i \}^\infty_{i=1} \subset \mathbb{R}_{\geq 0}$ satisfying $\lim_{i \rightarrow \infty} t_i = 0$ such that 
    $ v = \lim_{i \rightarrow \infty} (\phi(t_i) - \phi(0))/t_i \in F_{\phi}(x). $ 
    
    To show that $v \in F(x)$, thanks to the concept of solutions, we introduce a sequence $\{ \delta^j_i \}^{N_i}_{j = 1} \subset [0,t_i]$, with $N_i \in \{1,2,\ldots,\infty\}$, such that 
    \begin{itemize}
        \item $t \mapsto \phi(t)$ is differentiable on each interval $(\delta^j_i,\delta^{j+1}_i)$.
        \item $\dot{\phi}(s) \in F(\phi(s))$ for all $s \in (\delta^j_i,\delta^{j+1}_i)$.
        \item $\delta^{1}_{i} = 0$ and 
        $\delta^{N_i}_{i} = t_i$.
    \end{itemize}
     
    Now, using the Mean-Value Theorem, we conclude the existence of $c^j_i \in ( \delta^{j}_i, \delta^{j+1}_i)$ such that
    \begin{equation*}
        \begin{aligned}
            \phi(\delta^{j+1}_i) - \phi(\delta^{j}_i) & = 
             \dot{\phi}(c^j_i) \left[ \delta^{j+1}_i \hspace{-0.2cm} - \delta^{j}_i \right]  \in F(\phi(c^j_i)) \left[ \delta^{j+1}_i  \hspace{-0.2cm} - \delta^{j}_i \right].
        \end{aligned}
    \end{equation*}
    The latter entails, for $v_i := 
    (\phi(t_i) - \phi(0))/t_i$, that
    $$ v_i \in \bigcup_{j \in \{1,\ldots,N_i\}} \left\{ F(\phi(c^j_i)) \right\} \subset \bigcup_{s \in [0,t_i]} \left\{ F(\phi(s)) \right\}.  $$
    
    Now, since $F$ is upper semicontinuous and has closed images, it follows from~\cite[Proposition 1.4.8]{aubinSetValuedAnalysis2009} that the graph of $F$ is closed. This allows us to conclude that $\lim_{t_i \rightarrow 0} \bigcup_{s \in [0,t_i]} \left\{ F(\phi(s)) \right\}  \subset F(x)$. Hence, $ v = \lim_{i \rightarrow \infty}  (\phi(t_i) - \phi(0))/t_i  \in F(x)$.
    
    To show that $v \in T_{C \backslash K}(x)$, we use the fact that $\phi$ leaves the set $K$ immediately to conclude that for some $N \in \mathbb{N}$ big enough,
    \begin{equation}\label{eq_phi_ti_leaves}
        \phi(t_i) = x + t_i v_i  \in C \backslash K \qquad \forall i \in \{N,N+1,\ldots\},
    \end{equation}
    where $ v_i := (\phi(t_i)- \phi(0))/t_i$. \eqref{eq_phi_ti_leaves} is enough to conclude that $v \in T_{C \backslash K}(x) \subset T_{C}(x)$. 
    
    Now,  using \eqref{eq_clFphi_in_F}-\eqref{eq_clFphi_in_TCnotKe}, we prove that
    \[
        \text{cl}(F_\phi(x)) \subset T_{\partial K \cap C}(x) \backslash T_{\partial K \backslash \text{int}(C)}(x). 
    \]
    
    On one hand,~Assumption~\ref{ass_propagation} is used to obtain  
    \begin{equation}\label{eqprf} 
        \text{cl} (F_\phi(x)) \subset F(x) \subset T_{K}(x).  
    \end{equation}
    On the other hand,~\cite[Theorem 4.3.3]{Aubin1991Viability} entails
    \begin{equation}\label{eq_quincapoix} 
        T_{\partial K}(x) = T_{K}(x) \cap T_{\mathbb{R}^n \backslash K}(x).
    \end{equation}
    
    Hence, combining~\eqref{eq_clFphi_in_TCnotKe} and~\eqref{eqprf}-\eqref{eq_quincapoix},  we obtain 
    \begin{align}\label{eq_clFphi_subset_TparK}
        \text{cl} (F_{\phi}(x)) & \subset F(x) \cap T_{C \backslash K}(x) \\ &    
        \subset T_K(x) \cap T_{\mathbb{R}^n \backslash K}(x) =  T_{\partial K}(x).
    \end{align}
    
    Next, using Lemma~\ref{lemtrans} in the Appendix with $(K_1,K_2)$ therein replaced by $(\partial K, C)$, under Assumption~\ref{ass_transversality}, we conclude that 
$T_{\partial K}(x) \cap  T_{C}(x) 
            =  T_{\partial K \cap C}(x)$.

    This allows us to conclude from~\eqref{eq_clFphi_subset_TparK} and~\eqref{eq_clFphi_in_TCnotKe} that $\text{cl} (F_{\phi}(x)) \subset   T_{\partial  K \cap C}(x)$. Now, using Assumption~\ref{ass_decisiveness},
    we conclude that $ F(x) \cap T_{\partial K \cap \partial C}(x) = \emptyset $, leading to 
    \begin{align*}
        \text{cl} (F_{\phi}(x)) \subset   T_{\partial  K \cap C}(x) \backslash T_{\partial K \cap \partial C}(x).  
    \end{align*}
    At this point,~\eqref{eqMotiv1} would follow if we show that 
    \begin{equation}\label{eqtranscons} 
        T_{\partial  K \cap C}(x)   \cap T_{\partial K \backslash \text{int}(C)}(x) = T_{\partial K \cap \partial C}(x). 
    \end{equation} 
    To prove the latter equality, it is enough to show that 
    \begin{equation}\label{eqtranscons+} 
        T_{C}(x) \cap T_{\partial K \backslash \text{int}(C)}(x) = T_{\partial K \cap \partial C}(x). 
    \end{equation}
    For this, we use Lemma~\ref{lemtrans}, while replacing $(K_1,K_2)$ therein by $(C,\partial K \backslash \text{int}(C))$, to conclude that~\eqref{eqtranscons+} holds provided that $T_C(x) = T^a_C(x)$ and 
    \begin{itemize}
        \item there exist $c> 0$, $\alpha \in [0,1)$, and a neighborhood $\mathcal{N}(x)$ around $x$ such that
        \begin{equation*}
            \begin{split}
                & \left[ (x_1 -  x_2)  + \alpha |x_1 - x_2| \mathbb{B} \right] \cap  (v_2 - v_1) \neq \emptyset, 
                \\  & 
                \forall (x_1,x_2) \in \begin{aligned}[t]
                        &([\partial C \backslash \partial K] \cap \mathcal{N}(x)) \times (( \partial K \backslash C) \cap \mathcal{N}(x)),
                    \end{aligned}
            \end{split}
        \end{equation*}
        for some $(v_1,v_2) \in T_{C} (x_1) \times T_{\partial K \backslash \text{int}(C)} (x_2)$ with $|\left[v_1^\top, v_2^\top \right]| \leq c |x_1 - x_2|$.
    \end{itemize}
    The latter item is verified under ($\dagger$) in Assumption~\ref{ass_transversality} since 
    $T_{\partial K \backslash \text{int}(C)} (x_2) = T_{\partial K} (x_2)$  for all $x_2 \in  \partial K \backslash C$.
\end{proof}

\begin{proposition} \label{prop_leaving_implies_star}
Consider $\Sigma$ and a closed set $ K \subset \mathbb{R}^n$ such that Assumption~\ref{ass_transversality} (transversality) holds.  Then, \eqref{eqMotiv1} implies ($\star$) for every solution $\phi$ starting from $x \in \mathcal{P}$.
\end{proposition}

\begin{proof}
    We show that
    \begin{equation}
    \label{eq_sequence_ti} 
        \begin{aligned}
            \forall \{ t_i \}^\infty_{i=1} & \subset \mathbb{R}_{\geq 0} ~  \text{with} ~ \lim_{i \rightarrow \infty} t_i = 0, 
            \\
            \exists \{ r_i \}^\infty_{i=1} & \subset \mathbb{R}_{\geq 0} ~  \text{with} ~  \lim_{i \rightarrow \infty} r_i = 0: 
            \\
            \phi(t_i) & \in x + t_i F_{\phi}(x) + t_i r_i \mathbb{B} \quad \forall i \in \{1,2,\ldots\}. 
        \end{aligned}
    \end{equation}
    The latter would entail that 
   $\phi(t) \in x + t F_\phi(x) + t o(t) \mathbb{B}$,
    with $\lim_{t \rightarrow 0^+} o(t) = 0$. From here, Lemma~\ref{lem_first_order_leaving_implies_star} in the Appendix can be used to conclude that $(\star)$ holds.
    
    To show~\eqref{eq_sequence_ti}, we use the fact that $F$ is locally bounded an we show that every sequence $\{ t_i \}^\infty_{i=1}$ as in~\eqref{eq_sequence_ti} leads to a sequence 
    $\left\{ (\phi(t_i) - \phi(0))/t_i \right\}^\infty_{i=1}$ that is bounded. Indeed, by absolute continuity of solutions, we have 
    $\phi(t_i) - \phi(0) = \int^{t_i}_{0} \dot{\phi}(s) ds$ for all $t_i \geq 0$, 
    which leads to 
    $ |\phi(t_i) - \phi(0)| \leq  t_i ~ \text{esssup}_{[0,t_i]} |\dot{\phi}(s)| $   for all $t_i \geq 0$.   
    Hence,  by the concept of solutions,  for all $t_i \geq 0$,
    $$ \big|(\phi(t_i) - \phi(0))/t_i \big| \leq  \text{sup} \left\{ |v| : v \in F(\phi(s)), ~ s \in [0,t_i] \right\}.$$
    
    Furthermore, by definition of $F_\phi(x)$, all the converging subsequences of 
    $ \left\{ (\phi(t_i) - \phi(0))/t_i \right\}^{\infty}_{i=1}$ have their limit in $F_{\phi}(x)$. 
    Hence,  there exists a sequence $\{ r_i \}_i \subset \mathbb{R}_{\geq 0}$ converging to $0$ such that 
    $$ (\phi(t_i) - \phi(0))/t_i \in F_{\phi}(x) + r_i \mathbb{B} \qquad \forall i \in \{1,2,\ldots\}. $$
\end{proof}

\begin{proposition} \label{prop_star_implies_not_leaving}
   Consider system $\Sigma$ with $F$ 
   one-sided locally Lipschitz and a closed set $S \subset \mathbb{R}^n$. 
    Consider a solution $\phi$ such that, for some $T > 0$,  we have 
   \begin{align}  \label{eqS} 
   \forall t \in (0,T], ~ \exists y \in  \text{Proj}_S (\phi(t)) ~ \text{:} ~ F(y) \subset T_S(y).
   \end{align}
    Then, $\phi([0,T]) \subset S$.
\end{proposition}

\begin{proof}
   To achieve a contradiction, let $t_1\in [0, T)$ and $t_2 \in (t_1, T]$ such that $\phi([0,t_1]) \subset S \cap C$ and $\phi((t_1,t_2]) \subset C \backslash S$. Given $t \in [t_1,t_2]$, we use $y(t)$ to denote the projection of $\phi(t)$ on $S$ such that $F(y) \subset T_S(y)$, and we let $\delta(t) := |\phi(t)-y(t)|$.
     We recall the difference of squares identity
    \begin{equation}\label{eq_difference_squaes}
        \begin{split}
            (a+b)(a-b) & = a^2 - b^2 \qquad \forall a, b \in \mathbb{R},  \\
            (a+b)^\top (a-b) & = a^\top a - b^\top b \qquad \forall a, b \in \mathbb{R}^n.
        \end{split}    
    \end{equation}
    Using the previous identity, we conclude that, for any $h>0$ and $t \in (t_1, t_2)$ such that $t+h \in (t_1,t_2]$,  
\begin{equation*}
        \begin{aligned}
             \delta(t+h) - \delta(t) 
           &  = \ |\phi(t+h)-y(t+h)| - |\phi(t) - y(t)|\\ & =
            \frac{|\phi(t+h)-y(t+h)|^2 - |\phi(t) - y(t)|^2}{|\phi(t+h) - y(t+h)| + |\phi(t) - y(t)|}.
        \end{aligned}
    \end{equation*}
    Next, for almost all $t \in (t_1, t_2)$,
    we replace $\phi(t+h)$ by 
    \begin{equation}\label{eq_phi_taylor_exp}
        \phi(t+h) = \phi(t) + h \dot{\phi}(t) + o(h). 
    \end{equation} 
   Hence, using~\eqref{eq_phi_taylor_exp} and the inequality
    \begin{equation}
        \begin{split}
            |\phi(t+h) - y(t+h)| & \leq  |\phi(t+h) - y(t)|  \\
            & \leq |\phi(t) - y(t) + h \dot{\phi}(t)| + |o(h)|,
        \end{split}
    \end{equation}
    we obtain that for almost all $t \in (t_1, t_2)$ and $h>0$ such that $t+h \in (t_1,t_2]$,
    \begin{equation}
        \begin{split}
        &\frac{\delta(t+h) - \delta(t)}{h} \leq \\
        &\hspace{2em}\frac{(|\phi(t) - y(t) + h \dot{\phi}(t)| + |o(h)|)^2 - |\phi(t) - y(t)|^2}{h(|\phi(t+h) - y(t+h)| + |\phi(t) +y(t)|)}.
        \end{split}
    \end{equation}
    From the previous inequality, and expanding the numerator using identity~\eqref{eq_difference_squaes}, we conclude that for almost all $t \in (t_1, t_2)$ and $h>0$ and such that $t+h \in (t_1,t_2]$, 
    \begin{equation}
    \label{eq_expansion_delta}
        \begin{aligned}
            &\frac{\delta(t+h) - \delta(t)}{h} \leq\\
            &\frac{(2\phi(t) - 2y(t) + h \dot{\phi}(t))^\top (h\dot{\phi}(t))}{h(|\phi(t+h) - y(t+h)| + |\phi(t) +y(t)|)} + R(h),
        \end{aligned}
    \end{equation}
    \[
        R(h) := \frac{|2\phi(t) - 2y(t) + h \dot{\phi}(t)||o(h)| + |o(h)|^2}{h(|\phi(t+h) - y(t+h)| + |\phi(t) +y(t)|)}.
    \]
    We now argue that for all $t \in (t_1, t_2)$, $R(h) \to 0$ as $h \to 0^+$. Indeed, the term 
    \[
        h R_1(h) := |2\phi(t) - 2y(t) + h \dot{\phi}(t)||o(h)| + |o(h)|^2
    \]
    satisfies $R_1(h) \to 0$ as $h \to 0^+$ since $|o(h)|/h \to 0$ and $|2\phi(t) - 2y(t) + h \dot{\phi}(t)|$ remains bounded. For the term
    \[
        R_2(h) := 1/ \left[|\phi(t+h) - y(t+h)| + |\phi(t) +y(t)|\right],
    \]
    continuity of $\phi$ and of the projection $y$ ensure that $|\phi(t+h)-y(t+h)| \to |\phi(t)-y(t)| = \delta(t)$ as $h \to 0^+$. Since we consider $t$ such that $\delta(t)>0$, the denominator is bounded away from zero for $h$ sufficiently small. Hence $ R(h) \to 0$ as $h \to 0^+$.
    
    Taking the $\limsup$ of~\eqref{eq_expansion_delta} as $h \to 0^+$ and simplifying the first term, we obtain that for almost all $t \in (t_1,t_2)$,
    \begin{equation}\label{eq.ns2}
        \begin{split}
            &\limsup_{h \rightarrow 0^{+}} \frac{\delta(t+h) - \delta(t)}{h} \leq \frac{2(\phi(t)-y(t))^\top \dot{\phi}(t)}{|\phi(t)-y(t)|}. 
        \end{split}
    \end{equation}
    
    Next, we claim that 
    \begin{equation}\label{eqClm} 
       \hspace{-0.3cm} (\phi(t)-y(t))^\top \eta_y \leq 0 \quad  \forall \eta_y \in F(y(t)) \quad \forall t \in (t_1, t_2). 
    \end{equation}
    Using~\eqref{eqClm}, for each $\eta_y \in F(y(t))$, the term 
    \[
        - (2(\phi(t) - y(t))^\top \eta_y)/(|\phi(t)-y(t)|)
    \]
    can be added in~\eqref{eq.ns2}. As a result,  for almost all $t \in (t_1, t_2)$ 
    \begin{equation*}
        \begin{split}
            &\limsup_{h \rightarrow 0^{+}} \frac{\delta(t+h) - \delta(t)}{h} \leq \frac{2(\phi(t)-y(t))^\top (\dot{\phi}(t) - \eta_y)}{|\phi(t)-y(t)|}.
        \end{split}
    \end{equation*}    
    Since $\dot{\phi}(t) \in F(\phi(t))$ for almost all $t \in (t_1, t_2)$, 
    using the one-sided local Lipschitzness property, we conclude that it is always possible to find $\eta^*_y \in F(y(t))$ such that
    \[
        (\phi(t) - y(t))^\top (\dot{\phi}(t) - \eta^*_y) \leq L |\phi(t) - y(t)|^2.
    \]
    Consequently, we obtain that, for almost all $t \in (t_1, t_2)$, 
    \begin{align*} 
        \limsup_{h \rightarrow 0^{+}} \frac{\delta(t+h) - \delta(t)}{h} \leq & 2 L |\phi(t)-y(t)| \leq 2 L \delta(t).
    \end{align*}
    Since $\delta$ is nonnegative and absolutely continuous, and satisfies $\dot{\delta}(t) \leq 2 L \delta(t)$ for almost all $t \in (t_1, t_2)$, Grönwall's inequality~\cite[Theorem 1.12]{taoterenceNonlinearDispersiveEquations2006} implies that
    \[
        \delta(t) \leq \delta(t_1)e^{2 L (t-t_1)} \qquad \forall t \in (t_1, t_2).
    \]
   Hence, $\delta (t) = 0$ for all $t \in (t_1, t_2)$, i.e., $\phi((t_1,t_2)) \subset S$, which contradicts the fact that $\phi([t_1,t_2]) \subset C \backslash S$.
    
    To prove the inequality in~\eqref{eqClm}, we pick $t \in (t_1,t_2)$ and consider, for all $\eta_{y} \in F(y(t))$, the inequality
    \begin{equation}
        \begin{aligned}
            |\phi(t)-y(t)| - |y(t) + h \eta_{y}|_{S} & \leq |\phi(t) - y(t) - h \eta_y|.
        \end{aligned}
    \end{equation}
    To obtain the previous inequality, we used the fact that $|\phi(t)-y(t)| = |\phi(t)|_{S}$ and that  $|\cdot|_{S}$ is globally Lipschitz with a Lipschitz constant equal to $1$. Next, by taking the square in both sides of the last inequality and dividing by $h$, we obtain, for each $\eta_{y} \in F(y(t))$,
    \begin{align*}
        \frac{|\phi(t) - y(t)|^2}{h}  & \leq  \frac{|\phi(t)-y(t) - h \eta_{y}|^2}{h} + 
        \frac{|y(t) + h \eta_{y}|_{S}^2}{h}  \\
        &+ 2 |\phi(t) - y(t) - h \eta_{y}| |y(t) + h \eta_{y}|_{S}/h,
        \end{align*}
    which implies that
    \begin{align*} 
        &\frac{|\phi(t) - y(t)|^2}{h} \leq \frac{|\phi(t)-y(t)|^2}{h} + h |\eta_{y}|^2 -\\
        & 2 (\phi(t)-y(t))^T \eta_{y} + \frac{|y(t) + h \eta_{y}|^2_{S}}{h} +\\
        & 2 |\phi(t)-y(t) - h \eta_{y}| 
        |y(t) + h \eta_{y}|_{S}/h. 
    \end{align*}
    Finally, letting $h \rightarrow 0^{+}$ through a suitable sequence,~\eqref{eqClm} is proved using the fact that $y(t) \in S$  and $\eta_{y} \in F(y) \subset  T_{S}(y(t))$, leading to, see~\cite[Definition 1.1]{aubin1987differential},
    $\liminf_{h \rightarrow 0^{+}} |y(t) + h \eta_{y}|_{S}/h = 0. $
\end{proof}

\begin{remark}
    The proof of Proposition \ref{prop_star_implies_not_leaving} uses the same arguments as the proof of \cite[Theorem 2]{maghenemSufficientConditionsForward2021}, which itself builds upon the proof of 
\cite[Theorem 1]{redheffer1972theorems}. 
In particular, 
Proposition \ref{prop_star_implies_not_leaving} is a local version of \cite[Theorem 2]{maghenemSufficientConditionsForward2021} that relies on the local condition,  \eqref{eqS}, to prevent the solution $\phi$ from exiting $K$ immediately.
Indeed, \cite[Theorem 2]{maghenemSufficientConditionsForward2021} guarantees forward invariance provided that \eqref{eqconstNag}, Assumption \ref{ass_propagation}, and  one of the following hold. 
\begin{itemize}
\item $\forall x \in \left( \mathcal{N}(\mathcal{P}) \cap \partial K \right) \backslash C$, $F(x) \subset T_K(x)$. 

\item $\forall x \in \mathcal{N}(\mathcal{P}) \cap \partial (K \cap C) \cap \partial C$, $F(x) \subset T_{K \cap C}(x)$. 

\item  $C$ is convex and $F(x) \subset T_{K \cap C}(x)$ for all $ x \in \mathcal{P}$.
\end{itemize}
The items above force the projection on $K$ (or on $K \cap C$) of the solution $\phi$, restricted to a neighborhood $\mathcal{N}(\mathcal{P})$, to be within an area where $F$ is within $T_K$ (or $T_{K \cap C}$). 
\end{remark}

\ifitsdraft

\section{Discussion} \label{sec_discussions}
 
We start demonstrating the importance of each Assumption to establish sufficiency of \eqref{eqconstNag}; see Figure~\ref{fig_examples} for illustration. 

\begin{figure*}[t]
    \centering    
    \begin{subfigure}[t]{0.32\textwidth}
        \centering
        \caption{Example~\ref{example_no_continuity}}
        \includegraphics[width=\linewidth]{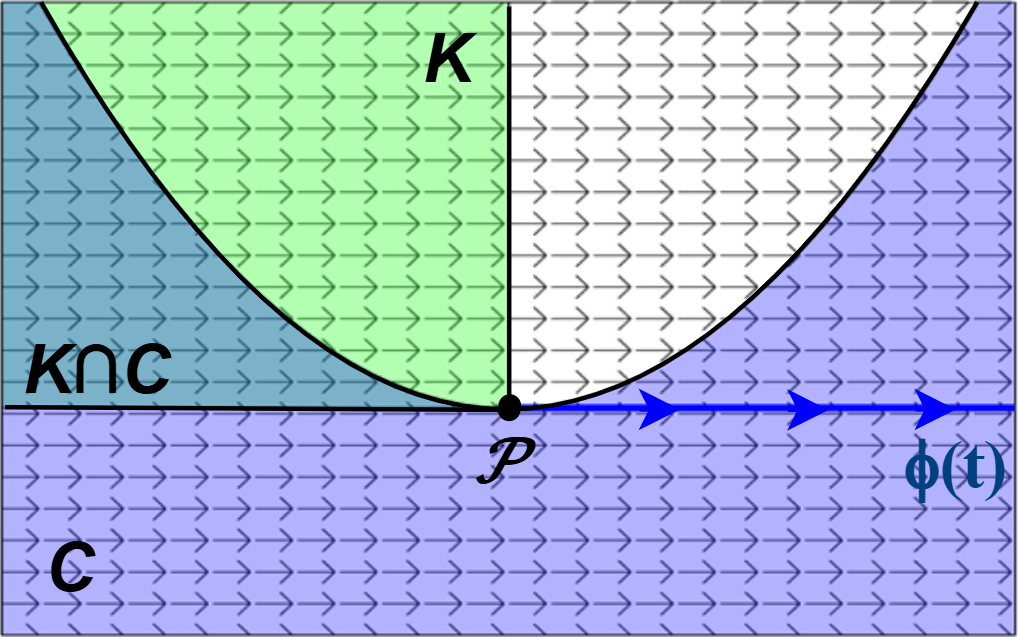}
    \end{subfigure}
    \hfill
    \begin{subfigure}[t]{0.32\textwidth}
        \centering
        \caption{Example~\ref{example_no_nontangenciality}}
        \includegraphics[width=\linewidth]{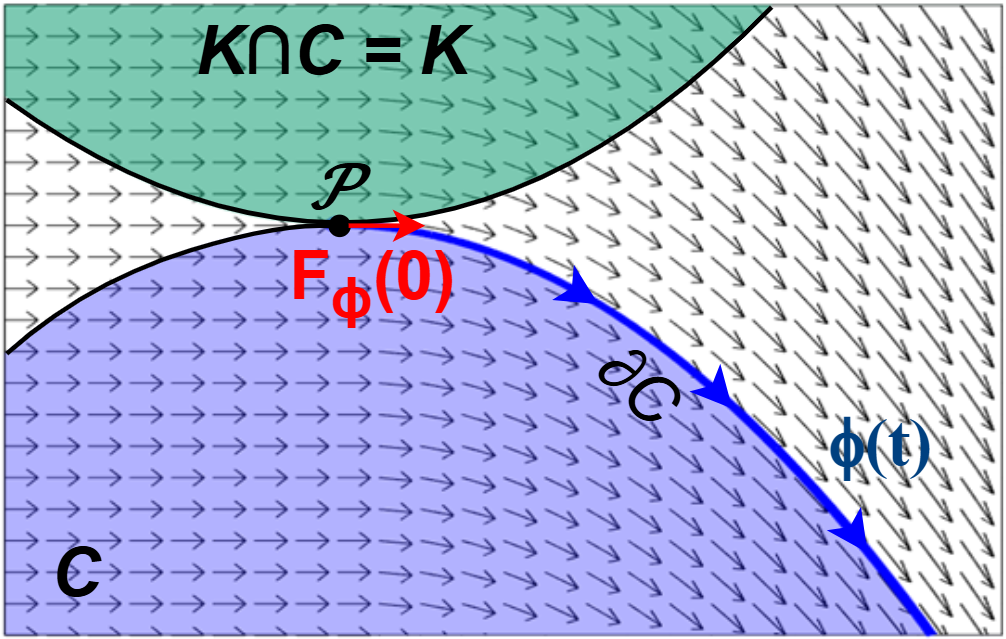}
    \end{subfigure}
    \hfill
    \begin{subfigure}[t]{0.32\textwidth}
        \centering
        \caption{Example~\ref{example_no_transversality}}
        \includegraphics[width=\linewidth]{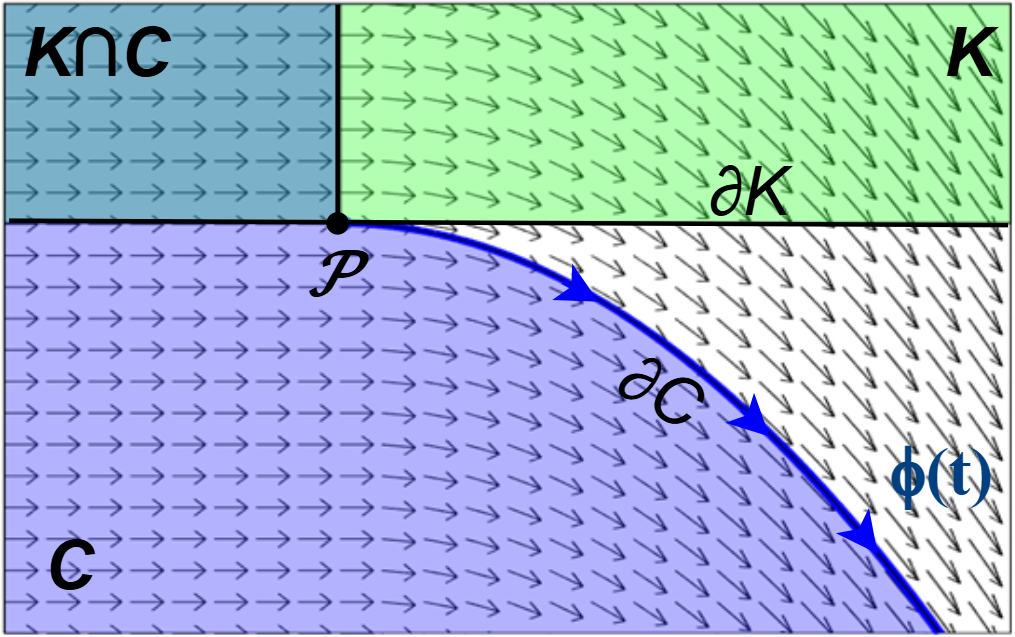}
    \end{subfigure}
    
    \caption{Illustration of the systems in Examples~\ref{example_no_continuity}-\ref{example_no_transversality}. In Example~\ref{example_no_continuity},  $\partial K \ni y \leadsto T_K(y)$ is dis-continuous at $0 \in \mathcal{P}$, $F(0) \not \in T_K(0)$, and the solution $\phi$ exits the set $K$. In Example~\ref{example_no_nontangenciality}, $F(\phi(0))$, 
    $\phi(0)= x \in \mathcal{P}$ is tangent to $\partial K \cap \partial C$. This causes the immediate exit of $K$. In Example~\ref{example_no_transversality}, the transversality assumption is not verified, and a solution starting from~$\mathcal{P}$ exits $K$ immediately.}
    \label{fig_examples}
\end{figure*}

\subsection{Example : Assumption~\ref{ass_propagation} is not verified } \label{example_no_continuity}

Consider the sets $C := \{ x \in \mathbb{R}^2 : x_2 \leq x_1^2 \}$ and $ K := \{ x \in \mathbb{R}^2 : x_1 \leq 0, ~ x_2 \geq 0\}$ and let $F(x) := \{[1  ~~ 0]^\top\}$. 

We can verify that $\mathcal{P} = 
\{0\} = \partial K \cap \partial C$ since $F(0) \not\in T_K(0) = \{v \in \mathbb{R}^2 : v_1 \leq 0, ~ v_2 \geq 0\}$. For the same reason,  Assumption~\ref{ass_propagation} is not satisfied. However, Assumptions \ref{ass_decisiveness}-\ref{ass_transversality} hold. Indeed, 
\[
    F(0) \cap T_{\partial K \cap \partial C}(0) = \{[1 ~~ 0]^\top\} \cap \{0\} = \emptyset,
\]
fulfilling Assumption~\ref{ass_decisiveness}. Moreover, all the items in Assumption~\ref{ass_transversality} hold. Indeed, the first item holds by Lemma~\ref{lem_prac_properties} in the Appendix. The second item holds by sleekness of the set $C$ (i.e., $T_C = T_K^a = T_K$). To verify $(\dagger)$ at $\mathcal{P}$, we note  that 
\begin{align*}
    (\partial K \backslash C) & = \{ x \in \mathbb{R}^2 : x_1 = 0, ~ x_2 > 0\},\\
    (\partial C \backslash \partial K) & = \{ x \in \mathbb{R}^2 : x_2 = x_1^2, ~ x_2 \neq 0\}.
\end{align*}
Now, for any neighborhood $\mathcal{N}(0)$, we let 
\begin{align*}
    z & := [0 ~~ z_{2}]^\top \in ((\partial K \backslash C) \cap \mathcal{N}(0)), \\
    y & := [y_{1} ~~ y_{1}^2]^\top \in ((\partial C \backslash \partial K) \cap \mathcal{N}(0)), 
\end{align*}
so that  
\begin{align*}
T_{\partial K}(z)  =  
\mathbb{R} \times \{0\}, 
~~
T_{C}(y) =  \{v \in \mathbb{R}^2 : v_2 \leq 2 y_{1} v_1 \}.
\end{align*}

Furthermore, for 
\begin{align*}
w_z  & := [0 ~~~ - y_{1}^2- z_2]^{\top} \in T_{\partial K}(z),\\ 
w_y & := [-y_{1} ~~ -2y_{1}^2]^{\top} \in T_{C}(y), 
\end{align*}
it holds that $y - z = w_{z} - w_{y}$. Moreover,
\[
    \begin{split}
      |[w_{z}^\top ~ w_{y}^\top]|^2 & = (y_{1}^2 + z_{2})^2 + y_{1}^2 + (2 y_{1}^2)^2 \\ & 
      = |y-z|^2 + 4 (y_{1}^2 + y_{1}^2 z_{2}).
    \end{split}
\]
Therefore, 
$|[w_z^\top ~ w_y^\top]| \leq \gamma |y - z|$  if  $\gamma^2 \geq 1 + \frac{4 (y_{1}^2 + y_{1}^2 z_{2})}{|y-z|^2}. $ 
Besides, for any $k > 0$, we have $ 4 (y_{1}^2 + y_{1}^2 z_{2})/(|y-z|^2) < k $ whenever $z_{2}$ and $y_{1}$ are small enough. This is enough to conclude that ($\dagger$) holds at $x = 0$, for $\alpha = 0$ and $c = \sqrt{1 + k}$, by choosing a small enough $\mathcal{N}(0)$. 

At the same time, \eqref{eqconstNag} holds since for all $x \in \partial K \cap \text{int}(C)$,
\[
    F(x) = \{[1 ~~  0]^\top\} \subset T_{K}(x) = \{v \in \mathbb{R}^2 : v_2 \geq 0\}.
\]
However, the solution $\phi(t) := [t ~~ 0]^\top$ leaves the set $K$.

\subsection{Example : Assumption~\ref{ass_decisiveness} is not verified}\label{example_no_nontangenciality}
    Consider the sets 
    \begin{equation*}
        C :=  \{x \in \mathbb{R}^2 : x_2 \geq x_1^2/2\} \cup \{x \in \mathbb{R}^2 : x_2 \leq - x_1^2/2\}\\
    \end{equation*}
    and $K := \{x \in C : x_2 \geq 0\}$, and let $F(x) := \{[1  ~~ -x_1]^\top\}$. We can see that $\mathcal{P} = \{0\} = \partial K \cap \partial C$ since the solution $\phi(t) := [t ~~  -t^2/2]^\top$ leaves the set $K$ immediately.
    
    We further note that
    \begin{align}
        &T_{K}(0) = \{v \in \mathbb{R}^2 : v_2 \geq 0\},\label{eq_TK_ex_tang}\\
        &T_{C}(0) = T_{C}^a(0) =  \mathbb{R}^2, \label{eq_TpparC_ex_tang}
    \end{align}
    \begin{equation}
    \label{eq_TparKcapparC_ex_tang}
        \begin{split}
            T_{\partial K}(0) & = T_{\partial K}^a(0) = T_{\partial K \cap \partial C}(0) =  \mathbb{R} \times \{0\}.
      \end{split}
    \end{equation}
    Since $F(0) = \{[1  ~~ 0]^\top\} \subset T_{\partial K \cap \partial C}(0)$, we conclude that  Assumption~\ref{ass_decisiveness} does not hold. 
    However, Assumption~\ref{ass_propagation} holds, under~\eqref{eq_TK_ex_tang}. Furthermore, all the items in Assumption~\ref{ass_transversality} hold as well. Indeed, The first two items follow directly form~\eqref{eq_TpparC_ex_tang} and~\eqref{eq_TparKcapparC_ex_tang}, and $(\dagger)$ holds trivially at $\mathcal{P}$ because $\partial K \backslash C = \emptyset$.  
    At the same time,~\eqref{eqconstNag} also holds trivially because $\partial K \cap \mathrm{int}(C) = \emptyset$.   

\subsection{Example: Assumption~\ref{ass_transversality} is not verified}\label{example_no_transversality}

    Consider the sets 
    \begin{equation*}
            C := ( \mathbb{R}_{\leq 0} \times \mathbb{R} ) \cup \{x \in \mathbb{R}^2 : x_2 \leq - x_1^2/2\}
    \end{equation*}
    and $K := \mathbb{R} \times \mathbb{R}_{\geq 0}$, and let
    \begin{equation*}
        F(x) := \begin{cases}
                \{[1 ~~ 0]^\top\} \quad \text{if } x_1 \leq 0, \\
                \{[1 ~~ -x_1]^\top\} \quad \text{otherwise}.
            \end{cases}
    \end{equation*}
    Note that $\mathcal{P} = \{0\} = \partial K \cap \partial C$ since the solution $\phi(t) := [t  ~~  -t^2/2]^\top$ leaves the set $K$ immediately, rendering $K$ not forward invariant. 

    That said, Assumption~\ref{ass_propagation} holds by 
    Lemma~\ref{lemAp} in the Appendix.  Moreover, Assumption~\ref{ass_decisiveness} holds because
    \[
        F(0) \cap T_{\partial K \cap \partial C}(0) = \{[1 ~~ 0]^\top\} \cap \{0\} = \emptyset.
    \] 
    
    However, we can show that the third item in Assumption~\ref{ass_transversality} is not verified. To do so, we pick any neighborhood $\mathcal{N}(0)$,  $c > 0$, and $\alpha \in [0,1)$, and show the existence of $z \in (\partial K \backslash C) \cap \mathcal{N}(0)$ and $y \in (\partial C \backslash \partial K) \cap \mathcal{N}(0)$ such that, for any $(v, w) \in T_{\partial K}(z) \times T_{C}(y)$,  
    \begin{align} \label{eqy-z}
    y-z \in v - w + \alpha (|z - y|)\mathbb{B} 
   \end{align}
   implies 
   \begin{align*}
   |[v^\top ~ w^\top]| > c |z - y|. 
    \end{align*}
    First, $z \in (\partial K \backslash C) \cap \mathcal{N}(0)$ and $y \in (\partial C \backslash \partial K) \cap \mathcal{N}(0)$ imply that, for some $\epsilon > 0$, we have  
    $$ z = [z_{1} ~~ 0]^{\top}, ~ y = [y_{1} ~~ -y_{1}^2/2]^{\top}   \text{ for some } z_{1}, y_{1} \in (0, \epsilon). $$ 
    Furthermore, we note that 
    \[
        \begin{split}
            T_{\partial K}(z) & =  \{v \in \mathbb{R}^2 : v_2 = 0\} \qquad \forall z \in \partial K \backslash C,\\
            T_{C}(y) & =  \{v \in \mathbb{R}^2 : v_2 \leq - y_{1} v_1 \} \qquad \forall y \in \partial C \backslash \partial K.
        \end{split}
    \]
    Next, we choose $z$ and $y$ such that the first components verify $y_{1} - z_{1} = 0$, which implies that $y - z =  \frac{1}{2} [0 ~~ -y_{1}^2]^{\top}$. 
    
    Having \eqref{eqy-z}, we conclude   that, for some $a_1, a_2 \in (-1,1)$, 
    \[
    y - z = \begin{bmatrix}
      0 \\ -y_{1}^2/2
    \end{bmatrix} = \begin{bmatrix}
      v_{1} - w_{1} \\ -v_{2}
    \end{bmatrix} + \frac{y_{1}^2}{2} \begin{bmatrix}
      a_1  \\ a_2  
    \end{bmatrix}. 
    \]
    This is possible only when 
    \begin{align*}
      v - w & = a_1 y_{1}^2/2, \quad w  = (1 + a_2) y_{1}^2/2, \\
      w & \leq - v_{2}/y_{1} = -(1 + a_2) y_{1}/2 < 0.
    \end{align*} 
    Hence, 
    \[
    \begin{split}
    |[v^\top ~ w^\top]|^2 & \geq |w|^2 \geq (1 + a_2)^2 \left(y_{1}^4 + y_{1}^2\right) /4 \\
    & = (1 + a_2)^2 
    \left(1 + 1/y_{1}^2\right) y_{1}^4/4
    \end{split}
    \]
    and, at the same time, $c^2 |z - y|^2 = c^2 z_{1}^4/4$.
    By choosing $y_{1}$ sufficiently small such that $(1 + a_2)^2 (1 + 1/y_{1}^2) > c^2$, it follows that 
    $|[v^\top ~ w^\top]|^2 > c^2 |z - y|^2$. 
    
    It is easy to see that~\eqref{eqconstNag} also holds since, for all $x \in \partial K \cap \mathrm{int}(C)$,
    $F(x) = \{[1 ~~ 0]^\top\} \subset T_K(x) = \{v \in \mathbb{R}^2 : v_2 \leq 0\}$. 

\newpage 
    
\subsection{Alternatives to ($\dagger$) in Assumption~\ref{ass_transversality}}
\label{subsec_alternatives}

 According to Lemma~\ref{lemtranssimple} in the Appendix, ($\dagger$) in Assumption~\ref{ass_transversality} is verified for a given $x \in \partial K \cap \partial C$ if 
\begin{itemize}
    \item[($\ddagger$)] $\exists \mathcal{N}(x)$, a neighborhood of $x$, such that 
    \begin{equation*} 
        \begin{aligned}
             x_2 - x_1 & \in 
            C_{\partial K} (x) - C_{C} (x) 
            \\ &
            := \{ v - w : (v, w) \in C_{\partial K} (x) \times C_{C} (x) \}  
          \\  \forall (x_1,x_2) & \in 
                ((\partial K \backslash C) \cap \mathcal{N}(x)) \times  ((\partial C \backslash \partial K) \cap \mathcal{N}(x)).
            \end{aligned}
    \end{equation*}
\end{itemize}

Although ($\ddagger$) is in general more restrictive than ($\dagger$), see Example~\ref{example_simpler_trans} below, it is used later to verify ($\dagger$) using a simpler and more tractable assumption when the sets $(C,K)$ are finitely represented.

\begin{example} \label{example_simpler_trans}
    Let
    \begin{align*} 
        K = \partial K & := \{ x \in \mathbb{R}^3 :x_2 \leq \sqrt{|x_1|}, ~ x_3 = 0\}, \\
        C = \partial C & :=  \{ x \in \mathbb{R}^3 : x_2 \geq  \sqrt{|x_1|}, ~ x_3 = 0\}.
    \end{align*}
    At $0 \in  \partial K \cap \partial C = K \cap C$, we can see that 
    \begin{align*}
        C_{C}(0) & = \{x \in \mathbb{R}^3 : x_1 = x_3 = 0, ~ x_2 \geq 0 \},\\
        C_{\partial K}(0) & = \{x \in \mathbb{R}^3 :  x_3 = 0, ~ x_2 \leq 0 \}. 
    \end{align*}
    Hence, $C_{\partial K}(0) - C_{C}(0) = \{ x \in \mathbb{R}^3 : x_3 = 0, ~ x_2 \leq 0 \}$. 

    At the same time, for any $\mathcal{N}(0)$, we will find $(x,z) \in (\partial K \backslash C) \cap \mathcal{N}(0) \times  (\partial C \backslash \partial K) \cap \mathcal{N}(0)$ such that 
    $$ x - z \notin C_{\partial K}(0) - C_{C}(0).  $$
    Indeed, let $z:= [0  ~~  \epsilon  ~~  0]^\top \in \partial C \backslash \partial K$, for $\epsilon > 0$, and let 
    $x := [5 \epsilon^2   ~~  2 \epsilon  ~~  0]^\top \in \partial K \backslash C$.  Now, for $\epsilon > 0$ is sufficiently small, it follows that
    \[
        x - z = [5\epsilon^2  ~~  \epsilon   ~~  0]^\top  \notin C_{\partial K}(0) - C_{C}(0).
    \]
    As a result, for any neighborhood $\mathcal{N}(0)$,
    \[
        [(\partial K \backslash C )\cap\mathcal{N}(0)] - [(\partial C \backslash \partial K) \cap\mathcal{N}(0)] \not\subset C_{\partial K}(0) - C_{C}(0),
    \]
    and $(\ddagger)$ does not hold. However, ($\dagger$) trivially holds  at $x = 0$, for $\alpha = 0$, $c = 1$, $v_2 = 0$, and $v_1 = x_2 -x_1$, since $T_{\partial K}(0) = \{ x \in \mathbb{R}^3 : x_3 = 0\}. $
\end{example}

\begin{remark}
We succeeded to show that, when ($\dagger$) in Assumption~\ref{ass_transversality} is verified at some $x \in \partial K \cap \partial C$ and $T_C(x) = T^a_C(x)$, then $T_{\partial K \cap C}(x) = T_{\partial K}(x) \cap T_{C}(x)$.  The latter equality is established in existing literature, under one of the following more-restrictive conditions.
\begin{enumerate}
    \item  The sets $\partial K$ and $C$ are sleek and  $T_{\partial K}(x) \cap \text{int} (T_{C}(x)) \neq \emptyset$; see~\cite{rockafellar1979clarke}.
    \item   There exists  $v_o \in C_{\partial K}(x)$ such that 
    $$ \alpha v_o + (1- \alpha) v \in D_{C}(x) \qquad \forall v \in T_{C}(x), \quad \forall \alpha \in (0,1];  $$ 
      see~\cite[Proposition 4.3.7]{aubinSetValuedAnalysis2009}.      
    \item   $T_C(x) = T^a_C(x)$ and  
    $ C_{\partial K}(x) - C_{C}(x) = \mathbb{R}^n$; see~\cite[Corollary 4.3.5]{aubinSetValuedAnalysis2009}. 
    \item   
   The sets  $(\partial K,C)$ are sleek and  
    $ N^L_{\partial K}(x) \cap N^L_{C}(x) = \{0\}; $
    see~\cite[Theorem 6.42]{rockafellarVariationalAnalysis2009},
    where $N^L_{K}(x) := \limsup_{y \rightarrow x} N^P_{K}(y)$ and  $N^P_{K}$ is the proximal normal cone. 
\end{enumerate}
The non-sleek sets $K=C=\{x\in\mathbb{R}^2: x_2\le |x_1|\}$ makes the above conditions unsatisfied at $x=0$, although conditions $(\ddagger)$ and $(\dagger)$ are satisfied.
\end{remark}

\fi

\section{Result for Finitely-Represented Sets} \label{sec_prac} 

We propose to restrict ourselves to finitely-represented sets $C$ and $K$.

\begin{definition} \label{defPrac}
     A set $\star \subset \mathbb{R}^n$ is a finitely-represented set if there exists $d_\star \in \mathbb{N}$ and continuously-differentiable functions $g_{\star}^i: \mathbb{R}^n \to \mathbb{R}$ for all $i \in [d_{\star}]$ such that
    \begin{equation} \label{eq_def_ke_ineqs}
        \star := \{x \in \mathbb{R}^n : g_{\star}^i(x) \leq 0 ~~~ \forall i \in [d_{\star}]\},
    \end{equation}
    where $[d_\star] := \{1, 2, \ldots, d_\star\}$.
\end{definition}
We further require the following assumption.
\begin{assumption} [Constraint qualification] \label{ass_cq_practical}
  The sets $K$ and $C$ are 
  finitely-represented sets and, 
  for all $\star \in \{C, K\}$, 
  it holds that 
    \begin{equation}\label{eq_cq_pra}
     \hspace{-0.2cm} \forall x\in \partial \star, ~ \exists v_o \in \mathbb{R}^n : \nabla g_{\star}^i(x)^{\top} v_o < 0 ~~~  \forall i \in \mathrm{Act}_{\star}(x),
    \end{equation}
    where $  \mathrm{Act}_{\star} (x)$, the set of active equality constraints at $x$, is given by 
    \[
        \mathrm{Act}_{\star}(x) :=  \begin{cases}
            \{i \in  [d_{\star}] : g_{\star}^{i}(x)= 0\} & \text{if } x \in \star, \\
            \emptyset & \text{if } x \notin \star.
        \end{cases}    
    \]
\end{assumption}

For a finitely-represented set $\star$, condition \eqref{eq_cq_pra} coincides with the 
Mangasarian-Fromovitz Constraint Qualification (MFCQ) used in~\cite{mangasarianFritzJohnNecessary1967}. It is named transversality in \cite[Proposition 4.3.7]{aubinSetValuedAnalysis2009} and~\cite[Assumption 1]{maghenemSufficientConditionsForward2021}, and allows the representation
\begin{equation}\label{eq_TK_practical}
   \hspace{-0.3cm} T_{\star}(x) = \left\{v \in \mathbb{R}^n ~:~ \nabla g_i (x) v \leq 0\quad \forall i \in  \mathrm{Act}_{\star}(x)\right\}.
\end{equation}
Finally, a finitely-represented set $\star$ verifying \eqref{eq_cq_pra} is also termed practical set in~\cite[Definition 4.9]{blanchini2008set}. 

We show later that, for finitely-represented sets $(K,C)$ verifying Assumption \ref{ass_cq_practical}, Assumption \ref{ass_decisiveness} (decisiveness) is verified under the following assumption.

\begin{assumption}
\label{ass_decisiveness_pra}
    {(\bfseries\upshape Decisiveness under finite representation})
    
        For each $x \in \mathcal{P}$ and for each $\eta \in F(x)$ such that
    \begin{equation}
    \label{eq_decis_prac1}
            \nabla g_{\star}^i(x)^\top \eta \leq 0 \quad \forall i \in  \mathrm{Act}_{\star}(x), \quad \forall \star \in \{K,C\},
        \end{equation}
        it holds that
$$ \eta^\top (\nabla g_{K}^i(x) + 
\nabla g_{C}^k(x))  < 0 ~~ 
 \forall (i,k) \in  \mathrm{Act}_{K}(x) \times \mathrm{Act}_{C}(x). $$
\end{assumption}

Furthermore, according to Lemma~\ref{lem_prac_properties} in the Appendix, the first item in Assumption~\ref{ass_transversality} holds without any additional assumption. However, for the remaining two items to hold, we use the following additional assumption. 
\begin{assumption}\label{ass_transversalityPra}
    {(\bfseries\upshape Transversality under finite representation)}
    For each $x \in \mathcal{P} \subset \partial K \cap \partial C$ and for each $j \in \mathrm{Act}_{K}(x)$, there exists $v_j \in \mathbb{R}^n$ such that 
    \begin{align} 
        \nabla g_{K}^j(x)^\top v_j & = 0, \label{eq1needed}\\ 
        \nabla g_{K}^i(x)^\top v_j & < 0 \quad \forall i \in \mathrm{Act}_{K}(x) \backslash\{j\},  \label{eq2needed}\\ 
        \nabla g_{C}^k(x)^\top v_j & < 0 \quad \forall k \in \mathrm{Act}_{C}(x). \label{eq3needed}
    \end{align}
\end{assumption}

Assumption \ref{ass_transversalityPra} guarantees, for a given $x \in \mathcal{P}$, the existence
of a vector $v_j$, for all $j \in \text{Act}_K(x)$, that is tangent
to the active boundary $j$ via \eqref{eq1needed}. Simultaneously,
$v_j$ points strictly inward of all the other active constraints
of $K$ and $C$ at $x$ by, respectively, satisfying \eqref{eq2needed}
and \eqref{eq3needed}. 
This forces a non-tangent separation between $\partial K$ and $\partial C$. 
Indeed, we are able to show that   condition~($\dagger$) in Assumption~\ref{ass_transversality} follows from Assumption \ref{ass_transversalityPra}; see the forthcoming Proposition \ref{prop_transpractical}. Moreover, we show that  conditions~\eqref{eq1needed} and~\eqref{eq2needed} imply $T_{\partial K}(x)=T^a_{\partial K}(x)$; see Lemma \ref{lem_prac_properties} in the Appendix.

\begin{remark}[Polytopic sets]
    Assumption~\ref{ass_transversalityPra} is more restrictive than Assumption~\ref{ass_transversality}. 
    Indeed, as we show in Lemma~\ref{lem_trans_poly} in the Appendix, Assumption~\ref{ass_transversality} holds for free when $K$ and $C$ are \textit{polytopic}, i.e.,
    \begin{equation}
    \label{eqPolitopKC}
    \begin{aligned} 
        K & := \{x \in \mathbb{R}^n  :  a^i_{K} x - b^i_{K} \leq 0 ~~~ \forall i \in [d_{K}]\}, \\
        C & := \{x \in \mathbb{R}^n  :  a^i_{C} x - b^i_{C} \leq 0 ~~~ \forall i \in [d_{C}]\}. 
   \end{aligned}
   \end{equation}
   This is clearly not the case of Assumption \ref{ass_transversalityPra}.
\end{remark}

\ifitsdraft

\begin{example}
    Another example of sets fulfilling Assumption~\ref{ass_transversality} but not~Assumption \ref{ass_transversalityPra} are those given in Example~\ref{example_assumptions1}.  Indeed, $K$ and $C$ therein are finitely-represented sets with $g_K^1(x) = x_1 - 1$, $g_C^1(x) = -x_1^2 - x_2^2 +1$, $g_C^2(x) = -x_1$ and $g_C^3(x) = x_2$. Since   $2 \nabla g_K^1(x) = \nabla g_C^1(x) = [-2  ~~  0  ~~ 0]^\top$  for all $x\in \mathcal{P}$,
     we can see that Assumption~\ref{ass_transversalityPra} does not hold because $\nabla g_K^1$ and $\nabla g_C^1$ are collinear on $\mathcal{P}$.
\end{example}

\else
Another example of sets fulfilling Assumption~\ref{ass_transversality} but not Assumption~\ref{ass_transversalityPra} is given in \cite[Example 5]{arXiv_version_auto}.
\fi

We are now ready to provide a  characterization of forward invariance in the case of finitely-represented sets.  

\begin{theorem} \label{thm2_practical}
    Let $K$ and $C$ be finitely-represented sets fulfilling Assumption~\ref{ass_cq_practical}.
    \begin{itemize}
    \item When $F$ in $\Sigma$ is lower semicontinuous, forward invariance of $K$ implies 
    \begin{equation} 
    \label{eqMotivprac}
        \begin{aligned} 
      \hspace{-0.4cm}  \nabla g^j_K(x)^\top \eta  \leq 0 \quad &  \forall j \in \text{Act}_K(x), \quad 
        \forall \eta \in F(x),
        \\ & \hspace{-0.2cm} \forall x \in \partial K \cap \text{int}(C).  
        \end{aligned}
    \end{equation}
    \item The opposite is true provided that $F$ is one-sided locally Lipschitz and one of the following holds.  \\
   -Assumptions~\ref{ass_decisiveness_pra} and~\ref{ass_transversalityPra}, and  $F$ lower semicontinuous. 
   \\
   -Assumptions~\ref{ass_propagation} and~\ref{ass_decisiveness_pra} and $(K,C)$ being polytopic.\\
   - Assumption~\ref{ass_transversalityPra}, $F$ lower semicontinuous, and density of $\partial K \cap \text{int}(C)$  in $\partial (K \cap C)$. 
   \end{itemize}
\end{theorem}

\begin{proof}
The necessity of~\eqref{eqMotivprac} follows from Theorem~\ref{theo_main_result} and the first item of Lemma~\ref{lem_prac_properties} in the Appendix.

To prove sufficiency of~\eqref{eqMotivprac}, we use Lemma~\ref{lem_transPrac_implies_eqTangCond} in the Appendix to obtain from Assumption~\ref{ass_transversalityPra}, continuity of $F$, and \eqref{eqMotivprac} that the propagation assumption  (Assumption~\ref{ass_propagation}) holds, i.e.,
\begin{align} \label{eqtouse}
F(x) \subset T_{K}(x) \qquad \forall x \in \mathcal{P}. 
\end{align}

Now, we verify Assumptions~\ref{ass_decisiveness} and~\ref{ass_transversality} under either one of the first two items in Theorem~\ref{thm2_practical}, which is enough to conclude sufficiency of~\eqref{eqMotivprac} under these two items, thanks to Theorem \ref{theo_main_result}. Indeed, using Proposition~\ref{prop_nontang_prac}, we conclude that the decisiveness assumption  (Assumption~\ref{ass_decisiveness}) is verified under Assumption~\ref{ass_decisiveness_pra}. Furthermore, Proposition~\ref{prop_transpractical} can be used to conclude that the transversality assumption 
(Assumption~\ref{ass_transversality}) is verified when Assumption~\ref{ass_transversalityPra} holds. Finally, Lemma~\ref{lem_trans_poly} in the Appendix can be used to conclude that Assumption~\ref{ass_transversality} holds naturally when $(K, C)$ are polytopic. This completes the proof under the first two items.

Now, under the third item in Theorem~\ref{thm2_practical}, we show that the decisiveness assumption (Assumption \ref{ass_decisiveness_pra}) is not needed and thus the already established propagation (Assumptions \ref{ass_propagation}) and the transversality assumption (Assumption \ref{ass_transversalityPra}) are enough. Indeed, we note that
$$ \partial (K \cap C) = (\partial K \cap \text{int}(C)) \cup (\partial K \cap \partial C) \cup (\text{int}(K) \cap \partial C). $$

Furthermore, by the density property, we conclude that $\text{int}(K) \cap \partial C = \emptyset$ and that $\mathcal{P} = \partial K \cap \partial C$. Hence,~\eqref{eqtouse} can be expressed as $F(x) \subset T_{K}(x)$ for all $x \in (\partial K \cap \partial C)$. In addition to~\eqref{eqconstNag}, we conclude that 
\begin{align*}
    F(x) & \subset T_{K}(x)
    \qquad  \forall x \in (\partial K \cap \partial C) \cup (\partial K \cap \text{int}(C)). 
\end{align*}
Note that 
$ T_{K}(x) = T_{K \cap C}(x)$  for all $x \in \partial K \cap \text{int}(C). $ 
Hence, if we show that 
\begin{align} \label{eqtoshw}
    F(x)  \subset T_{K \cap C}(x)  \qquad \forall x \in \partial K \cap \partial C. 
\end{align}
Then, we would conclude that 
$F(x) \subset T_{K \cap C}(x)$ for all $x \in \partial (K \cap C)$,
which proves forward invariance of the set $K$. To verify~\eqref{eqtoshw}, it is enough to show that 
\begin{align} \label{eqtoshw1}
    F(x) \subset T_{C}(x)  \qquad \forall x \in \partial K \cap \partial C. 
\end{align}
Indeed, under Assumption~\ref{ass_transversalityPra}, we have shown in Proposition~\ref{prop_transpractical} that ($\dagger$) in Assumption~\ref{ass_transversality} is verified. At the same time, ($\dagger$) entails the existence of a neighborhood of $x$, $\mathcal{N}(x)$, $c> 0$, and $\alpha \in [0,1)$ such that
\begin{align*}
     x_2 - x_1 & \in v_1 - v_2 + \alpha (|x_1 - x_2|) \mathbb{B} 
        \\ 
        \forall (x_1,x_2) & \in 
             ((\partial K \backslash C) \cap \mathcal{N}(x)) \times
            ((\partial C \backslash K) \cap \mathcal{N}(x))
        \end{align*}
for some $(v_1,v_2) \in  T_{K} (x_1) \times T_{C} (x_2)$ such that $|\left[v_1^\top ~~ v_2^\top\right]| \leq c |x_1 - x_2|$.
Hence, according to Lemma~\ref{lemtrans} in the Appendix, we conclude that 
\begin{align*} 
 T_{C \cap K}(x) = T_C(x) \cap T_K(x) \qquad \forall x \in \partial K \cap 
 \partial C. 
\end{align*}

So, let us now verify~\eqref{eqtoshw1} using contradiction. That is, we pick $x \in \partial K \cap 
\partial C$ and we assume the existence of $\eta \in F(x)$ such that $\eta \notin T_C(x)$. As a result, there exists $L_x \subset  \text{Act}_C(x)$ such that 
$\nabla g^i_C(x)^\top \eta  > 0 ~~ \forall i \in L_x$ and 
$\nabla g^i_C(x)^\top \eta  \leq  0 ~~ \forall i \in \text{Act}_C(x) \backslash L_x$. 
Next, we introduce the vector function  

$$ v(\lambda) :=  \eta + \lambda  \sum_{j \in \text{Act}_K(x)}  v_j,  \quad \lambda > 0, $$
where $\{ v_j\}_{j \in \text{Act}_K(x)}$ is introduced in Assumption~\ref{ass_transversalityPra}. 

Using continuity of $v$, we conclude the existence of $\lambda^* > 0$ such that, for $v^* := v(\lambda^*)$, we have 
\begin{equation}\label{eqalmst}
\begin{aligned}
   \nabla g^i_C(x)^\top v^* & \leq 0 \qquad \forall i \in  \text{Act}_C(x),
   \\
     \nabla g^l_C(x)^\top v^* & = 0 \qquad \text{for some}  ~ l \in  L_x.
\end{aligned}
\end{equation}
 
Additionally, since $\eta \in T_K(x)$, we conclude that
\begin{align} \label{eqalmst1}
\nabla g^j_K(x) v^* & < 0 \qquad \forall j \in  \text{Act}_K(x).
\end{align}
Finally, we show that a consequence of~\eqref{eqalmst}-\eqref{eqalmst1} is that $\partial C \cap \text{int}(K) \neq \emptyset$, which yields to a contradiction. Indeed, by Lemma~\ref{lem_prac_properties} in the Appendix 
\eqref{eqalmst1} implies that $v^* \in D_K(x)$ and \eqref{eqalmst} implies 
that $v^* \in T_{\partial C}(x)$. Hence, by definition of $T_{\partial C}$,  we conclude the existence of $\{t_i\}^\infty_{i = 1} \subset \mathbb{R}_{>0}$,  $\{v_i \}^\infty_{i=1} \subset \mathbb{R}^n$ such that  $t_i \rightarrow 0^+$,  $v_i \rightarrow v^*$, and $y_i := x + t_i v_i \in \partial C$ for all $i \in \{1,2,...\}$.   At the same time,  by definition of Dubovitskiy cone $D_K$,  we conclude that $y_i \in \text{int}(K) \cap \partial C$ for all $i$ sufficiently large.  
\end{proof}

   \ifitsdraft

Figure~\ref{fig_implications} provides a schematic view of the implications presented throughout this section and Section~\ref{sec_main_result}.
    
    \begin{figure*}[t]
        \makebox[\textwidth][c]{%
            \hspace*{-2.5cm}\input{implications}
        }
        \caption{Summary of the results presented throughout this article.}
        \label{fig_implications}
    \end{figure*}
\else

A  schematic figure summarizing  the implications presented so far is in \cite{arXiv_version_auto}.

\fi

\subsection{Intermediate Results}
\label{subsec_intermediate_results_prac}

\ifitsdraft

For finitely-represented sets $K, C \subset \mathbb{R}^n$ fulfilling Assumption~\ref{ass_cq_practical}, 
the following proposition establishes that Assumption~\ref{ass_decisiveness_pra} implies Assumption~\ref{ass_decisiveness}. 
The converse is established in Lemma~\ref{lem_nontang_prac} in the Appendix, under an additional condition.

\else

    The following proposition establishes that for 
    finitely-represented sets $K, C \subset \mathbb{R}^n$ fulfilling Assumption~\ref{ass_cq_practical},  Assumption~\ref{ass_decisiveness_pra} implies Assumption~\ref{ass_decisiveness}. The converse is established in~\cite[Lemma 23]{arXiv_version_auto}, under an extra condition.
\fi

\begin{proposition} \label{prop_nontang_prac}
    For finitely-represented sets $K, C \subset \mathbb{R}^n$ fulfilling Assumption~\ref{ass_cq_practical}, Assumption~\ref{ass_decisiveness_pra} $\Rightarrow$ Assumption~\ref{ass_decisiveness}.
\end{proposition}

\begin{proof}
 The proof is based on characterizing the maps $\partial K \ni x \leadsto T_{\partial K}(x)$  and  $\partial C \ni  x \leadsto T_{\partial C}(x)$ in terms of  $\{ g^i_K \}^{d_K}_{i = 1}$ and $\{ g^k_C \}^{d_C}_{k = 1}$, respectively.  
Indeed, we already know that for all $x \in \mathcal{P}$,
\begin{equation*}
    \begin{aligned} 
       &  T_{\partial K \cap \partial C}(x)   \subset T_{\partial K}(x) \cap T_{\partial C}(x) 
        =
        \\ &
        \left\{ v \in T_K(x) \cap T_C(x) : 
        \exists (i,k) \in  \mathrm{Act}_{K}(x) \times \mathrm{Act}_{C}(x) 
        : \right. 
        \\ &
        \left. 
        \left( \nabla g_{K}^i(x) + \nabla g_{C}^k(x) \right)^\top v =0   \right\}.
    \end{aligned}
\end{equation*}
With Assumption~\ref{ass_decisiveness_pra}, we guarantee Assumption~\ref{ass_decisiveness}.
\end{proof}

\begin{proposition}  \label{prop_transpractical}
   Under the hypothesis of Proposition \ref{prop_nontang_prac},  Assumption~\ref{ass_transversalityPra} $\Rightarrow$  Assumption~\ref{ass_transversality}. 
\end{proposition}

\begin{proof}  
    Given $x \in \mathcal{P}$,  the first and the second items 
    in Assumption~\ref{ass_transversality} are verified using the second and the sixth items in Lemma~\ref{lem_prac_properties} in the Appendix. 

    \ifitsdraft
Next, to verify ($\dagger$) in Assumption~\ref{ass_transversality},  we propose to verify ($\ddagger$) in Section~\ref{subsec_alternatives} for the sets $(\partial K_j,C)$, instead of the sets $(\partial K, C)$, for all $j \in \text{Act}_K(x)$, where $K_j$ was defined in~\eqref{eq_def_Kj}. 
    \else
        Next,  to verify ($\dagger$) in Assumption~\ref{ass_transversality},  we propose to verify ($\ddagger$) for the sets $(\partial K_j,C)$, instead of the sets $(\partial K, C)$, for all $j \in \text{Act}_K(x)$, where $K_j$ was defined in~\eqref{eq_def_Kj}.
    \fi
    Then,  we  apply Lemma~\ref{lemtranssimple} in the Appendix to conclude that  ($\dagger$) holds for $(\partial K_j,C)$, instead of $(\partial K, C)$,  for all $j \in \text{Act}_K(x)$.  
    In particular,   for all $j \in \text{Act}_K(x)$,  we conclude the existence of a neighborhood $\mathcal{N}(x)$, 
    $c_j> 0$,  and $\alpha_j \in [0,1)$ such that,  
    $$  \forall (x_{1j},x_2) \in  ((\partial K_j \backslash C) \cap \mathcal{N}(x)) \times  ((\partial C \backslash \partial K) \cap \mathcal{N}(x)),  $$
    there exists $(v_j,w_j) \in  T_{\partial K_j} (x_{1j}) \times T_{C} (x_2)$ such that 
    \begin{align*}
      x_2 - x_{1j} & \in  v_j - w_j + \alpha_j (|x_{1j} - x_2|) \mathbb{B},\\
      |[v_j^\top ~ w_j^\top]| & \leq c_j |x_{1j} - x_2|.
    \end{align*}
    Now, to verify ($\dagger$) for $(\partial K, C)$, we note that when $\mathcal{N}(x)$ is sufficiently small,  then,  for all $x_{1} \in ((\partial K \backslash C) \cap \mathcal{N}(x))$,   there exists $j_1 \in \text{Act}_K(x)$ such that $x_1 \in (\partial K_{j_1} \backslash C) \cap \mathcal{N}(x)$.  Hence, by letting $ \alpha := \max_{j \in \text{Act}_K(x)} \{\alpha_j \}$ ~~ and $c := \max_{j \in \text{Act}_K(x)} \{c_j \}$, 
     we conclude the existence of 
    $$ (v,w) \in  T_{\partial K_{j_1}} (x_{1}) \times T_{C} (x_2) \subset T_{\partial K} (x_{1}) \times T_{C} (x_2) $$ 
    such that $|[v^\top ~~ w^\top]| \leq c |x_{1} - x_2|$ and 
    $ x_2 - x_{1} \in  v - w + \alpha (|x_{1} - x_2|) \mathbb{B}$. 

    \ifitsdraft
        Let us now verify ($\ddagger$) in Section~\ref{subsec_alternatives} for the sets  $(\partial K_j,C)$,  instead of the sets $(\partial K, C)$,  for all $j \in \text{Act}_K(x)$.  To do so,  we note that the set $\partial K$ can be expressed as $\partial K = \bigcup_{j=1}^{d_{K}} \partial K_j$.  
        Next,  given $x \in \partial K$, we note that $T_{\partial K}(x)  = \bigcup_{j=1}^{d_{K}} T_{\partial K_j}(x) = \bigcup_{j\in \text{Act}_K(x)} T_{\partial K_j}(x)$. Furthermore, given $x \in \mathcal{P}$ and $j \in \text{Act}_K(x)$, under~\eqref{eq1needed}-\eqref{eq2needed}, we apply the fifth item of Lemma~\ref{lem_prac_properties} in the Appendix to conclude that
        
    \else
    
        Let us now verify ($\ddagger$) for the sets  $(\partial K_j,C)$,  instead of the sets $(\partial K, C)$,  for all $j \in \text{Act}_K(x)$.  To do so,  we note that the set $\partial K$ can be expressed as $\partial K = \bigcup_{j=1}^{d_{K}} \partial K_j$.  
        Next,  given $x \in \partial K$, we note that $T_{\partial K}(x)  = \bigcup_{j=1}^{d_{K}} T_{\partial K_j}(x) = \bigcup_{j\in \text{Act}_K(x)} T_{\partial K_j}(x)$. Furthermore, given $x \in \mathcal{P}$ and $j \in \text{Act}_K(x)$, under~\eqref{eq1needed}-\eqref{eq2needed}, we apply the fifth item of Lemma~\ref{lem_prac_properties} in the Appendix to conclude that
    \fi   
    \begin{align} \label{eqTC} 
        T_{\partial K_j}(x) = C_{\partial K_j}(x) \qquad \forall j \in \text{Act}_K(x).   
    \end{align}
        
    Now,  we consider $v_j$ from Assumption~\ref{ass_transversalityPra} and a sufficiently small neighborhood, $\mathcal{N}(x)$,  around $x$ such that
    \begin{align*} 
    \nabla g_{C}^k(x)^\top v_j & - \nabla g_{C}^k(x)^\top (x_2 - x_{1j}) \leq 0 
    \quad 
    \forall k \in\mathrm{Act}_{C}(x), 
    \\ 
    \forall (x_{1j},x_2) & \in ((\partial K_j \backslash C) \cap \mathcal{N}(x)) \times ((\partial C \backslash \partial K) \cap \mathcal{N}(x)). 
    \end{align*}
      Such a neighborhood $\mathcal N(x)$ exists by continuity  and since  
    $\nabla g_{C}^k(x)^\top v_j < 0$  for all $k \in \mathrm{Act}_{C}(x)$.  Furthermore, using~\eqref{eqTC}, we conclude that $v_j \in  T_{\partial K_j}(x) = C_{\partial K_j}(x)$ and 
    $w_j := v_j - (x_2 - x_{1j}) \in  C_{C}(x)$.
    Hence, 
    \[ \forall (x_{1j},x_2) \in ((\partial K_j \backslash C) \cap \mathcal{N}(x)) \times ((\partial C \backslash \partial K) \cap \mathcal{N}(x)), \]
    we are able to find
    $(w_j,v_j) \in C_{C}(x) \times C_{\partial K_j}(x)$ such that 
    $ x_2 - x_{1j} = v_j - w_j$. 
\end{proof}

\section{LMI Certificates}
\label{sec:app-lmi}

In this section, we focus on polytopic sets $C$ and $K$ and polytopic set-valued maps $F$. We seek a 
numerically-tractable characterization of forward invariance in the form of finite LMIs.  To this end, we specify
\begin{equation} \label{Eq45}
\Sigma_P ~ : ~	\dot x \in F(x) := \mathrm{co}\{A_j x\}_{j=1}^N \qquad x \in C,
\end{equation}
where $A_j \in \mathbb{R}^{n\times n}$ for all $j \in [N]$, $K := \{ x \in \mathbb{R}^n : A_K x \le b_K \}$, and 
$ C := \{ x \in \mathbb{R}^n : A_C x \le b_C \}$,
	with 
	\begin{align*}
    A_K & := [a^{1}_K \ \cdots \ a^{d_K }_K]^\top, \quad 
	b_K := (b^1_K,\dots,b^{d_K}_K)^\top, 
	\\
    A_C & := [a^{1}_C \ \cdots\ a^{d_C }_C]^\top,
    \quad 
	b_C := (b^1_C,\dots,b^{d_C}_C)^\top,
    \end{align*}
and $a^i_K, a^k_C \in \mathbb{R}^n\setminus\{0\}$ for all $i \in [d_K]$ and 
	$k \in [d_C]$.

Finally, for $i \in [d_K]$ and $k \in [d_C]$, we let 
\begin{align*}
\partial K_i & := \{x\in K : a^{i\top}_K x = b^i_K\}, 
\\
\partial C_k & := \{x\in C : a^{k \top}_C x = b^k_C \}.
\end{align*}

Note that Assumption~\ref{ass_transversality}
(transversality) holds in this case by virtue of Lemma~\ref{lem_trans_poly} in the Appendix. 

We show later that Assumption \ref{ass_decisiveness_pra} holds in this case under the following property.
\begin{itemize}
\item[$(\spadesuit)$]
		For every $i\in[d_K]$ and $k\in[d_C]$ such that
		$(\partial K_i\cap\partial C_k)\cap\mathcal P \neq \emptyset$,
		there exist $\alpha,\beta\in\mathbb R^N$,
		$\lambda\in\mathbb R_{\geq0}^{Nd_K}$,
		$\mu\in\mathbb R_{\geq0}^{Nd_C}$, and
		$\sigma,\delta,\rho \geq 0$
		such that, for $\mathcal A:=[A_1~\cdots~A_N]$,
		\begin{align}
		\hspace{-0.4cm}	0={}&(\alpha\otimes I_n)a_K^i
			+(\beta\otimes I_n)a_C^k +(I_N\otimes A_K^\top)\lambda
			\nonumber\\
		\hspace{-0.4cm}	&
			+(I_N\otimes A_C^\top)\mu +\sigma \mathcal A^\top a_K^i
			+\delta \mathcal A^\top a_C^k
			\nonumber \\
            \hspace{-0.4cm}
			&-\rho \mathcal A^\top(a_K^i+a_C^k),
			\label{eq:decisiveness_exact_eq}\\
			\hspace{-0.4cm} 0>{}&\alpha b_K^i+\beta b_C^k
			+(I_N\otimes b_K^\top)\lambda
			+(I_N\otimes b_C^\top)\mu.
			\label{eq:decisiveness_exact_ineq}
		\end{align}
\end{itemize}

  The process of selecting the pairs $(i,k)$ such that $(\partial K_i \cap \partial C_k) \cap \mathcal{P} \neq \emptyset$ renders $(\spadesuit)$ not numerically tractable. This is particularly due to the latter two negated propositions in the definition of $\mathcal{P}$. 
  That said, replacing the latter sets by an 
  over-approximation, considering the latter two negated propositions for specific values of $\epsilon$ and specific choices of $S$ (such us $S= K$, $S = K \cap C$, $S= \partial K \cap \partial C$), or even not considering them at all, would lead to a numerically tractable formulation of  $(\spadesuit)$, although more restrictive in this case.

We are now ready to derive a  Nagumo-type characterization, as a finite set of LMIs. In particular, we show that the proposed LMI condition ensures both Assumption \ref{ass_propagation} and the 
Nagumo-type condition \eqref{eqconstNag}. 
On the other side, when 
\begin{align} \label{eqextra}
\partial K_i \cap\operatorname{int}(C)\neq\emptyset \qquad \forall
	i\in[d_K],
 \end{align}   
we show that  \eqref{eqconstNag} ensures  Assumption \ref{ass_propagation} and  the proposed LMIs certificate.

\begin{theorem}
	\label{thm3_LMI}
	Consider $\Sigma_P$ in (\ref{Eq45}) and the polytopes $K$ and $C$. Suppose that
	condition~$(\spadesuit)$ and \eqref{eqextra} hold. Then, $K$ is forward invariant for $\Sigma_P$ if
	and only if 
    \begin{itemize}
    \item  $\forall j\in[N]$, $\exists M_j \in \mathbb R_{\geq0}^{d_K\times d_C}$ and 
	$H_j\in\mathbb R^{d_K\times d_K}$, with ${H_j}_{ik} \geq 0$ for all $i \neq k$, such that
	\begin{equation}
		\hspace{-0.6cm} A_KA_j=H_jA_K+M_jA_C,\quad 
		H_jb_K+M_jb_C\leq0.
		\label{eq:lmi_matrix_ineq}
	\end{equation}
    \end{itemize}
\end{theorem}

\begin{proof}
    The proof relies on Theorem~\ref{theo_main_result}. First, the inequality
	$d_H(F(x),F(y))
	\leq\max_{j\in[N]}\|A_j\|\|x-y\|$
	shows that $F$ in \eqref{Eq45} is globally Lipschitz in the Hausdorff distance.
	In particular, $F$ is lower semicontinuous and one-sided locally
	Lipschitz.
	
	\emph{Necessity.}
	Suppose that $K$ is forward invariant for $\Sigma_P$. Fix
	$i\in[d_K]$ and $j\in[N]$. By the necessity part of
	Theorem~\ref{theo_main_result},
	$F(x)\subset T_K(x)$ for every
	$x\in\partial K\cap\operatorname{int}(C)$. Since
	$A_jx\in F(x)$, the tangent-cone representation of $K$ gives~\cite[Example 5.2.6-b)]{hiriart_urrutyFundamentalsConvexAnalysis2004}
	\begin{equation}
		(a_K^i)^\top A_jx\leq0
		\qquad \forall 
		x \in \partial K_i \cap \operatorname{int}(C).
		\label{eq:necessity-interior-facet}
	\end{equation}
	
	Choose $x_i^\circ\in
	\partial K_i\cap\operatorname{int}(C)$. For
	$x\in\partial K_i\cap C$ and $\tau\in(0,1]$, define
	$x_\tau:=(1-\tau)x+\tau x_i^\circ$. By convexity of the set $C$ and $\partial K_i$, $x_\tau\in\partial K_i\cap\operatorname{int}(C)$, and
	$x_\tau\to x$ as $\tau\to0^+$. Therefore,
	\eqref{eq:necessity-interior-facet} and continuity yield
	$(a_K^i)^\top A_jx \leq 0$ for all $x \in \partial K_i\cap C$.  Equivalently,
	$$ \left[ a_K^{i \top} x=b_K^i, ~ A_K x \leq b_K, ~
	A_C x \leq b_C \right] \Rightarrow a_K^{i\top} A_j x \leq 0. $$
	
	By Farkas' lemma
	\cite[Thm.~4.6, p.~165]{bertsimasIntroductionLinear1997},
	there exist $\Gamma_{ji}\in\mathbb R$,
	$\Lambda_{ji}\in\mathbb R_{\geq0}^{d_K}$, and
	$N_{ji}\in\mathbb R_{\geq0}^{d_C}$ such that
	\begin{align}
		A_j^\top a_K^i
		& =\Gamma_{ji}a_K^i
		+A_K^\top\Lambda_{ji}
		+A_C^\top N_{ji},
		\label{eq:facet-dual-equality}\\
		0  & \geq \Gamma_{ji}b_K^i
		+\Lambda_{ji}^\top b_K
		+N_{ji}^\top b_C.
		\label{eq:facet-dual-inequality}
	\end{align}
	Define
	$\Gamma_j:=\operatorname{diag}_{i\in[d_K]}\{\Gamma_{ji}\}$,
	$\Lambda_j:=[\Lambda_{j1}\ \cdots\ \Lambda_{jd_K}]$, and
	$N_j:=[N_{j1}\ \cdots\ N_{jd_K}]$. Stacking
	\eqref{eq:facet-dual-equality}--%
	\eqref{eq:facet-dual-inequality} over $i\in[d_K]$ gives
	$A_KA_j=(\Gamma_j+\Lambda_j^\top)A_K+N_j^\top A_C$ and $(\Gamma_j+\Lambda_j^\top)b_K+N_j^\top b_C\leq0$.
	Finally, set $H_j:=\Gamma_j+\Lambda_j^\top$ and
	$M_j:=N_j^\top$. Then, $M_j\geq0$, and the off-diagonal entries
	of $H_j$ are non-negative because $\Lambda_j\geq0$. Thus, \eqref{eq:lmi_matrix_ineq} holds.
	
	\emph{Sufficiency.}
	Suppose that \eqref{eq:lmi_matrix_ineq} holds. Let
	$x\in \partial K \cap C$. 
	Fixing $i\in\operatorname{Act}_K(x)$, for each
	$j\in[N]$, we let $H_{ji}$ and $M_{ji}$ denote the $i$th row of
	$H_j$ and $M_j$. Since $(a_K^i)^\top x=b_K^i$,
	$A_Kx\leq b_K$, and $H_j$ is Metzler, we have
	$H_{ji}A_Kx\leq H_{ji}b_K$. Moreover, $M_j\geq0$ and
	$A_Cx\leq b_C$ imply $M_{ji}A_Cx\leq M_{ji}b_C$. Hence,
	$$ a_K^{i\top} A_jx
	=H_{ji}A_Kx+M_{ji}A_Cx
	\leq(H_jb_K+M_jb_C)_i\leq0. $$
	
	Where the equality follows from \eqref{eq:lmi_matrix_ineq}. Now consider and $\eta\in F(x)$. Writing $\eta=\sum_{j=1}^N\theta_jA_jx$, with
	$\theta\in\mathbb R_{\geq0}^N$ and
	$\sum_{j=1}^N\theta_j=1$, gives
	$(a_K^i)^\top\eta
	=\sum_{j=1}^N\theta_j a_K^{i\top} A_jx\leq0$.
	Since this holds for every
	$i\in\operatorname{Act}_K(x)$, the tangent-cone representation
	of $K$ gives
	\begin{equation}
		F(x)\subset T_K(x)
		\qquad \forall x\in K\cap C.
		\label{eq:F-contained-TK}
	\end{equation}
	In particular,
	$F(x)\cap T_C(x)\subset T_K(x)$ for every $x\in\mathcal P$,
	so Assumption~\ref{ass_propagation} holds. Moreover,
	\eqref{eq:F-contained-TK} implies the Nagumo-type
	condition~\eqref{eqconstNag}. By Proposition~\ref{prop:decisiveness-iff},
	condition~$(\spadesuit)$ implies Assumption~\ref{ass_decisiveness_pra}, and  thus 
	Assumption~\ref{ass_decisiveness} follows by Proposition \ref{prop_nontang_prac}. Finally, Assumption~\ref{ass_transversality} holds in this case
	by Lemma~\ref{lem_trans_poly}. 
    Therefore, forward invariance of $K$ for $\Sigma_P$ follows using Theorem~\ref{theo_main_result}.
\end{proof}

\begin{remark}
Using Nagumo's theorem for unconstrained systems with $F(x) := A_1 x$, it is stated in \cite[Proposition 2.1]{castelan1993invariant} that a polytopic set $K$ is forward invariant if and only if there exists a Metzler  matrix, a matrix in which all the 
off-diagonal components are nonnegative, and $H_1 \in \mathbb{R}^{d_K \times d_K}$ such that $A_K A_1 = H_1 A_K$ and $H_1  b_K \le 0$.	This is exactly \eqref{eq:lmi_matrix_ineq} when $A_C = b_C = 0$. 
\end{remark}

\ifitsdraft

\else
An example illustrating $(\spadesuit)$ and the proposed 
finite-dimensional characterization is in \cite{arXiv_version_auto}. 
\fi

\subsection{Intermediate Result}

\begin{proposition}	\label{prop:decisiveness-iff}
Consider polytopic sets $K, C \subset \mathbb{R}^n$ and the constrained system $\Sigma_P$ in (\ref{Eq45}). 

Then, Assumption \ref{ass_decisiveness_pra} (decisiveness) holds if $(\spadesuit)$ holds. 
\end{proposition}

\begin{proof}
    Suppose that
Assumption~\ref{ass_decisiveness_pra} does not hold. 
Then, there exist $x \in\mathcal P$, 
$(i,k)\in\mathrm{Act}_K(x)\times\mathrm{Act}_C(x)$, and 
$\theta \in \Delta_N := 
\{\theta\in [0,1]^N:
\sum_{j=1}^N\theta_j=1\}$ such that
\[
a_K^{\ell \top} A_\theta x\leq0 ~
\forall\ell\in\mathrm{Act}_K(x),
~~ a_C^{q\top} A_\theta x \leq 0 
~ \forall q\in\mathrm{Act}_C(x),
\]
and $(a_K^i+a_C^k)^\top A_\theta x  \geq0$, where
$A_\theta:=\sum_{j=1}^N \theta_j A_j$.

Hence, negating Assumption \ref{ass_decisiveness_pra} suggests the existence of $x \in \mathcal{P}$, $(i,k)\in [d_K]\times [d_C]$, and $\theta \in \Delta_N$  such that
\begin{equation}
	\label{eq:negated-system}
	\begin{aligned}
	\hspace{-0.4cm}	a_K^{i\top} x & =b_K^i,\;
		a_C^{k\top} x=b_C^k,\;
		A_Kx\leq b_K,\; A_Cx\leq b_C,
        \\
	\hspace{-0.4cm}	a_K^{i\top} A_\theta x & \leq 0,
        \;
		a_C^{k\top} A_\theta x\leq0,\;
		-(a_K^i+a_C^k)^\top A_\theta x\leq0.
	\end{aligned}
\end{equation}
The rest of the proof shows that \eqref{eq:negated-system} is actually infeasible for all $x \in \mathbb{R}^n$, $(i,k)\in [d_K]\times [d_C]$ such that
		$(\partial K_i\cap\partial C_k)\cap\mathcal P \neq \emptyset$, and $\theta \in \Delta_N$. Thus,
	Assumption~\ref{ass_decisiveness_pra} holds.

System~\eqref{eq:negated-system} can be written as
$Ex=f$ and $G_\theta x\leq h$, 
$
E:=
\begin{bmatrix}
	a_K^i &
	a_C^k
\end{bmatrix}^\top$,
$
f:=
\begin{bmatrix}
	b_K^i &
	b_C^k
\end{bmatrix}^\top$, $h:=
\begin{bmatrix}
	b_K&
	b_C&
	0&
	0&
	0
\end{bmatrix}^\top$,
\[
G_\theta:=
\begin{bmatrix}
	A_K^\top &~
	A_C^\top &~
	A_\theta^\top a_K^i &~
	A_\theta^\top a_C^k &~
	-A_\theta^\top(a_K^i+a_C^k)
\end{bmatrix}^\top.
\]
	By 
    \ifitsdraft
    Lemma~\ref{lem:farkas-mixed},
    \else
    \cite[Lemma 24]{arXiv_version_auto}, which is a slight generalization of Theorem~4.6 in \cite[page 165]{bertsimasIntroductionLinear1997},
    \fi
    this system is infeasible (i.e., there is no $x \in \mathbb{R}^n$ such that \eqref{eq:negated-system} holds) if and only if 
	there exist $y_\theta\in\mathbb R^2$ and
	$z_\theta\in\mathbb R_{\geq0}^{d_K+d_C+3}$ such that
	\begin{equation}
		E^\top y_\theta+G_\theta^\top z_\theta=0,
		\qquad
		f^\top y_\theta+h^\top z_\theta<0.
		\label{eq:Farkas-theta}
	\end{equation}
	
	We now show that condition~$(\spadesuit)$ provides such
	multipliers for every $\theta\in\Delta_N$ and for every $(i,k)\in [d_K] \times [d_C]$ such that $(\partial K_i\cap\partial C_k) \cap \mathcal P \neq \emptyset$. 
    
    We let  $\alpha=\operatorname{col}_{j\in[N]}(\alpha_j)$,
	$\beta=\operatorname{col}_{j\in[N]}(\beta_j)$,
	$\lambda=\operatorname{col}_{j\in[N]}(\lambda_j)$, and
	$\mu=\operatorname{col}_{j\in[N]}(\mu_j)$, where
	$\lambda_j\in\mathbb R_{\geq0}^{d_K}$ and
	$\mu_j\in\mathbb R_{\geq0}^{d_C}$.
		The $j$-th block of
	\eqref{eq:decisiveness_exact_eq} and the $j$-th component of
	\eqref{eq:decisiveness_exact_ineq} give, respectively,
	\begin{align}
		0 & =
		\alpha_j a_K^i+\beta_j a_C^k
		+A_K^\top\lambda_j+A_C^\top\mu_j
		+\sigma A_j^\top a_K^i
		\nonumber\\
		& + \delta A_j^\top a_C^k-\rho A_j^\top(a_K^i+a_C^k),
		\label{eq:vertex-Farkas-equality}
        \\
		0 & >
		\alpha_j b_K^i+\beta_j b_C^k
		+b_K^\top\lambda_j
		+b_C^\top\mu_j.\label{eq:vertex-Farkas-inequality}
	\end{align}
	
	Define
	$\alpha_\theta:=\sum_{j=1}^N\theta_j\alpha_j$,
	$\beta_\theta:=\sum_{j=1}^N\theta_j\beta_j$,
	$\lambda_\theta:=\sum_{j=1}^N\theta_j\lambda_j$, and
	$\mu_\theta:=\sum_{j=1}^N\theta_j\mu_j$.
	Because $\theta_j\geq0$, $\lambda_j\geq0$, and $\mu_j\geq0$,
	we have $\lambda_\theta\geq0$ and $\mu_\theta\geq0$.
	Multiplying \eqref{eq:vertex-Farkas-equality} and \eqref{eq:vertex-Farkas-inequality} by $\theta_j$,
	summing over $j\in[N]$, and using
	$\sum_{j=1}^N\theta_j=1$ yield
	\begin{align}
		0={}&
		\alpha_\theta a_K^i+\beta_\theta a_C^k
		+A_K^\top\lambda_\theta+A_C^\top\mu_\theta
		+\sigma A_\theta^\top a_K^i
		\nonumber\\
		&+\delta A_\theta^\top a_C^k-\rho A_\theta^\top(a_K^i+a_C^k). \label{eq:theta-Farkas-equality}
		\\
		0>&\alpha_\theta b_K^i+\beta_\theta b_C^k +\lambda_\theta^\top b_K+\mu_\theta^\top b_C. \label{eq:theta-Farkas-inequality}
	\end{align}
	Setting
	$
	y_\theta :=
	\begin{bmatrix}
		\alpha_\theta
        &
		\beta_\theta
	\end{bmatrix}^\top$ and 
	$z_\theta:=
	\begin{bmatrix}
		\lambda_\theta
        &
		\mu_\theta
        &
		\sigma_{ik}
        &
		\delta_{ik}
        &
		\rho_{ik}
	\end{bmatrix}^\top$,
	we have $z_\theta\geq0$, and
	\eqref{eq:theta-Farkas-equality}--%
	\eqref{eq:theta-Farkas-inequality} imply
	\eqref{eq:Farkas-theta}. 
\end{proof}

\ifitsdraft

\subsection{Numerical example}
\label{sec:app-lmi-numeric}

We illustrate Theorem~\ref{thm3_LMI} using 
$\Sigma_P$ with $n=N=2$, and
\begin{equation}
\label{eq:num-vertices}
 A_1=
 \begin{bmatrix}
  -\frac12 & 1\\
  1 & -\frac12
 \end{bmatrix},
 \qquad
 A_2=
 \begin{bmatrix}
  -1&0\\
  0&-1
 \end{bmatrix}.
\end{equation}
Thus,
\[
 F(x)=\operatorname{co}\{A_1x,A_2x\},
 \qquad x\in C.
\]
For the extension in the definition of $\mathcal{P}$, we use  
$\widehat F(x):=\operatorname{co}
\{A_1x,A_2x\}$ for all $x \in\mathbb R^2$. 

Consider the polytopic sets $K$ and $C$ illustrated in Fig. \ref{fig:Polytopes} and defined as :
\begin{align}
    K := \{x\in\mathbb R^2:A_Kx\leq b_K\}, \; C := \{x\in\mathbb R^2:A_Cx\leq b_C\},
\end{align}
where
\begin{align}
A_K = \begin{bmatrix}
1&0\\0&1\\0&-1\\-1&0
\end{bmatrix},&
b_K=\mathbf 1_4,
~
A_C = \begin{bmatrix}
1&1\\-1&-1\\1&0\\0&1\\-1&0\\0&-1
\end{bmatrix}, ~
b_C=\begin{bmatrix} 
\tfrac32 \\ \tfrac32 \\ 2 \\ 2 \\ 2\\ 2 \end{bmatrix}.
\label{eq:num-C-data}
\end{align}
Equivalently,
$ K=[-1,1] \times [-1,1]$
and
\[
 C=
 \left\{
 x\in[-2,2] \times [-2,2]:
  |x_1+x_2|\leq\frac32
 \right\}.
\]


We further consider 
\[
 \bar x^1=
 \begin{bmatrix}1\\0\end{bmatrix},
 \bar x^2=
 \begin{bmatrix}0\\1\end{bmatrix},
 \bar x^3=
 -\begin{bmatrix}0\\1\end{bmatrix},
 \bar x^4=
 -\begin{bmatrix}1\\0\end{bmatrix}.
\]
These points satisfy
$
 \bar x^i\in\partial K_i\cap\operatorname{int}(C)$
 for all $i\in[4]$.
Consequently,
$\partial K_i\cap\operatorname{int}(C)\neq \emptyset$
for all $i\in[4]$,
and condition~\eqref{eqextra} holds.


We observe that
$
 \partial K\cap\partial C
 =
 \{x^{11},x^{21},x^{32},x^{42}\},
$
where
\begin{equation}
\label{eq:num-boundary-points}
 \begin{aligned}
 x^{11}=
 \begin{bmatrix}
  1\\ \frac12
 \end{bmatrix},
~
 x^{21} =
 \begin{bmatrix}
  \frac12\\1
 \end{bmatrix},
 ~ x^{32} =
 \begin{bmatrix}
  -\frac12\\-1
 \end{bmatrix},
~
 x^{42} =
 \begin{bmatrix}
  -1\\-\frac12
 \end{bmatrix}.
 \end{aligned}
\end{equation}
The superscript \(ik\) indicates that the \(i\)-th facet of \(K\)
and the \(k\)-th facet of \(C\) are active. Each of these points also
belongs to \(\operatorname{cl}(C\backslash K)\).


At \(x^{11}\) and \(x^{21}\), the first constraint of \(C\) is
active, and hence
$
 T_C(x)
 =
 \left\{
 v\in\mathbb R^2:v_1+v_2\leq0
 \right\}$. 
At \(x^{32}\) and \(x^{42}\), the second constraint of \(C\) is
active, and therefore
$
 T_C(x)
 =
 \left\{
 v\in\mathbb R^2:v_1+v_2\geq0
 \right\}$.
We observe that
\[
 A_1x^{ik}\notin T_C(x^{ik}),
 \qquad
 A_2x^{ik}\in\operatorname{int}\bigl(T_C(x^{ik})\bigr),
\]
for every \(x^{ik}\) in~\eqref{eq:num-boundary-points}. Since
\(F(x^{ik})\) is convex, it contains the direction
$
 \eta^{ik}
 :=
 \frac{2}{3}A_1x^{ik}
 +\frac{1}{3}A_2x^{ik}$.
Hence,
\[
\begin{aligned}
 \eta^{11}
 =
 \begin{bmatrix}
  -\frac13\\[1mm]\frac13
 \end{bmatrix},
~
 \eta^{21}
 =
 \begin{bmatrix}
  \frac13\\[1mm]-\frac13
 \end{bmatrix}, ~
 \eta^{32}
 =
 \begin{bmatrix}
  -\frac13\\[1mm]\frac13
 \end{bmatrix},
~
 \eta^{42}
 =
 \begin{bmatrix}
  \frac13\\[1mm]-\frac13
 \end{bmatrix}.
\end{aligned}
\]
Consequently,
$\eta^{ik}_1+\eta^{ik}_2=0$
at all four points. Hence,
\[
 \eta^{ik}
 \in
 F(x^{ik})\cap T_{\partial C}(x^{ik})
 \subseteq
 F(x^{ik})\cap T_C(x^{ik}),
\]
which proves that
$\overline{\left[
 F(x^{ik})\subset
 \mathbb R^2\backslash T_C(x^{ik})
 \right]}$
for all \(x^{ik}\).

We now examine the second condition in the definition of
\(\mathcal P\). Each vector
\(A_1x^{ik}\) is tangent to the corresponding active facet of \(K\):
\[
\begin{aligned}
 \begin{bmatrix}1&0\end{bmatrix}A_1x^{11}&=0,
&
 \begin{bmatrix}0&1\end{bmatrix}A_1x^{21}&=0,
\\
 \begin{bmatrix}0&-1\end{bmatrix}A_1x^{32}&=0,
&
 \begin{bmatrix}-1&0\end{bmatrix}A_1x^{42}&=0.
\end{aligned}
\]
Since exactly one constraint of \(K\) is active at each of these
points, it follows that
\[
 A_1x^{ik}
 \in
 T_{\partial K}(x^{ik})
 =
 T_K(x^{ik})\backslash D_K(x^{ik}).
\]
In particular,
$A_1x^{ik}\notin D_K(x^{ik})$,
and, because \(A_1x^{ik}\in F(x^{ik})\), we conclude that
$
 F(x^{ik})\nsubseteq D_K(x^{ik})$.

We next verify the third condition in the definition of
\(\mathcal P\). For each \(x^{ik}\), let
$ z^{ik}(t):=e^{A_1t}x^{ik} $
denote the unique solution of
$\dot z^{ik}=A_1 z^{ik}$,  $z(0)=x^{ik}$.
Define
\[
 q(t):=\frac34e^{t/2}+\frac14e^{-3t/2},
 \quad
 r(t):=\frac34e^{t/2}-\frac14e^{-3t/2}.
\]
Direct calculation gives
\[
\begin{aligned}
 z^{11}(t)&=
 \begin{bmatrix}q(t)\\r(t)\end{bmatrix},
&
 z^{21}(t)&=
 \begin{bmatrix}r(t)\\q(t)\end{bmatrix},
\\
 z^{32}(t)&=
 \begin{bmatrix}-r(t)\\-q(t)\end{bmatrix},
&
 z^{42}(t)&=
 \begin{bmatrix}-q(t)\\-r(t)\end{bmatrix}.
\end{aligned}
\]
Since
$ q(t)=1+\frac38t^2+O(t^3), $
the coordinates corresponding to the active facets of \(K\) satisfy
\begin{equation}\label{eq:num-immediate-leaving}
    \begin{aligned}
     &z^{11}_1(t)=1+\frac38t^2+O(t^3), ~~ z_2^{21}(t)=1+\frac38t^2+O(t^3),\\
     & z_2^{32}(t)=-1-\frac38t^2+O(t^3),~~ z_1^{42}(t)=-1-\frac38t^2+O(t^3).
    \end{aligned}
\end{equation}
Therefore, for every sufficiently small \(t>0\), the first two
trajectories violate an upper-facet inequality of \(K\), whereas the
last two violate a lower-facet inequality. Hence, all four trajectories
leave \(K\) immediately.
Therefore, for every $x^{ik}$, there is  no such pair \((\epsilon,S)\) such that
\[
 \overline{\left[
 \begin{aligned}
  &\exists\epsilon>0,\ \exists S\subseteq K\text{ closed}:\\
  &x^{ik}\in S,\qquad
  \widehat F(y)\subseteq T_S(y)
  \quad
  \forall y\in\mathbb B_\epsilon(x^{ik})\cap S
 \end{aligned}
 \right]}
\]
holds at all four points.

Finally, \(C\) is convex at every \(x^{ik}\). However, since
$A_1x^{ik}\in T_K(x^{ik})$ and $A_1x^{ik}\notin T_C(x^{ik})$,
we have $F(x^{ik})\nsubseteq T_{C\cap K}(x^{ik})$.
Therefore, the last negated condition in the definition of
\(\mathcal P\) also holds at every \(x^{ik}\). We conclude that
\begin{equation}
\label{eq:num-P}
 \mathcal P
 =
 \{x^{11},x^{21},x^{32},x^{42}\}.
\end{equation}


By~\eqref{eq:num-P}, the active facet pairs to be checked in
condition~\((\spadesuit)\) are
$ (i,k) \in \{ (1,1),(2,1),(3,2),(4,2) \}$.

Let $\mathcal A=[A_1\ A_2]$.
For every one of these four pairs, choose
\begin{equation}
\label{eq:num-spadesuit-certificate}
 \alpha=
 \begin{bmatrix}
  -\frac92 & -3
 \end{bmatrix}^\top,
 \quad
 \beta=
 \begin{bmatrix}
  \frac73 & \frac43
 \end{bmatrix}^\top,
 \quad
 \lambda=0_8,
 \quad
 \mu=0_{12},
\end{equation}
together with
$
 \sigma=0$, $\delta=\frac{13}{3}$,
 $\rho=3$. Direct substitution gives
\begin{align*}
    0_4={}& (\alpha\otimes I_2)a_K^i+(\beta\otimes I_2)a_C^k+(I_2\otimes A_K^\top)\lambda+\\
    & (I_2\otimes A_C^\top)\mu +\sigma\mathcal A^\top a_K^i
     +\delta\mathcal A^\top a_C^k
     -\rho\mathcal A^\top(a_K^i+a_C^k),
\end{align*}
\[
 \alpha b_K^i+\beta b_C^k
 +(I_2\otimes b_K^\top)\lambda
 +(I_2\otimes b_C^\top)\mu
 =
 \begin{bmatrix}-1\\-1\end{bmatrix}
 <0.
\]
Thus, condition~\((\spadesuit)\) holds for every active facet pair
intersecting \(\mathcal P\).


For \(j=1\), consider
$
 H_1=-\frac32 I_4 $
and
$
 M_1=
 \begin{bmatrix}
  1&0&0&0&0&0\\
  1&0&0&0&0&0\\
  0&1&0&0&0&0\\
  0&1&0&0&0&0
 \end{bmatrix}$.
Then \(M_1\geq0\), \(H_1\) is Metzler, and
\[
 A_KA_1=H_1A_K+M_1A_C,
 \qquad
 H_1b_K+M_1b_C=0_4.
\]

For \(j=2\), take
$
 H_2=-I_4$
and 
 $M_2=0_{4\times6}$.
Again, \(M_2\geq0\), \(H_2\) is Metzler, and
\[
 A_KA_2=H_2A_K+M_2A_C,
 \quad
 H_2b_K+M_2b_C=-\mathbf 1_4<0.
\]

Since~\((\spadesuit)\)
and~\eqref{eqextra} also hold, Theorem~\ref{thm3_LMI} guarantees that
\(K\) is forward invariant for the constrained system 
\(\Sigma_P\).

\begin{figure}
    \centering
    \includegraphics[width=1\linewidth]{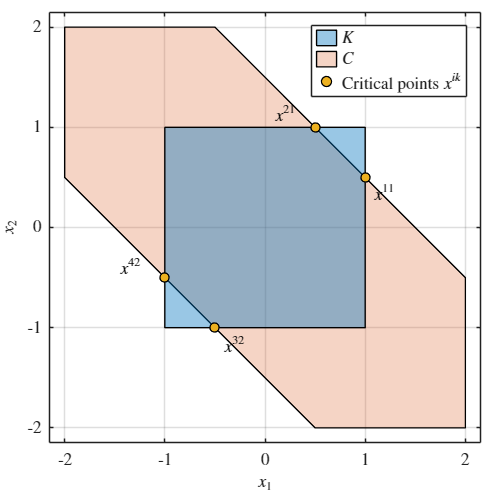}
    \caption{Polytopes \(K\) and \(C\), and the critical points \(x^{ik}\) associated with their active facet pairs.}
    \label{fig:Polytopes}
\end{figure}

\fi

\section{Conclusion}

This paper provided a fine extension of Nagumo’s invariance criterion for forward invariance for constrained systems. For the proposed condition to imply forward invariance, in the general case, three  assumptions are imposed. 
The main result is further specialized to finitely-represented sets, which allowed us to drop the propagation assumption and to provide simpler, and more tractable, formulation of the remaining ones. In addition to a finite-LMI characterization under polytopic sets and dynamics. For future work, we envision dropping the decisiveness assumption (Assumption \ref{ass_decisiveness}) at the expense of hardening the transversality condition (Assumption \ref{ass_transversality}) for further specifying the dynamics $F$.

\bibliography{biblio}
\bibliographystyle{unsrt}

\section*{Appendix}


\ifitsdraft

\textit{Proof of the first item in Theorem \ref{theo_main_result}:}
    Let $v \in F(x)$ for some $x \in \partial K \cap \text{int}(C)$.  Using the continuity of $F$ plus the property of its images (the one-sided Lipschitzness
    property is not needed here), we show the existence of $T > 0$ and a solution $\phi$ starting from $x$ such that 
    $$ \phi(t) = x + t (v + r(t)) \qquad \forall t \in [0,T], $$
    for some continuous function $t \mapsto r(t)$ satisfying $ \lim_{t \rightarrow 0} r(t) = 0$.
    Indeed, using Michael's selection theorem~\cite{michael1956continuous}, we deduce the existence of a continuous map $y \mapsto u(y) \in F(y)$, with $\text{dom}\, u$ containing a neighborhood of $x$,  such that $u(x) = v$. 
    Hence, there exists $T > 0$ and a $\mathcal{C}^1$ function $\phi : [0,T] \rightarrow \mathbb{R}^n$  verifying 
     $$ \dot{\phi}(t) = u(\phi(t)) \in F(\phi(t))
     \quad \forall t \in [0,T] \quad \text{and} \quad \phi(0) = x.  $$  
    Since $x \in \text{int}(C)$, we conclude that, for $T>0$ small enough, we have $\phi([0,T]) \subset \text{int}(C)$ and $\phi$ is a solution to $\Sigma$. Furthermore, being continuously differentiable, the first-order approximation of $\phi$ verifies
    $$ \phi(t) = x + t (v + r(t))  \quad \forall t \in [0,T], $$
    for a continuous function $t \mapsto r(t)$ satisfying $ \lim_{t \rightarrow 0} r(t) = 0$. Finally, since $K$ is forward invariant, we conclude that 
    $$ \phi(t) = x + t (v + r(t)) \in K \cap C \subset K \qquad \forall t \in [0,T]. $$
    The latter is enough to conclude that $v \in T_{K}(x)$ by definition of the contingent cone. 
\hfill $\blacksquare$

\fi


\begin{lemma} \label{lem_first_order_leaving_implies_star}
    Consider closed subsets $C$, $K \subset \mathbb{R}^n$ such that Assumption~\ref{ass_transversality} holds. Let $x \in \partial K \cap \partial C$ and $\phi : \mathbb{R}_{\geq 0} \rightarrow \mathbb{R}^n$ be a function verifying
    \begin{equation}\label{eqsequence}
         \phi(t) \in \{ x + t (y  + O(t) z) : y \in \Pi,  z \in \mathbb{B} \}
           \quad \forall t \geq 0, 
    \end{equation}
    with $\lim_{t \rightarrow 0^+} O(t) = 0$ and the compact set $\Pi$ verifies 
    \begin{equation} \label{eq_pio_in_TparKcapCNotTparKcapparC}
        \Pi \subset T_{\partial K \cap C}(x) \backslash T_{\partial K \backslash \mathrm{int}(C)}(x). 
    \end{equation}
    Then ($\star$) is verified. 
\end{lemma}

\begin{proof}
    To find contradiction,  we assume the existence of  
    $$
    \{ t_i \}^\infty_{i=1} \subset \mathbb{R}_{> 0} \quad  \text{with} \quad  \lim_{i \rightarrow \infty} t_i = 0^+
    $$
    and a sequence $\{y_i\}_{i=1}^\infty$  such that
    \begin{equation}\label{eq_project_phi}
        y_i \in \text{Proj}_{\partial K}(\phi(t_i)) \subset \partial K \backslash \text{int}(C) ~~ \forall i \in \{1,2,\ldots\}.
    \end{equation}
    
    Using~\eqref{eqsequence},  we conclude the existence of $\{\nu_i\}_i \subset \Pi$ and $\{z_i\}_i \subset \mathbb{B}$ such that 
    \[ \phi(t_i) =  x + t_i \nu_i + t_i o(t_i) z_i   \qquad \forall i \in \{1,2,\ldots\}.\]
    Subtracting $y_i$ from both sides of the previous inclusion and taking the absolute value, we obtain
    \begin{equation}
    \label{eq_abs_phi_min_yi}
    \begin{split}
    & |\phi(t_i) - y_i| \geq  t_i \left|\left(\nu_i - \frac{y_i-x}{t_i}\right)\right| - t_i o(t_i).
    \end{split}
    \end{equation}   
    Through an appropriate sub-sequence, we let  
    $$ w := \lim_{i \rightarrow \infty} \left[ w_i := \frac{y_i - x}{t_i} \right], \quad \nu := \lim_{i \rightarrow \infty}  \nu_i.  $$
    We will show that $w \neq \nu$. In fact,  we assume that $w$ is finite, since otherwise $w \neq \nu$ would hold trivially. Now, since $x + t_i w_i = y_i \in \partial K \backslash \text{int}(C)$ for all $i \in \{1,2,\ldots\}$, we conclude that $w \in T_{\partial K \backslash \text{int}(C)}(x)$. 
    In parallel, from~\eqref{eq_pio_in_TparKcapCNotTparKcapparC}, we conclude that  
    $\nu \subset T_{\partial K \cap C}(x) \backslash T_{\partial K \backslash \text{int}(C)}(x)$. 
    
    As a result, for 
    $\epsilon := | \nu - w | > 0$, we can write 
    $$  |\phi(t_i) - y_i| \geq t_i \epsilon - r_i t_i \qquad \forall i \in 
    \{1,2,\ldots\}, $$
    $$ r_i :=  o(t_i) + | w - w_i | + |\nu - \nu_i| \qquad \forall i \in \{1,2,\ldots\}. $$
    Clearly, the  sequence 
    $\{r_i\}_{i}$ converges to $0$. 
    
    Now, we use Assumption~\ref{ass_transversality} and Lemma~\ref{lemtrans} to obtain
    \begin{equation}
    \label{eqderivvv}
        \begin{split}
        T^a_{\partial K \cap C}(x) &= T^a_{\partial K}(x) \cap T^a_{C}(x) \\
                & = T_{\partial K}(x) \cap T_{C}(x) = T_{\partial K \cap C}(x).
        \end{split}
    \end{equation}
    Hence,  since $\nu \in T_{\partial K \cap C}(x) $, we conclude that 
    $\nu \in T^a_{\partial K \cap C}(x)$. Thus, there exists $\{\hat{\nu}_i \}_{i} \subset \mathbb{R}^n$ such that, for the same sequence $\{t_i\}_i$, we have 
    $$ \lim_{i \rightarrow \infty} \hat{\nu}_i = \nu  \text{ and } x_i := x + t_i \hat{\nu}_i \in \partial K \cap C \quad  \forall i \in \{1,2,\ldots\}. $$
    That is, we established the existence of $\{x_i \}_i \subset \partial K \cap C$ such that
    $$ |\phi(t_i) - x_i| \leq t_i |\nu_i - \hat{\nu}_i|  +  t_i o(t_i) \qquad \forall i \in 
    \{1,2,\ldots\} $$
    and, at the same time, for the sequence $\{y_i\}_i \subset  \partial K \backslash \text{int}(C)$ verifying~\eqref{eq_project_phi}, 
    we have 
    $$  |\phi(t_i) - y_i| \geq t_i \epsilon - r_i t_i \quad \forall i \in \{1,2,\ldots\}. $$
    As a consequence, for $i$ large enough, it holds that 
    $$ |\phi(t_i) - x_i| < t_i \epsilon/2 \quad \text{and} \quad |\phi(t_i) - y_i| \geq (\epsilon/2) t_i. $$
    Hence, for $i$ large enough, it follows that
    $$ |\phi(t_i) - x_i| <  |\phi(t_i) - y_i| = |\phi(t_i)|_{\partial K}.  $$
    This contradicts~\eqref{eq_project_phi}, since
    \[
    |\phi(t_i)|_{\partial K} = |\phi(t_i) - y_i| \leq |\phi(t_i) - u| \qquad \forall u \in \partial K.
    \]
\end{proof}

The proof of the following lemma  leverages some   arguments used in the proof of the constrained inverse-function theorem; see~\cite[Sec. 3.4.2 and Thm. 3.4.5]{aubinSetValuedAnalysis2009}.

\begin{lemma}\label{lemtrans}
    Consider closed subsets $K_1$, $K_2 \subset \mathbb{R}^n$  and let $x \in \partial K_1 \cap \partial K_2$ such that $K_1$ is derivable at $x$, i.e.,
    \begin{equation}\label{eqderrive}
        T_{K_1}(x) = T^a_{K_1}(x). 
    \end{equation}
    We further assume that 
    \begin{itemize}
        \item[($\diamond$)]  $\exists \mathcal{N}(x)$ a neighborhood of $x$, $c>0$, and $\alpha \in [0,1)$  such that
    \begin{equation*}
    \begin{split}
        x_2 - x_1 \in & \ (v_1 - v_2) +  \alpha |x_1 - x_2| \mathbb{B} \\ 
        & \hspace{-2cm}
        \forall (x_1,x_2) \in \begin{aligned}[t]
      & ((\partial K_1 \backslash K_2) \cap \mathcal{N}(x)) \times ((\partial K_2 \backslash K_1)  \cap \mathcal{N}(x)),
            \end{aligned}
    \end{split}
    \end{equation*}
    for some $(v_1,v_2) \in  T_{K_1} (x_1) \times T_{K_2} (x_2)$ such that $|[v_1^\top ~ v_2^\top]| \leq c |x_1 - x_2|$.
    \end{itemize}
    Then, 
    \begin{align}
        T_{K_1}(x) \cap T_{K_2}(x) & = T_{K_1 \cap K_2}(x), \label{eq1lemtrans} \\
        T^a_{K_1}(x) \cap T^a_{K_2}(x) & = T^a_{K_1 \cap K_2}(x). \label{eq2lemtrans}
    \end{align}
\end{lemma}
\begin{proof}
    It is straightforward to show that $T_{K_1 \cap K_2}(x) \subset T_{K_1}(x) \cap T_{K_2}(x)$ and 
    $
    T^a_{K_1 \cap K_2}(x) \subset T^a_{K_1}(x) \cap T^a_{K_2}(x)$. 
    Hence, to establish~\eqref{eq1lemtrans}-\eqref{eq2lemtrans}, it is enough to show that 
    \begin{align}
    T_{K_1}(x) \cap T_{K_2}(x) \subset T_{K_1 \cap K_2}(x),  \label{eqtoprove1} \\
    T^a_{K_1}(x) \cap T^a_{K_2}(x) \subset T^a_{K_1 \cap K_2}(x). \label{eqtoprove2} 
    \end{align}
    In the sequel, we only verify~\eqref{eqtoprove1}. \eqref{eqtoprove2} can be verified using the exact same arguments. 
    
    Consider a vector $v \in T_{K_1}(x) \cap T_{K_2}(x)$ for some $x \in \partial K_1 \cap \partial K_2$. By definition of the contingent cone, there exist $\{t_i\}_i \rightarrow 0$ and 
    $\{w_i\}_i \rightarrow v$ such that $x_{2i} := x + t_i w_i \in K_2$. Now, using~\eqref{eqderrive}, we conclude that, for the same sequence $\{t_i\}_i$, there exists $ \{z_i\}_i \rightarrow v$ such that $x_{1i} := x + t_i z_i \in K_1$.   
    
    The next step consists in showing the existence of $\{ (x^*_{1i}, x^*_{2i} ) \}_i \subset K_1 \times  K_2$, $l >0$, and $\sigma : 
    \{1,2,\ldots\} \rightarrow \{1,2\}$ 
    such that, through an appropriate subsequence,  
    \begin{equation}\label{ineqinterm--}
         x^*_{\sigma(i)i}  \in K_1 \cap K_2 \qquad \forall i \in \{1,2, \ldots \},
    \end{equation}  
    \begin{equation}\label{ineqinterm}
    \begin{split}
        | x_{1i} - x^*_{1i} | + |x_{2i} - x^*_{2i} |  & \leq l |x_{1i} - x_{2i}|\\
        &\forall i \in \{1,2, \ldots \}.  
    \end{split}
    \end{equation}
    The latter would imply that
    $\lim_{i \rightarrow + \infty} v_i \left( := -\frac{x - x^*_{\sigma(i) i} }{t_i} \right) = v$
    and $x^*_{\sigma(i) i} = x + t_i v_i \in K_1 \cap K_2$, which is enough to show that $v \in T_{K_1 \cap K_2}(x)$. 
    
    Note that, when $x_{1i} \in K_2 \cap K_1$ (resp., $x_{2i} \in K_1 \cap K_2$), then we just take $x^*_{1i} = x^*_{2i} = x_{1i}$ (resp., $x^*_{1i} = x^*_{2i} = x_{2i}$) to conclude that the inequality in~\eqref{ineqinterm} holds for 
    $l = 1$ and $\sigma(i) = 1$ (resp., $\sigma(i) = 2$).
    
    The remainder of the proof is devoted to establishing the following auxiliary claim.
    \begin{claim}  \label{clmm1}
        Given 
        $\{x_i := (x_{1i},x_{2i})\}_{i} \subset (K_1 \backslash K_2) \times (K_2 \backslash K_1)$  such that $x_{i} \rightarrow (x,x)$, there exists $\{ (x^*_{1i}, x^*_{2i}) \}_i \subset K_1 \times K_2$, $l >0$, and  $\sigma : 
        \{1,2,\ldots\} \rightarrow 
        \{1,2\}$ 
        such that~\eqref{ineqinterm--} and~\eqref{ineqinterm} hold at least through a certain subsequence. 
    \end{claim}
    \underline{Proof of Claim~\ref{clmm1}:}
    We introduce the function 
    $$ f(x) := x_1 - x_2  \qquad \forall x:= (x_1,x_2) \in K_1 \times K_2. $$
    Then, applying the Ekeland Variational Principle (see, e.g., \cite[Theorem 3.3.1]{aubinSetValuedAnalysis2009}), 
    to the function $|f(\cdot)|$, we conclude that, for each $\varepsilon > 0$,  there exists  a sequence $\{ \bar{x}_i := (\bar{x}_{1i}, \bar{x}_{2i}) \} \subset K_1 \times K_2$ such that, for all $i \in \{1,2,\ldots\} $, we have, for all $y_i (\in K_1 \times K_2) \neq \bar{x}_i$,  
    \begin{equation}\label{equse-1}
        \begin{cases}
            |f(\bar{x}_i)| - |f(x_i)| \leq -\varepsilon | \bar{x}_i - x_i |, \\
            |f(\bar{x}_i)| - |f(y_i)| < \varepsilon | \bar{x}_i - y_i |. 
        \end{cases}
    \end{equation}
    
    Using the first inequality, we conclude that for $l :=  2/ \varepsilon$,
    $$ |\bar{x}_{1i} - x_{1i}| + |\bar{x}_{2i} - x_{2i}| \leq 2| \bar{x}_i - x_i | \leq 2|f(x_i)| / \varepsilon = l |x_{1i} - x_{2i}| $$
    Hence, we conclude that $\bar{x}_i \rightarrow (x,x)$. 
    
    The proof is completed if we show , for a certain choice of $\varepsilon > 0$, the existence of $\sigma : 
    \{1,2,\ldots\} \rightarrow \{1,2\}$ such that, after passing to an appropriate subsequence, we find
    $\{ \bar{x}_{\sigma(i) i} \}_i \subset K_1 \cap K_2$.  To find a contradiction, we assume that, for any choice of $\varepsilon > 0$, we have  
    \begin{equation}\label{eqcontradict}
    (\bar{x}_{1i}, \bar{x}_{2i}) \in (K_1 \backslash K_2) \times (K_2 \backslash K_1) \quad \forall i ~ \text{large}. 
    \end{equation}
    
    Next, we show the existence of $\bar{c} > 0$ and $\bar{\alpha} \in [0,1)$ such that,  for $ \bar{x}_i \in  (K_1 \backslash K_2) \times (K_2 \backslash K_1) $, there is
     $v_i := (v_{1i}, v_{2i})  \in T_{K_1}(\bar{x}_{1i}) \times  T_{K_2}(\bar{x}_{2i})$ and $w_i \in \mathbb{R}^n$ such that 
     \begin{equation}\label{equse} 
    - f(\bar{x}_i) = \nabla f(\bar{x}_i)^\top v_i + w_i, 
    \end{equation}
    \begin{equation}\label{equse1}
     |[v_{1i}^\top ~~ v_{2i}^\top]| \leq \bar{c} |f(\bar{x}_i) |, ~~~~ |w_i| \leq \bar{\alpha} |f(\bar{x}_i)|.   
    \end{equation}
     Indeed, the transversality property ($\diamond$) allows us to conclude, when $ \bar{x}_i \in (\partial K_1 \backslash K_2) \times (\partial K_2 \backslash K_1), $
    the existence of $v_i := (v_{1i}, v_{2i})  \in T_{K_1}(\bar{x}_{1i}) \times  T_{K_2}(\bar{x}_{2i})$ and $w_i \in \mathbb{R}^n$ such that $- f(\bar{x}_i) = \nabla f(\bar{x}_i)^\top v_i + w_i$ with $|[v_{1i}^\top ~ v_{2i}^\top]| \leq c |f(\bar{x}_i)|$ and $|w_i| \leq \alpha |f(\bar{x}_i)|.$
    Otherwise, when $\bar{x}_{1i} \in \text{int}(K_1) \backslash K_2$ (resp. $\bar{x}_{2i} \in \text{int}(K_2) \backslash K_1$), 
    since $T_{K_1}(\bar{x}_{1i}) = \mathbb{R}^n$ (resp., $T_{K_2}(\bar{x}_{2i}) = \mathbb{R}^n$), 
    we can just take $v_i := (\bar{x}_{2i} - \bar{x}_{1i}, 0)$ (resp. $v_i := (0, \bar{x}_{1i} - \bar{x}_{2i})$) so that~\eqref{equse} holds with  
    $ |[v_{1i}^\top ~~ v_{2i}^\top]| = |f(\bar{x}_i) |$ and $|w_i| = 0$. 
    Hence, whenever $ \bar{x}_i \in  (K_1 \backslash K_2) \times (K_2 \backslash K_1) $, there exists 
     $v_i := (v_{1i}, v_{2i})  \in T_{K_1}(\bar{x}_{1i}) \times  T_{K_2}(\bar{x}_{2i})$ and $w_i \in \mathbb{R}^n$ such that 
    ~\eqref{equse}-\eqref{equse1} hold with $\bar{c} := \max \{c,1\}$ and $\bar{\alpha} := \max \{0, \alpha\}$. 
    
    By definition of the contingent cone, we conclude the existence of $h_p \rightarrow 0$ and $e_p \rightarrow 0$ such that 
    $ y_i :=  \bar{x}_i + h_p ( v_i + e_p ) \in K_1 \times K_2$
    and 
    $ f(y_i) = f(\bar{x}_i) + h_p \nabla f(\bar{x}_i)^\top (v_i + e_p) + h_p o(h_p).  $
    Using~\eqref{equse}, we obtain 
    $ f(y_i) = f(\bar{x}_i) + h_p ( - f(\bar{x}_i) - w_i + \nabla f(\bar{x}_i)^\top e_p) + h_p o(h_p).  $
    Hence, 
    \begin{align*}  
    f(y_i) & = (1-h_p) f(\bar{x}_i) - h_p w_i + h_p \nabla f(\bar{x}_i)^\top  e_p + h_p o(h_p). 
    \end{align*}
    
    Now, using the second inequality in~\eqref{equse-1}, we obtain 
    $|f(\bar{x}_i)| - \varepsilon 
    |\bar{x}_i - y_i | < |f(y_i)|$.
    Hence, 
    \begin{align*}
        |f(\bar{x}_i)|  \leq & \ (1-h_p) |f(\bar{x}_i)| + h_p |w_i| + h_p |\nabla f(\bar{x}_i)^\top  e_p + o(h_p)|\\
        & + \varepsilon h_p | v_i + e_p|.
    \end{align*} 
    The latter allows us to conclude that 
    $ h_p |f(\bar{x}_i)|  \leq h_p |w_i| + h_p |\nabla f(\bar{x}_i)^\top  e_p + o(h_p)|  + \varepsilon h_p | v_i + e_p |.   $
    Dividing by $h_p$ and letting $p \rightarrow \infty$, we obtain 
    $$  |f(\bar{x}_i)|  \leq |w_i| + \varepsilon  | v_i | \leq (\bar{\alpha} + \varepsilon \bar{c}) |f(\bar{x}_i)|.   $$
    Now, if $\varepsilon > 0$ is such that $\bar{\alpha} + \varepsilon \bar{c} < 1$, we conclude that  $|f(\bar{x}_i)| = 0$ and thus $\bar{x}_{i1} = \bar{x}_{i2}$,  contradicting~\eqref{eqcontradict}.
\end{proof}

\ifitsdraft
In the following lemma, we replace ($\diamond$) by a simpler, and yet  more conservative, condition.

\fi

\begin{lemma}\label{lemtranssimple}
    Consider closed subsets $K_1$, $K_2 \subset \mathbb{R}^n$, and let $x \in \partial K_1 \cap \partial K_2$. 
    The property ($\diamond$) in Lemma~\ref{lemtrans} is verified provided that 
    \begin{itemize}
        \item [($\diamond \diamond$)] $\exists$ a neighborhood $\mathcal{N}(x)$ around $x$ such that
        \begin{equation*}
            \begin{aligned}
              \hspace{-0.7cm} & 
              x_2 - x_1 \in  
              C_{K_1} (x) - C_{K_2} (x)
              \\
               \hspace{-0.7cm} & \forall (x_1,x_2) \in 
                ((\partial K_1 \backslash K_2) \cap \mathcal{N}(x)) \times ((\partial K_2 \backslash K_1) \cap \mathcal{N}(x)).
                \end{aligned}   
        \end{equation*}
    \end{itemize}
\end{lemma}
\begin{proof}   
A key argument to prove the lemma consists of using the Open Mapping Theorem~\cite[Theorem 2.2.1]{aubinSetValuedAnalysis2009} to  conclude the existence of $l > 0$ such that, for any $(x_1,x_2) \in ((\partial K_1 \backslash K_2) \cap \mathcal{N}(x)) \times ((\partial K_2 \backslash K_1) \cap \mathcal{N}(x))$, there exists $(c_1,c_2) \in C_{K_1}(x) \times C_{K_2}(x)$ such that 
$$ (x_2 - x_1) = (c_1 - c_2) ~~ \text{and} ~~ |[c_1^\top ~ c_2^\top]| \leq l |c_2 - c_1|. $$ 

Indeed, we note that the map $F : C_{K_1}(x) \times C_{K_2}(x) \rightarrow C_{K_1}(x) - C_{K_2}(x)$ given by $F(c_1,c_2) := c_1 - c_2$
admits a closed and convex graph which is a cone. Hence, by the Open Mapping Theorem \cite[Theorem 2.2.1]{aubinSetValuedAnalysis2009}, we conclude that $F^{-1}$ is locally Lipschitz. 
Next, we let 
$$ K := ((\partial K_1 \backslash K_2) \cap \mathcal{N}(x) ) - ((\partial K_2 \backslash K_1) \cap \mathcal{N}(x)). $$
We can see, under property ($\diamond \diamond$), that $y_o := 0 \in \text{int}_{K}( C_{K_1}(x) - C_{K_2}(x))$, where $\text{int}_{K}(\cdot)$ stands for the interior of $(\cdot)$ relative to $K$. In other words, there exists $\epsilon > 0$ such that $(y_o + \epsilon \mathbb{B}) \cap K \subset C_{K_1}(x) - C_{K_2}(x)$. 
Next, we note that 
$(c_{1o}, c_{2o}) := (0,0) \in F^{-1}(y_o). $ 
The Lipschitz property of $F^{-1}$ allows us to conclude that, there exists $\epsilon > 0$ sufficiently small and $l > 0$ such that, $\forall y \in (y_o + \epsilon \mathbb{B} ) \cap K \subset C_{K_1}(x) - C_{K_2}(x)$, 
we can find $(c_1,c_2) \in C_{K_1}(x) \times C_{K_2}(x)$ such that $F((c_1,c_2)) = c_1 - c_2 = y$ and 
$ |[c_1^\top ~ c_2^\top] - [c_{1o}^\top ~ c_{2o}^\top]| \leq l |y-y_o| = l|y| = l |c_2 - c_1|$. 

Next, since the maps $x \leadsto C_{K_1}(x)$ and $x \leadsto C_{K_2}(x)$ are lower semicontinuous, we conclude that, for each $\varepsilon > 0$, there exists a neighborhood $\mathcal{N}(x)$ such that, for each $(c_1,c_2) \in C_{K_1}(x) \times C_{K_2}(x)$  and for each $(x_1,x_2) \in ((\partial K_1 \backslash K_2) \cap \mathcal{N}(x)) \times ((\partial K_2 \backslash K_1) \cap \mathcal{N}(x))$,  there exists $(v_1,v_2)  \in T_{K_1}(x_1) \times T_{K_2}(x_2)$
such that 
$ |[v_1^\top ~ v_2^\top] - [c_1^\top ~ c_2^\top]| \leq \varepsilon |[c_1^\top ~ c_2^\top]|$.
As a result, 
$$ x_2 - x_1 = (v_1 - v_2) + w, \quad w := (c_1-c_2) - (v_1 - v_2),  $$ 
where  
$|w| \leq \varepsilon l |x_2 - x_1|$   and  
$|[v_1^\top ~ v_2^\top]| \leq  (1 + \varepsilon) l |x_2 - x_1|$. 
Finally, we choose $\varepsilon$ small so that $\varepsilon l < 1$.  
\end{proof} 

\ifitsdraft

The following lemma allows us to verify Assumption~\ref{ass_propagation} and the first item in Assumption~\ref{ass_transversality} by imposing further regularities on the sets $K$ and $C$. 

\begin{lemma} \label{lemAp}
    Consider closed subsets $K$, $C \subset \mathbb{R}^n$, and let $\mathcal{P} \subset \partial K \cap \partial C$ such that the following property holds.
    \begin{itemize}
    \item[(Pr1)] The map $T_{K} : \partial K \rightrightarrows \mathbb{R}^n$ is continuous on $\mathcal{P}$, and 
        $$ x \in \text{cl}(\partial K \cap \text{int}(C)) \qquad  \forall x \in \mathcal{P}. $$
    \end{itemize}
    Additionally, we let $F \cap T_C : C \rightrightarrows \mathbb{R}^n$ 
    be continuous at every $x \in \mathcal{P}$ and~\eqref{eqconstNag} holds. Then, for all $x \in \mathcal P$, $F(x) \cap T_C(x) \subset T_K(x)$ and $T^a_{\partial K}(x) = T_{\partial K}(x)$.
\end{lemma}

\begin{proof}
  Given $x \in \mathcal{P}$, (Pr1) entails the  existence of 
$\{x_i\}^{\infty}_{i=1} \subset \partial K \cap \text{int}(C)$ such that 
$x = \lim_{i} x_i$. Moreover, based on~\eqref{eqconstNag},  we know that 
$ F(x_i) = F(x_i) \cap T_C(x_i) \subset T_{K}(x_i)$ for all  $i \in \{1,2,\ldots\}$.  
Hence, since both $T_{K}$ and $F \cap T_C$ are continuous at $x$, it follows that  $F(x) \cap T_C(x) \subset T_{K}(x)$.  

 To prove $T^a_{\partial K}(x) = T_{\partial K}(x)$, it is enough to show that the map $T_{\partial K} : \partial K \rightrightarrows \mathbb{R}^n$ is continuous at every $x \in \mathcal{P}$. Indeed, we note that the set $K$ is sleek; hence, $C_{K}(x) = T_{K}(x)$. Moreover, from~\cite[Theorem 2]{rockafellar1979clarke}, we know that $ \mathrm{int}(C_{K}(x)) \subset D_{K}(x)$. 
On the other hand, taking into account that $D_{K}(x) \subset T_{K}(x)$, and that $D_{K}(x)$ is open, it holds that $D_{K}(x) \subset \mathrm{int}(T_{K}(x))$. Hence, $D_{K}(x) = \mathrm{int}(T_{K}(x))$. Combining the latter to~\cite[Theorem 4.3.3]{Aubin1991Viability}, we obtain 
\begin{equation*}
    \begin{split}
       & \partial T_{K}(x) = T_{K}(x)\backslash \mathrm{int}(T_{K}(x)) = T_{K}(x) \backslash D_{K}(x) =
       \\ & 
     T_{K}(x) \backslash (\mathbb{R}^n \backslash T_{\mathbb{R}^n \backslash K}(x)) = T_{K}(x) \cap T_{\mathbb{R}^n \backslash K}(x) = T_{\partial K}(x).
    \end{split}
\end{equation*}
Finally,  since $T_{K} : \partial K \rightrightarrows \mathbb{R}^n$ is continuous, then $ \partial T_{K}\nolinebreak(\nolinebreak=\nolinebreak T_{\partial K})  : \partial K \rightrightarrows \mathbb{R}^n$ is also continuous. 
\end{proof}
\fi



\begin{lemma} \label{lem_prac_properties}
    Let $K \subset \mathbb{R}^n$ be a finitely-represented set fulfilling \eqref{eq_cq_pra} with $\star = K$ and let $x \in \partial K$. Then, 
    \begin{itemize}
        
\item[1)] The map $K \ni x \leadsto T_K(x)$ is lower semicontinuous at $x$ and
        $
        T_K(x) = T^a_K(x) = C_K(x).  
        $
        
        \item[2)] $ D_{K}(x)  = \left\{ v \in \mathbb{R}^n :  \nabla g_{K}^j(x)^\top v < 0 ~~ \forall j \in \mathrm{Act}_{K}(x)  \right\}$. 
        
        \item[3)] 
        $ T_{\partial K}(x)  = \partial T_{K}(x) = $
        \\
        $ \left\{ v \in T_K(x) : \exists j \in \mathrm{Act}_{K}(x) :  \nabla g_{K}^j(x)^\top v = 0 \right\}$. 
    \end{itemize}
    
    Additionally,  provided that
    \begin{itemize}
        \item [(Cd)] $\forall j \in \mathrm{Act}_{K}(x)$, there exists $v_j \in \mathbb{R}^n$ such that 
        \begin{align*} 
            \hspace{-0.4cm} \nabla g_{K}^j(x)^\top v_j  = 0,  ~  \nabla g_{K}^i(x)^\top v_j < 0 ~~ \forall i \in \mathrm{Act}_{K}(x) \backslash\{j\},
        \end{align*}
    \end{itemize}
    \begin{itemize}
        \item[4)] For each $j \in \text{Act}_K(x)$,  $\partial K_j \ni x \leadsto T_{\partial K_j}(x)$ is lower semicontinuous at $x$ and $T_{\partial K_j}(x) = C_{\partial K_j}(x)$, where
        \begin{equation}\label{eq_def_Kj}
            \partial K_j := \{x \in K~:~ g_K^j(x) = 0\}.
        \end{equation}
        \item[5)] $T_{\partial K}(x) = T^a_{\partial K}(x)$. 
    \end{itemize}
\end{lemma}

\begin{proof}
   To prove item 1),  we consider $x \in \partial K$ and a neighborhood $\mathcal{N}(x)$ around $x$. We already know that, for any $y \in \mathcal{N}(x) \cap  \partial K$,  
    \[ 
    T_{K}(y) = \{v \in \mathbb{R}^n  :  \nabla g_{K}^j(y)^\top v \leq 0 ~~ \forall j \in \mathrm{Act}_{K}(y)\}.
    \]
    Furthermore, note that for any $i \in [d_K]\backslash \mathrm{Act}_{K}(x)$, it holds that $g_K^i(x) < 0$. Hence, by continuity of $g_K^i$, there exists a neighborhood, $\mathcal{N}_i(x)$ such that
        $ g_K^i(y) < 0$ for all $y \in \mathcal{N}_i(x)$.
        By taking
       $ \mathcal{N}(x) := \bigcap_{i \in [d_K]\backslash \mathrm{Act}_{K}(x)} \mathcal{N}_i(x)$,
        it follows that
          $ \mathrm{Act}_{K}(y) \subset \mathrm{Act}_{K}(x)$ for all $y \in \mathcal{N}(x) \cap  \partial K$. 
    Combining the latter to continuity of $\nabla g_{K}^i$ for all $i \in \mathrm{Act}_{K}(x)$, we conclude that, for any sequence $\{y_i\}^{\infty}_{i=1} \subset \mathcal{N}(x) \cap  \partial K$ that converges to $x$, we have
    \begin{align*}
     \lim_{i \rightarrow \infty}  \left\{ \nabla g_{K}^{j}(y_i), j \in \mathrm{Act}_{K}(y_i) \right\} \subset 
       \left\{  \nabla g_{K}^k(x),~  k \in \mathrm{Act}_{K}(x) \right\}. 
    \end{align*}
    This plus the fact that  
    $\nabla g_{K}^k(x) \neq 0$ for all $k \in \mathrm{Act}_{K}(x)$ 
    are enough to conclude lower semicontinuity of $T_K$ at $x$.  
    To prove item 2), we use sleekness of the set $K$ to conclude that 
    $C_{K}(x) = T_{K}(x)$ for all $x \in \partial K$.    Moreover, from~\cite[Theorem 2]{rockafellar1979clarke},  we know that $ \mathrm{int}(T_{K}(x)) = \mathrm{int}(C_{K}(x)) \subset D_{K}(x)$.  On the other hand, taking into account that $D_{K}(x) \subset T_{K}(x)$ and that $D_{K}(x)$ is open, we conclude that $D_{K}(x) \subset \mathrm{int}(T_{K}(x))$. Hence, $D_{K}(x) = \mathrm{int}(T_{K}(x))$.   Finally,  since $\nabla g_{K}^j(x) \neq 0$ for all $j \in \mathrm{Act}_{K}(x)$,  we conclude that  
    \begin{align*}
    \mathrm{int}(T_{K}(x))  = \left\{ v \in \mathbb{R}^n :  \nabla g_{K}^j(x)^\top v < 0 ~~ \forall j \in \mathrm{Act}_{K}(x)  \right\}. 
    \end{align*}
    
    For 3), combining  item 2) with~\cite[Theorem 4.3.3]{Aubin1991Viability} gives 
    \begin{align*}
             \partial T_{K}(x) & = T_{K}(x)\backslash \mathrm{int}(T_{K}(x)) = T_{K}(x) \backslash D_{K}(x)  \\
            & = T_{K}(x) \backslash (\mathbb{R}^n \backslash T_{\mathbb{R}^n \backslash K}(x)) 
            \\ & = T_{K}(x) \cap T_{\mathbb{R}^n \backslash K}(x) = T_{\partial K}(x).
    \end{align*}
    
    To prove item 4), we pick $j \in \text{Act}_K(x)$ and show the existence of a neighborhood $\mathcal{N}(x)$ such that, for each $y \in \mathcal{N}(x) \cap \partial K_j$, there exists $\bar{v}_{j}(y) \in \mathbb{R}^n$ such that 
    \begin{equation}
    \label{eqvjbar}
    \begin{aligned} 
     \nabla g_{K}^j(y)^\top \bar{v}_{j}(y) & = 0,  
     \\
     \nabla g_{K}^i(y)^\top \bar{v}_{j}(y) & < 0 \qquad \forall i \in \mathrm{Act}_{K}(y) \backslash\{j\}. 
    \end{aligned}
    \end{equation}
    First, for $\mathcal{N}(x)$ sufficiently small, we use a continuity argument to conclude
    \begin{align} \label{eqinclus} 
    \text{Act}_K(y) \subset \text{Act}_K(x) \qquad \forall y \in \mathcal{N}(x) \cap \partial K. 
    \end{align}
    Note that $v_j \in \text{Ker} \left( \nabla g_{K}^j(x)^\top \right)$. Since $\nabla g_{K}^j(\cdot)^{\top}$ is continuous and $\nabla g_{K}^j(x)^\top \neq 0$, we conclude that for $\mathcal{N}(x)$ sufficiently small, there exists a continuous map 
    $y \mapsto \bar{v}_j(y) \in \text{Ker} \left( \nabla g_{K}^{j}(y)^\top \right)$ defined on $\mathcal{N}(x) \cap \partial K_j$ that converges to $v_j$ as $y$ converges to $x$. At the same time, by continuity, when $\mathcal{N}(x)$  small enough we conclude that 
    $$  \nabla g_{K}^i(y)^\top \bar{v}_{j}(y)  < 0 \quad \forall i \in \mathrm{Act}_{K}(x) \backslash\{j\}, ~ \forall y \in \mathcal{N}(x) \cap \partial K_j, $$
    which is enough to conclude~\eqref{eqvjbar} under~\eqref{eqinclus}. 
    
    After that, we use~\cite[Proposition 4.3.7]{aubinSetValuedAnalysis2009} under~\eqref{eqvjbar} to conclude that, for each $y \in \partial K_j \cap \mathcal{N}(x)$, we have 
    \begin{align*}
      T_{\partial K_j}(y) & = \left\{ v \in T_K(x) : \; \nabla g_{K}^j(y)^\top v = 0 \right\}.
    \end{align*}
    Furthermore, since the set $K$ is fulfills \eqref{eq_cq_pra} with $\star= K$,
    it follows by the first item that $T_{K}(x) = C_{K}(x) $ for all $x \in \partial K$, and $\partial K \ni y \leadsto T_K(y)$ is lower semicontinuous. Now,  using continuity of $\nabla g_{K}^j$ and since  $\nabla g_{K}^j(y) \neq 0$ for all $y \in \partial K_j \cap \mathcal{N}(x)$, then $\partial K_j \ni y \leadsto T_{\partial K_j}(y)$ is lower semicontinuous at $x$, thus, $T_{\partial K_j}(x) = C_{\partial K_j}(x)$.

    Finally, to prove  item 5), we let $v \in T_{\partial K}(x)$ for some $x \in \partial K$ and we show that $v \in T^a_{\partial K}(x)$.  According to the fourth item,  there exists $j \in \text{Act}_K(x)$ such that $\nabla g^j_K(x) v = 0$.  Let us define the vector function $w_{j}(\beta) := \beta v + (1 - \beta) v$ for all $\beta \in [0,1]$.
    
    Note that,  for all $\beta \in (0,1]$,  we have $ \nabla g_{K}^j(x)^\top w_{j}(\beta)  = 0$  and
       $ \nabla g_{K}^i(x)^\top w_{j}(\beta)  < 0$  for all $i \in \mathrm{Act}_{K}(x) \backslash\{j\}$. 
    Furthermore,  since $T^a_{\partial K}(x)$ is closed,  then $v \in T^a_{\partial K}(x)$ provided that 
    $ w_{j}(\beta) \in T^a_{\partial K}(x)$  for all $\beta \in (0,1]$.   
    Hence,  to ease the notation,  we pick $w \in T_{\partial K}(x)$ that satisfies the same properties as $w_{j}(\beta)$ when $\beta \in (0,1]$, i.e., 
$ \nabla g_{K}^j(x)^\top w = 0$, 
            $
            \nabla g_{K}^i(x)^\top w < 0$ ~ 
            for all $i \in \mathrm{Act}_{K}(x) \backslash\{j\}$, 
    and we  show that $w \in T^a_{\partial K}(x)$.  
    We do so by considering an arbitrary sequence $\{h_k  \} \subset \mathbb{R}_{>0}$ such that $h_k  \rightarrow 0$,  and we let $y_k := \text{Proj}_{\tilde{K}_j}(x + h_k w)$ be the projection of $x_k := x + h_k w$ onto the set 
    $ \tilde{K}_j :=  \{ x \in \mathbb{R}^n : g^j_K(x) = 0 \}.  $
    
    We verify the following two items, to ensure $w \in T^a_{\partial K}(x)$. 
    \begin{itemize}
        \item $y_k \in \partial K$  when $k$ is sufficiently large. 
        \item $w_k := (y_k - x)/h_k    \rightarrow w$.
    \end{itemize}
     
    To show the second item, we start showing that the sequence $\{w_k \}_k$ is uniformly bounded.  Indeed, for $k \in \mathbb{N}$,    
    $
\text{Proj}_{\tilde{K}_j}(x + h_k w)  - x  = 
        \left[ \text{Proj}_{\tilde{K}_j}(x + h_k w)  - (x + h_k w) \right] + h_k w$.
    Hence, 
    $ |\text{Proj}_{\tilde{K}_j}(x + h_k w)  - x| 
        \leq  |x + h_k w|_{\tilde{K}_j}  + h_k |w|$
    and, at the same time,  
    $|(x + h_k w)|_{\tilde{K}_j} \leq  |x + h_k w -x|  = h_k |w|$, leading to $|w_k| = |\text{Proj}_{\tilde{K}_j}(x + h_k w)  - x|/h_k  \leq 2 |w| $ for all $k \in \{1,2,...\}$.
    
    Next, since $y_k$ is the projection of $x_k$ onto $\tilde{K}_j$, we conclude, for $k$ large enough, that 
    $| \nabla g^j_K(y_k)^\top  (x_k - y_k) | = |\nabla g^j_K(y_k)|  \cdot  |x_k - y_k|.$  
    As a result,  we obtain
    \begin{align*} 
        0 = g^j_K(y_k) & =  g^j_K(x) +   \nabla g^j_K(x)^\top (y_k - x)  + o(y_k - x)  \\ & 
         = \nabla g^j_K(x)^\top (y_k - x)  + o(y_k - x).
    \end{align*}
     
    Furthermore, we express $y_k$ as 
    $ y_k := x + h_k w_k = x + h_k w + h_k (w_k - w)$, and use $\nabla g^j_K(x)^\top w = 0$ to obtain
    \begin{align*}
        0  & =   \nabla g^j_K(x)^\top (w_k - w) + \frac{o(h_k w_k)}{h_k}  
        =   \nabla g^j_K(y_k)^\top (w_k - w) \\ & 
        \;\;\; +  \left( \nabla g^j_K(x) - \nabla g^j_K(y_k)  \right)^\top (w_k - w) + o(h_k w_k)/h_k \\ & 
        \in \left( \nabla g^j_K(x) - \nabla g^j_K(y_k)  \right)^\top (w_k - w) + o(h_k w_k)/h_k \\ &
        \;\;\; + (\nabla g^j_K(y_k)|  |w_k - w|)\mathbb{B}.
    \end{align*}
    Using continuity of $ \nabla g^j_K$,  we conclude that
    $$\left[ \nabla g^j_K(x) - \nabla g^j_K(y_k) \right]  \rightarrow 0,  ~~  \nabla g^j_K(y_k) \rightarrow \nabla g^j_K(x) \neq 0.$$
    Moreover, by boundedness of the sequence $\{w_k \}_k$, we conclude that
    $(w_k - w) \rightarrow 0$ since 
    $
    o(h_k w_k)/h_k   \rightarrow 0$ and  
    $
     \left( \nabla g^j_K(x) - \nabla g^j_K(y_k)  \right)^\top (w_k - w)  \rightarrow 0$.

    Finally,  it remains to prove the first item, i.e., to show that $y_k = x + h_k w_k \in \partial K$ for large values of $k$.  Indeed,  it holds by definition that 
    $ g^j_K(y_k) = 0$ for all $k$.  Furthermore,  given $i \in \text{Act}_K(x) \backslash \{j\}$,  we have 
    \begin{align*} 
     g^i_K(y_k) & =  g^i_K(x) +  \nabla g^i_K(x)^\top (y_k - x)  + o(y_k - x) 
    \\ & =  \nabla g^i_K(x)^\top  w + \nabla g^i_K(x)^\top (w_k - w) + o(h_k w_k)/h_k. 
    \end{align*}
    By the properties of $w$,  we get $ \nabla g^i_K(x)^\top  w <  0$,  which implies that  $g^i_K(y_k) \leq  0$ when $k$ is large,  under boundedness of $\{w_k\}_k$ and since $(w_k - w) \rightarrow 0$.  Hence,  $y_k \in \partial K$ when $k$ is sufficiently large.
\end{proof}

We next verify Assumption~\ref{ass_propagation} when $(K,C)$ are finitely-represented sets, Assumptions~\ref{ass_cq_practical} and~\ref{ass_transversalityPra}, and \eqref{eqMotivprac} hold.

\begin{lemma}\label{lem_transPrac_implies_eqTangCond}
    Consider system $\Sigma$ such that $F : C \rightrightarrows \mathbb{R}^n$ 
   is additionally continuous. Let 
    finitely-represented
    sets $K, C \in \mathbb{R}^n$ such that Assumption~\ref{ass_cq_practical} holds. Then, under \eqref{eqMotivprac},  Assumption~\ref{ass_transversalityPra} implies Assumption~\ref{ass_propagation}.
\end{lemma}

\begin{proof}
    Given $x \in \mathcal{P}$ and $i \in \text{Act}_K(x)$, we show that 
    \begin{align} \label{eqthm21} 
    \hspace{-0.6cm} \partial K_i \cap \text{int}(C) \cap \mathcal{N}(x) \neq \emptyset  \qquad 
     \forall \mathcal{N}(x),   
    \end{align} 
    where $\partial K_i := \{ x \in K :  g_{K}^i(x) = 0\}$.
    
    Under~\eqref{eqthm21},  we conclude the existence of a sequence $\{ y^i_k \}^{\infty}_{k=1} \subset \partial K_i \cap \text{int}(C)$ such that $y^i_k \rightarrow x$ as $k \rightarrow \infty$.  
    Furthermore, under~\eqref{eqconstNag},  we conclude that 
    $$  \nabla g^i_K(y^i_k)  \eta \leq 0 \qquad \forall \eta \in F(y_k),  \quad \forall k \in \{1,2,...\}.   $$
    A continuity argument allows us to conclude that 
    $\nabla g^i_K(x)  \eta \leq 0$  for all $\eta \in F(x)$, for all $i \in \text{Act}_K(x)$.  
     The latter is enough to verify Assumption~\ref{ass_propagation}  according to the first bullet in Lemma~\ref{lem_prac_properties}. 
    
    Now,  to verify~\eqref{eqthm21} for a given $i \in \text{Act}_K(x)$, we start using the fourth bullet in Lemma~\ref{lem_prac_properties} 
    to conclude that
    \begin{align*} 
    T_{\partial K_i}(x) & = 
    \left\{ v \in T_K(x) :   \nabla g_{K}^i(x)^\top v = 0 \right\}.
    \end{align*}
    Furthermore, from Assumption~\ref{ass_transversalityPra}, there exists a vector $v_i \in T_{\partial K_i}(x)$
    such that $\nabla g_{K}^j(x)^\top v_i 
             < 0$ 
            for all $j \in \mathrm{Act}_{K}(x) \backslash\{i\}$ and 
     $\nabla g_{C}^k(x)^\top v_i 
             < 0$ for all $k 
            \in \mathrm{Act}_{C}(x)$. 
    
    As a result, there exists $\{ w^i_k \}_k \subset \mathbb{R}^n$ such that $w^i_k \rightarrow v_i$ as $k \rightarrow \infty$ and $\{h_k\}_k \subset \mathbb{R}_{>0}$ such  that 
    $h_k \rightarrow 0$ as $k \rightarrow \infty$ such that 
    $y_k := x + h_k  w^i_k \in \partial K_i$ for all $k$.
  Also, 
     \begin{align*} 
         g^j_C(y_k) =  h_k  \nabla g_{C}^j(x)^\top v_i  +  h_k  \nabla g_{C}^j(x)^\top  (  w^i_k -  v_i )
             + o(h_k)  
      \end{align*}
     for all $j \in \mathrm{Act}_{C}(x)$. As a result,  for $k$ sufficiently large, we conclude that  
     $ g^j_C(y_k) < 0$ for all $j \in \mathrm{Act}_{C}(x)$,
     and thus 
     $y_k \in \text{int}(C)$  for all $j \in \mathrm{Act}_{C}(x)$,
      proving~\eqref{eqthm21}.
\end{proof}

 \ifitsdraft

In the following lemma, when $(K,C)$ are finitely-represented sets fulfilling Assumption~\ref{ass_cq_practical} 
and (Cd) from Lemma~\ref{lem_prac_properties} holds with $K$ therein replaced by $C$, we establish the equivalence between Assumptions~\ref{ass_decisiveness} and~\ref{ass_decisiveness_pra}.  

\begin{lemma} \label{lem_nontang_prac}
     Consider finitely-represented 
    sets $K, C \in \mathbb{R}^n$ such that Assumption~\ref{ass_cq_practical} holds and $\mathcal{P} \subset \partial K \cap \partial C$ such that Assumption~\ref{ass_transversalityPra} holds. Assume that for each $x \in \mathcal{P}$ and $j \in \mathrm{Act}_{C}(x)$, there exists $v_j \in \mathbb{R}^n$ such that 
    \begin{align*} 
        \hspace{-0.4cm} \nabla g_{C}^j(x)^\top v_j  = 0,  \quad  \nabla g_{C}^i(x)^\top v_j < 0 ~~ \forall i \in \mathrm{Act}_{C}(x) \backslash\{j\}. 
    \end{align*}
    Then Assumption~\ref{ass_decisiveness} is equivalent to Assumption~\ref{ass_decisiveness_pra}.
\end{lemma}

\begin{proof}
The proof is based on showing that  
\begin{equation}\label{eqTpartialKC}
    \begin{aligned} 
         T_{\partial K \cap \partial C}(x)  & = T_{\partial K}(x) \cap T_{\partial C}(x) \qquad \forall x \in \mathcal P.
    \end{aligned}
\end{equation}

To verify~\eqref{eqTpartialKC}, we show that ($\dagger$)  is verified while replacing $C$ therein is by $\partial C$.  To do so, we show the following property.
    \begin{itemize}
        \item ($\ddagger$) in Section~\ref{subsec_alternatives} holds while replacing $(\partial K, C)$ therein by $(\partial K_j, \partial C_i)$ for all $(j,i) \in \text{Act}_K(x) \times  \text{Act}_C(x)$.
    \end{itemize}

Then,  applying Lemma~\ref{lemtranssimple},  we conclude the existence of a neighborhood $\mathcal{N}(x)$,  $c_{ji}> 0$,  and $\alpha_{ji} \in [0,1)$ such that
$$ \forall (x_{1j},x_{2i}) \in ((\partial K_j \backslash C) \cap \mathcal{N}(x)) \times ((\partial C_i \backslash \partial K) \cap \mathcal{N}(x)),$$
there exists $(\bar{v}_j,\bar{w}_i) \in  T_{\partial K_j} (x_{1j}) \times T_{\partial C_i} (x_{2i})$ such that 
$|[\bar{v}_j^\top~\bar{w}_i^\top]| \leq c_{ji} |x_{1j} - x_{2i}|$ and 
\begin{equation*}
    x_{2i} - x_{1j} \in \  \bar{v}_j - \bar{w}_i + \alpha_j (|x_{1j} - x_{2i}|) \mathbb{B}.
\end{equation*}
Next,  for $\mathcal{N}(x)$ sufficiently small, we conclude that 
\[
    \forall x_{1} \in (\partial K \backslash C) \cap \mathcal{N}(x), ~~ \forall x_{2} \in (\partial C \backslash K) \cap \mathcal{N}(x),
\]
there exists $j \in \text{Act}_K(x)$ and $i \in \text{Act}_C(x)$ such that 
$$ x_1 \in (\partial K_j \backslash C) \cap \mathcal{N}(x) ~~ \text{and} ~~ x_2 \in (\partial C_i \backslash K) \cap \mathcal{N}(x). $$ 
Hence, by letting 
\begin{align*}
    \alpha  := \max_{\begin{matrix} 
    j \in \text{Act}_K(x),  \\ i \in \text{Act}_C(x) \end{matrix}}
     \{\alpha_{ji} \}, ~~
    c  := \max_{ \begin{matrix} j \in \text{Act}_K(x),  \\ i \in \text{Act}_C(x) \end{matrix}} \{c_{ji} \},  
\end{align*} 
we conclude the existence of 
$$(v,w)  \in  T_{\partial K_j} (x_{1}) \times T_{\partial C_i} (x_2)   \subset T_{\partial K} (x_{1}) \times T_{\partial C} (x_2) $$ 
such that $|[v^\top ~ w^\top]| \leq c |x_{1} - x_2|$ and 
\begin{equation*}
    \begin{split}
        x_2 - x_{1} \in \ & v - w + \alpha (|x_{1} - x_2|) \mathbb{B}, 
    \end{split}
\end{equation*}
which shows ($\dagger$). 

Let us now verify the property in the latter bullet.  Given $x \in \mathcal{P}$ and $j \in \text{Act}_K(x)$,  by considering the vector $v_j$ in Assumption~\ref{ass_transversalityPra} and a sufficiently small neighborhood $\mathcal{N}(x)$  around $x$, we have

\begin{equation}\label{equseclim}
    \begin{aligned}
        \hspace{-0.4cm} \nabla g_{C}^k(x)^\top v_j - \nabla g_{C}^k(x) (x_2 - x_{1j}) \leq 0 \quad  \forall k\in\mathrm{Act}_{C}(x), \\
        \hspace{-0.4cm} \forall (x_{1j},x_2) \in ((\partial K_j \backslash C) \cap \mathcal{N}(x)) \times ((\partial C \backslash \partial K) \cap \mathcal{N}(x)). 
    \end{aligned}
\end{equation}
Furthermore, by  item 5) in Lemma~\ref{lem_prac_properties}, 
we obtain that $v_j \in T_{\partial K_j}(x) = C_{\partial K_j}(x)$ and that
\[ w_j := v_j - (x_2 - x_{1j}) \in  T_C(x) = C_{C}(x). \]

Next, we introduce the vector 
\[ w_{j}(\beta) := \beta v_j - (x_2 - x_{1j})  \qquad \forall \beta \in [0,1], \] 
and we show the existence of  $\beta_j  \in [0,1]$ such that 
\begin{align*} 
w_{j} (\beta_j) & \in T_{\partial C}(x) = \partial T_C(x) 
\\ & = \left\{ v \in T_C(x) : \exists j \in \mathrm{Act}_{C}(x) :  \nabla g_{C}^j(x)^\top v = 0 \right\}. 
\end{align*}
 
\begin{claim}\label{clam44}
    There exists $\mathcal{N}(x)$ sufficiently small such that, 
    \[ \forall (x_{1j},x_2) \in ((\partial K_j \backslash C) \cap \mathcal{N}(x)) \times ((\partial C \backslash \partial K) \cap \mathcal{N}(x)),  \] 
    there exist $\beta_j  \in [0,1]$ and $l_j \in \text{Act}_C(x)$ such that 
    $w_{j}(\beta_j) \in T_C(x)$ and 
    $\nabla g_C^{l_j}(x)^\top w_{j}(\beta_j) = 0$.
\end{claim}
    
Under the latter claim,  since 
\begin{align*} 
    T_{\partial C}(x)  = \bigcup_{j=1}^{d_{C}} T_{\partial C_i}(x) = \bigcup_{i \in \text{Act}_C(x)} T_{\partial C_i}(x),
\end{align*}
we conclude the existence of $i \in \text{Act}_C(x)$ such that $w_{j} (\beta_j) \in T_{\partial C_i}(x)$. Furthermore, under~\eqref{eq_nontangenciality_pra} we can apply the fifth bullet of Lemma~\ref{lem_prac_properties} 
 (while replacing $K$ therein by $C$) to obtain 
\begin{align*} 
T_{\partial C_i}(x) & = 
\left\{ v \in T_C(x) :   \nabla g_{C}^i(x)^\top v = 0 \right\} = C_{\partial C_i}(x). 
\end{align*}
As a consequence, 
\[ \forall (x_{1j},x_{2i}) \in ((\partial K_j \backslash C) \cap \mathcal{N}(x)) \times ((\partial C_i \backslash \partial K) \cap \mathcal{N}(x)), \]
we are able to find
$(w^*_i,v^*_j) := (w_{j} (\beta_j),  \beta_j v_j) \in C_{\partial C_i}(x) \times C_{\partial K_j}(x)$ such that 
$x_{2i} - x_{1j} = v^*_j - w^*_i$.
\\
\underline{Proof of Claim~\ref{clam44}:}
To find a contradiction, we assume that,  for each $\mathcal{N}(x)$ arbitrary small,  we can find 
\[ (x_{1j},x_2) \in ((\partial K_j \backslash C) \cap \mathcal{N}(x)) \times ((\partial C \backslash \partial K) \cap \mathcal{N}(x)) \]
such that 
\begin{equation} 
\label{equseII} 
\nabla g_C^{i}(x)^\top w_{j}(\beta) \neq  0 \qquad \forall \beta \in [0,1],  \quad \forall i \in \text{Act}_C(x).  
\end{equation}
In particular,  we have 
$$ \nabla g_C^{i}(x)^\top w_{j}(0) \neq  0 \qquad   \forall i \in \text{Act}_C(x).  $$
At the same time,  using the fact that $x_{1j} \notin C$,  we conclude the existence of $k \in \text{Act}_C(x)$  such that 
\begin{align*} 
& g^k_C(x_{1j})  = \nabla g^k_C(x_2)^\top (x_{1j} - x_2)   +  O(x_{1j} - x_2)^\top  (x_{1j} - x_2)  
\\ &  =  \nabla g^k_C(x)^\top (x_{1j} - x_2) 
\\ &  +  \left[ \nabla g^k_C(x_2) - \nabla  g^k_C(x) +  O(x_{1j} - x_2)  \right]^\top (x_{1j} - x_2)   > 0,
\end{align*}
where the vector function $O(\cdot)$ satisfies $ \lim_{\cdot \rightarrow 0} O(\cdot) = 0$.

Now,  since 
$$\nabla g^k_C(x)^\top (x_{1j} - x_2)  =  \nabla g_C^{k}(x)^\top w_{j}(0) \neq  0,  $$
and 
$$   \left[ \nabla g^k_C(x_2) - \nabla  g^k_C(x) + O(x_{1j} - x_2) \right]^\top   $$
can be made arbitrarily small when the neighborhood $\mathcal{N}(x)$ is so,  we conclude that 
$$\nabla g^k_C(x)^\top (x_{1j} - x_2)  =  \nabla g_C^{k}(x)^\top w_{j}(0) >  0.  $$
Hence,  we just showed that 
$w_{j} (\beta = 0) \notin  T_C(x)$. 

On the other hand, we know from~\eqref{equseclim} that
$$   w_{j} (\beta = 1) = v_j - (x_2 - x_{1j})   \in T_C(x).   $$
As a result,  using closeness of $T_C(x)$ and continuity of $\beta \mapsto  w_{j} (\beta)$,  we conclude the existence of 
$\beta_j \in (0,1]$ such that 
\begin{align*} 
w_{j} (\beta_j)& \in \partial T_C(x) = T_{\partial C}(x) 
\\ &
= \left\{ v \in T_C(x) : \exists i \in \mathrm{Act}_{C}(x) :  \nabla g_{C}^i(x)^\top v = 0 \right\}.   
\end{align*}  
 This contradicts~\eqref{equseII}. 

\end{proof}

\fi



\begin{lemma}\label{lem_trans_poly}
    Let $K$ and $C$ be a polytopic sets. Then, for all $x \in \partial K \cap \partial C$,
        $T_{\partial K}(x) = T^a_{\partial K}(x)$,  $T_{K}(x)  = T^a_{K}(x)$,
  and $(\dagger)$ in Assumption~\ref{ass_transversality} holds.
\end{lemma}
\begin{proof}   
     We will denote the $j$-th row of $A_K$ as $a_K^j \in \mathbb{R}^{1\times n}$, and the corresponding $j$-th element of $b_K$ as $b_K^j\in \mathbb{R}$. 
    \ifitsdraft
      First, using Lemma~\ref{lemma_Tk_TparK_poly} hereafter, it is known that
    \else
       We start using~\cite[Example 5.2.6-b)]{hiriart_urrutyFundamentalsConvexAnalysis2004} to obtain
    \fi
    \begin{align}
        T_{K}(x) & = \{v \in \mathbb{R}^n : a_K^j v \leq 0 \quad \forall j \in \mathrm{Act}_{K}(x)\},
    \label{eq_Tstar_proof_trans_poly}
 \\
        T_{\partial K}(x) & = \{v \in T_{K}(x) : \exists j \in \mathrm{Act}_{K}(x): a_K^j v = 0 \}.
\label{eq_TpartialK_proof_trans_poly}
    \end{align}
     Furthermore, we note that, for any $j \in [m_k]\backslash \mathrm{Act}_{K}(x)$, it holds that $a_K^j x < b_K^j$. Hence, there exists a neighborhood  $\mathcal{N}_j(x)$ such that
     $a_K^j y < b_K^j$ for all $y \in \mathcal{N}_j(x)$.
     By taking
        $\mathcal{N}(x) := \bigcap_{j \in [m_k]\backslash \mathrm{Act}_{K}(x)} \mathcal{N}_j(x)$,
        it follows that
        \begin{equation} \label{eq_acty_sub_actx}
            \mathrm{Act}_{K}(y) \subset \mathrm{Act}_{K}(x) \qquad \forall y \in \mathcal{N}(x).
        \end{equation}
    This combined with~\eqref{eq_Tstar_proof_trans_poly} can be used to obtain that
     $T_{K}(x) \subset T_{K}(y)$ for all $y \in \mathcal{N}(x) \cap K$.
    This proves lower semicontinuity of $x \leadsto T_{K}(x)$ on $K$.

   To prove that $T_{\partial K}(x) = T_{\partial K}^a(x)$, we let $v \in T_{\partial K}(x)$ and show that $v \in T_{\partial K}^a(x)$. For this, we choose a continuous function $\beta: \mathbb{R}_{> 0} \to [1,+\infty)$ such that $\lim_{t \to 0^+} \beta(t) = 1$ and the ratio $t/\beta(t)$ is small enough to verify $x + v t/\beta(t)  \in \mathcal{N}(x)$ for all $t \geq 0$.
    Since $v \in T_K(x)$, using~\eqref{eq_Tstar_proof_trans_poly} and~\eqref{eq_acty_sub_actx}, we conclude that 
    $
    a_{K}^j \left(x + t v/\beta(t) \right) \leq b_K^j$ for all $j \in [m_K]$, for all $t \geq 0$.
 Hence, $x + v t/\beta(t) \in K$ for all $t \geq 0$. 
Furthermore, since $v \in T_{\partial K}(x)$, then there exists $L_x \subset \text{Act}_{K}(x)$ such that
 $a_{K}^l \left(x +  t v / \beta(t) \right) = b_K^l$  for all $l \in L_x$, for all $t \geq 0$.
Hence, $x + vt/\beta(t)  \in \partial K$ for all $t \geq 0$.
    Keeping this in mind, for any sequence $\{t_i\}_{i =0}^\infty$ that converges to $0$, the sequence $\{v_i\}_{i =0}^\infty$, defined as $v_i := v/\beta(t_i)$, converges to $v$ and satisfies $x + v_i t_i \in \partial K$ for all $i \in \{0, 1, \ldots\}$.
    This is enough to conclude that $v \in T_{\partial K}^a(x)$.
    Hence, the first item holds.

    To prove $(\dagger)$, we pick a point $x \in \partial K \cap \partial C$ and a small enough neighborhood $\mathcal{N}(x)$ such that~\eqref{eq_acty_sub_actx} holds. Hence,
    \begin{equation}
    \label{eq_sist_aC_aK_poly}
        \begin{cases}
            a_K^j x = b_{K}^j & \quad \forall j \in \mathrm{Act}_{K}(x) \neq \emptyset, \\
            a_C^l x = b_{C}^l & \quad \forall l \in \mathrm{Act}_{C}(x) \neq \emptyset. 
        \end{cases} 
    \end{equation}
      
    Given $x_1 \in (\partial K \backslash C) \cap \mathcal{N}(x)$ and $x_2 \in (\partial C \backslash \partial K) \cap \mathcal{N}(x)$, we claim the existence of $v_2 \in \mathbb{R}^n$ such that
    \begin{equation} \label{eq_syst_v2_poly}
        \begin{cases}
            a_K^j v_2 = -  a_K^j(x_2 - x_1) & \forall j \in  \mathrm{Act}_{K}(x_1), \\
            a_C^l v_2 = 0 & \forall l \in  \mathrm{Act}_{C}(x_2).
        \end{cases}
    \end{equation}
    Note that $v_2 \in T_{C}(x_2)$ in this case and, by taking
    $v_1 := (x_2 - x_1) + v_2$,
    we conclude that $a_K^j v_1 = 0$ for all $j \in \mathrm{Act}_{K}(x_1)$; hence, $v_1 \in T_{\partial K}(x_1) \subset T_{K}(x_1)$. 
    Also, since 
    $(x_2 - x_1) = v_1 - v_2$.
    We conclude that $(\dagger)$ holds.
    
    We now prove the existence of $v_2$ solution to~\eqref{eq_syst_v2_poly}. Indeed, let assume the existence of the scalar sequences $ \left( \{ \lambda_j\}_{j \in \text{Act}_K(x_1)}, \{\sigma_l \}_{l \in \text{Act}_C(x_2)} \right) $
    such that 
    \begin{align} \label{eqneed} 
    \sum_{j \in \text{Act}_K(x_1)} \lambda_j a_K^j - \sum_{l \in \text{Act}_C(x_2)} \sigma_l a_C^l 
    = 0.  
    \end{align}
    ($v_2$ would exist trivially otherwise). 
    Hence, we show that 
    \begin{align} \label{eqneed2}
    \sum_{j \in \text{Act}_K(x_1)} \lambda_j a_K^j (x_2 - x_1)  = 0.  
    \end{align}
    Indeed, from~\eqref{eq_sist_aC_aK_poly}, we note that 
    \begin{align*}
    \sum_{j \in \text{Act}_K(x_1)}  \hspace{-0.4cm} \lambda_j a_K^j x_1 & = \hspace{-0.4cm} \sum_{j \in \text{Act}_K(x_1)} \hspace{-0.4cm} \lambda_j b_K^j 
   = \hspace{-0.4cm} \sum_{l \in \text{Act}_C(x_2)} \hspace{-0.4cm} \sigma_l b_C^l = \hspace{-0.4cm} \sum_{l \in \text{Act}_C(x_2)} \hspace{-0.4cm} \sigma_l a_C^l x_2. 
    \end{align*} 
    Finally, using~\eqref{eqneed}, we obtain
    $ \sum_{j \in \text{Act}_K(x_1)} \lambda_j a_K^j x_1 = \sum_{j \in \text{Act}_K(x_1)} \lambda_j a_K^j x_2.  $
    Hence,~\eqref{eqneed2} is verified. 
\end{proof}

\ifitsdraft

The following lemma recalls a characterization of the contingent cone to polytopic sets; see~\cite[Example 5.2.6-b)]{hiriart_urrutyFundamentalsConvexAnalysis2004}.

\begin{lemma}\label{lemma_Tk_TparK_poly}
    Let $K \subset \mathbb{R}^n$ be defined as
    \begin{equation}
          K = \{x \in \mathbb{R}^n  :  A x \leq b\}, ~~ A \in \mathbb{R}^{m \times n}, ~ b \in \mathbb{R}^{m}.
    \end{equation}
    Then, 
    \begin{align}
    T_{K}(x) & = \{v \in \mathbb{R}^n : a^{j} v \leq 0 \quad \forall j \in \mathrm{Act}_{K}(x)\},
    \label{eq_Tk_polytopic} 
    \\
    T_{\partial K}(x) & = \{v \in T_{K}(x) : \exists j \in \mathrm{Act}_{K}(x): a^{j} v = 0 \}, \quad\label{eq_TparK_polytopic}
    \end{align}
    where
    \[ \mathrm{Act}_{K}(x) := \begin{cases}
            \{j \in [m] : a^j x = b^j\} & \text{if } x \in K \\
            \emptyset & \text{if } x \notin K,
        \end{cases} \]
     and $a^j \in \mathbb{R}^{1\times n}$ and $b^j\in \mathbb{R}$ denote the $j$-th row of $A$ and $b$, respectively. 
\end{lemma}

\begin{proof}
    We will only prove~\eqref{eq_Tk_polytopic}. The proof of~\eqref{eq_TparK_polytopic} follows similar arguments, taking into account that
    \[
        \partial K = \{x \in K: \mathrm{Act}_{K}(x) \neq \emptyset\}.
    \]   
    To show~\eqref{eq_Tk_polytopic}, we will prove the two inclusions:
    \begin{align}
        &T_{K}(x) \subset \{v \in \mathbb{R}^n : a^{j} v \leq 0 \quad \forall j \in \mathrm{Act}_{K}(x)\},\quad\label{eq_subset_poly}\\
        &T_{K}(x) \supset \{v \in \mathbb{R}^n : a^{j} v \leq 0 \quad \forall j \in \mathrm{Act}_{K}(x)\}.\label{eq_supset_poly}
    \end{align}
    To prove~\eqref{eq_subset_poly}, we pick $w \in T_{K}(x)$ and show that
    \begin{equation}
        w \in \{v \in \mathbb{R}^n : a^{j} v \leq 0 \quad \forall j \in \mathrm{Act}_{K}(x)\}.
    \end{equation}
    To do so, we  use~\cite[Definition 1.1.3]{Aubin1991Viability} to obtain that
    \begin{equation} \label{eq_def_tgnt_cone_lim}
    \liminf_{h \to 0^+} 
    \frac{d_K(x + h w)}{h} = 0,
    \end{equation}
    where $d_K(x + h w) = \inf_{z \in K} |x + h w - z|$ is the distance from $x + h w$ to the set $K$. From~\eqref{eq_def_tgnt_cone_lim} and the polytopic nature of the set $K$, we conclude that there exists $h > 0$ such that 
    $d_K(x + h w) = 0$, i.e., $x + h w \in K$. Hence,
    \[ A(x + h w) \leq b, \]
    which implies that
    \[
     a^{j} w  \leq 0 \qquad \forall j \in \mathrm{Act}_{K}(x).
    \]
    This proves~\eqref{eq_subset_poly}.
    
    To prove~\eqref{eq_supset_poly}, we pick 
    \begin{equation} \label{eq_v_st_akv_leq_0}
        w \in \{v \in \mathbb{R}^n : a^{j} v \leq 0 \quad \forall j \in \mathrm{Act}_{K}(x)\},
    \end{equation}
    and show that $w \in T_{K}(x)$. Since $x\in \partial K$, we have that \begin{equation}\label{eq_x_in_par_K_poli}
        \begin{split}
            &a^j x - b^j = 0 \quad \forall j \in \mathrm{Act}_{K}(x),\\
            &A^j x - b^j < 0 \quad \forall j \in [m_K]\backslash\mathrm{Act}_{K}(x).
        \end{split}
    \end{equation}
    Expressions~\eqref{eq_v_st_akv_leq_0}-\eqref{eq_x_in_par_K_poli} can be combined to conclude that, for every $h>0$ sufficiently small, we have 
    \[
        \begin{split}
            & a^j (x + hv) - b^j \leq 0 \qquad \forall j \in \mathrm{Act}_{K}(x), 
            \\
            & a^j (x + hv) - b^j < 0 \qquad \forall j \in [m_K]\backslash\mathrm{Act}_{K}(x).
        \end{split}
    \]
    Hence, $(x + hv) \in K$ for small enough $h > 0$, and thus $v \in T_{K}(x)$.
\end{proof}

\fi

\ifitsdraft

The following lemma is a slight generalization of Theorem~4.6 in \cite[page 165]{bertsimasIntroductionLinear1997}. 
	
	\begin{lemma}
		\label{lem:farkas-mixed}
		Given $m_1$, $m_2$, and $n \in \mathbb{N}$, we let $E\in\mathbb{R}^{m_1\times n}$, $G\in\mathbb{R}^{m_2\times n}$, $f\in\mathbb{R}^{m_1}$,
		$h\in\mathbb{R}^{m_2}$. One of the following statements must hold:
		\begin{itemize}
			\item[(a)] There exists $x\in\mathbb{R}^n$ such that $Ex=f$ and  $Gx\le h$.
			\item[(b)] There exist $y\in\mathbb{R}^{m_1}$ and $z\in\mathbb{R}^{m_2}_{\geq 0}$ such that 
            $$ E^\top y+G^\top z=0 \quad  \text{and} \quad f^\top y+h^\top z<0. $$
		\end{itemize}
	\end{lemma}
	
	\begin{proof}
		Write $x=x^+-x^-$ with $x^+,x^-\geq 0$, and introduce a slack $ s\geq 0$ so that
		$Gx \le h \iff Gx+s=h$. Then (a) holds if and only if there exists
		$w:=(x^+,x^-,s) \geq 0$ with $Aw=b$, where
		\begin{equation}
			A := \begin{bmatrix} E & -E & 0 \\ G & -G & I \end{bmatrix}, \qquad
			b := \begin{bmatrix} f \\ h \end{bmatrix},
		\end{equation}
		which is alternative (a) of Theorem~4.6 in \cite[page 165]{bertsimasIntroductionLinear1997} applied to $(A,b)$. By Theorem~4.6 in \cite[page 165]{bertsimasIntroductionLinear1997}, this fails if and
		only if there exists $p=(y,z)$, with $y\in\mathbb{R}^{m_1}$, $z\in\mathbb{R}^{m_2}$,
		such that $p^\top A\ge 0^\top$ and $p^\top b<0$. Writing out $p^\top A \ge 0^\top$ blockwise: the
		$x^+$-columns and $x^-$-columns give $y^\top E+z^\top G\ge0$ and
		$-(y^\top E+z^\top G)\ge0$ simultaneously, forcing the equality
		$y^\top E+z^\top G=0$; the $s$-columns give $z \geq 0$, which is alternative (b). Since Theorem~4.6 in \cite[page 165]{bertsimasIntroductionLinear1997} guarantees exactly
		one of its two alternatives holds, exactly one of (a)/(b) in our statement holds.
	\end{proof}

 \fi

\ifitsdraft

\else
\fi

\balance

\end{document}